\documentclass[english]{jnsao}

\usepackage[ruled,vlined,linesnumbered,algo2e,algosection]{algorithm2e}

\usepackage{cleveref}

\usepackage{booktabs}

\usepackage{graphicx}
\usepackage{subcaption}
\usepackage{xcolor}

\usepackage{enumitem}
\setlist[enumerate,1]{label=\textup{(\roman*)}}

\renewcommand\d{\mathop{}\!\mathrm{d}}

\newcommand\TV{\text{TV}}
\newcommand\BV{\text{BV}}

\newcommand\Dmax{\underline{\Delta}}
\newcommand\pred{\operatorname{pred}}
\newcommand\ared{\operatorname{ared}}

\newcommand\du{\left\langle}
\newcommand\al{\right\rangle}

\newcommand\swpz{\bar M_0}
\newcommand\swp{\bar M}
\newcommand\wprx{w_{\text{prox}}}

\DeclareRobustCommand{\rchi}{{\mathpalette\irchi\relax}}
\newcommand{\irchi}[2]{\raisebox{\depth}{$#1\chi$}} % inner command, used by \rchi
\newcommand{\trrprox}{\rho}

\newcommand{\BB}{\mathcal{B}}
\newcommand{\conjugate}{^\star}
\newcommand{\dom}{\operatorname{dom}}
\newcommand{\dt}{\d t}
\newcommand{\epi}{\operatorname{epi}}
\newcommand{\Id}{\operatorname{Id}}
\newcommand{\N}{\mathbb{N}}
\newcommand{\prox}{\operatorname{prox}}
\newcommand{\R}{\mathbb{R}}
\newcommand{\sgn}{\operatorname{sgn}}
\newcommand{\Uad}{U_{\mathrm{ad}}}
\newcommand{\weakly}{\rightharpoonup}
\newcommand{\weaklystar}{\stackrel\star\rightharpoonup}

\DeclareMathOperator*{\argmin}{arg\,min}

\DeclarePairedDelimiter\abs{\lvert}{\rvert}
\DeclarePairedDelimiter\norm{\lVert}{\rVert}

\DeclarePairedDelimiter\parens()
\DeclarePairedDelimiter\bracks[]

\DeclarePairedDelimiterX\dual[2]{\langle}{\rangle}{#1,#2}

\DeclarePairedDelimiter\ceil{\lceil}{\rceil}
\DeclarePairedDelimiter\floor{\lfloor}{\rfloor}

\providecommand\given{\nonscript\;\delimsize|\nonscript\;\mathopen{}}
\DeclarePairedDelimiterX\set[1]\{\}{#1}

\newenvironment{msc}%
{\par\noindent{\sectfont MSC:}\ }{\par}

\newcommand{\mscLink}[1]{\href{http://www.ams.org/mathscinet/msc/msc2010.html?t=#1}{#1}}

\newenvironment{keywords}%
{\par\noindent{\sectfont Keywords:}\ }{\par}

\manuscriptsubmitted{2025-08-15}
\manuscriptaccepted{2026-08-17}
\manuscriptvolume{6}
\manuscriptnumber{16354}
\manuscriptyear{2026}
\manuscriptdoi{10.46298/jnsao-2026-16354}

\manuscriptlicense{CC-BY-NC-ND 4.0}

\manuscripteprinttype{arXiv}
\manuscripteprint{2508.10692}

\begin{document}
\title{A trust-region method for optimal control of ODEs with continuous-or-off controls and TV regularization}

\shorttitle{Control of ODEs with continuous-or-off controls}

\author{%
	Markus Friedemann%
	\thanks{%
		Brandenburgische Technische Universität Cottbus--Senftenberg,
		Institute of Mathematics,
		03046 Cottbus,
		Germany,
		\email{markus.friedemann@b-tu.de},
		\url{https://www.b-tu.de/fg-optimale-steuerung}%
	},
	Gerd Wachsmuth%
	\thanks{%
		Brandenburgische Technische Universität Cottbus--Senftenberg,
		Institute of Mathematics,
		03046 Cottbus,
		Germany,
		\email{gerd.wachsmuth@b-tu.de},
		\url{https://www.b-tu.de/fg-optimale-steuerung},
		ORCID: \href{https://orcid.org/0000-0002-3098-1503}{0000-0002-3098-1503}%
	}
}

\shortauthor{Friedemann, Wachsmuth}

\acknowledgements{%
	This research has been funded by the Deutsche Forschungsgemeinschaft (DFG, German research foundation) via grant number 458468407.%
}

\maketitle

\begin{abstract}
	A solution algorithm for a special class of optimal control problems subject
	to an ordinary differential equation is proposed. The controls
	possess a continuous-or-off structure and
	are priced by a convex function.
	Additionally, a total variation regularization is applied to penalize switches.
	Our solution method combines
	a trust-region method
	and a proximal gradient method.
	The subproblems are solved via Bellman's optimality principle.
	Convergence with respect to a criticality measure is proven.
	As a numerical example, we solve a simple optimal control problem involving an SIR model.
\end{abstract}


\begin{msc}
	\mscLink{90C11}, \mscLink{49M37}
\end{msc}

\section{Introduction}
We consider the problem
\begin{equation}
	\label{eq:prob}
	\min_{u\in \Uad} J(u) \coloneqq F(u) + \int_{t_0}^{T}   g(u(t))\dt + \TV(\sgn(u))
	.
\end{equation}
The objective function of \eqref{eq:prob} contains
a differentiable part  $F$ that encodes the control-state
system.
The control $u \colon [t_0,T] \to \set{0} \cup [a, b]$
at each time instance $t$
is either off ($u(t) = 0$) or takes values in the interval $[a,b]$ with $0 < a < b$.
The total variation term $\TV(\sgn(u))$ penalizes each switch
between these two modes.
Finally, the function $g$ acts as a pricing function
and is assumed to be strongly convex on $[a,b]$.
In \cref{sec:full_problem},
we study a variant of \eqref{eq:prob} with vector-valued controls.
The precise assumptions on the data are given in \cref{sec:full_problem}.

Our investigation of \eqref{eq:prob} was motivated by the study of control
problems for pandemic models. Tracing back to \cite{Kermack1927}, the underlying processes can be approximated by
nonlinear ordinary differential equations, see also \cite{Brauer2019}. The
$\TV$ term is introduced as a penalization for the number of switches.
Problem \eqref{eq:prob} has further applications. One of them is the optimal control of
gas networks \cite{Goettlich2021,Goettlich2018,Bachmann2019,Schmidt2020}.
Another application is
the control of heatpumps, see \cite{Huang2018}, where effective cycle planning
for multiple units can reduce peak electricity demands while ensuring sufficient heating.
It can also be applied to the cycle optimization of a single unit to
minimize energy consumption while increasing the durability.

A problem very similar to \eqref{eq:prob}
was studied in \cite{Bock2025},
although with completely different techniques. 
Up to this contribution,
we are not aware of publications related to \eqref{eq:prob}.

We propose to solve \eqref{eq:prob}
by a combination of a proximal gradient method
with a trust-region method.
These algorithms are well established
for solving nonlinear optimization problems,
see
\cite{Beck2017,Conn2000}.
In
\cite{Hahn2022,Leyffer2022,Manns2023,Manns2023b,Sharma2020}
it was shown
that the trust-region method can be applied to solve
control problems with a switching structure
which, in its simplest form, corresponds to problem \eqref{eq:prob}
with $a = b = 1$, i.e., the control can only take the values $0$ and $1$,
and each switch between these values is penalized.
An abstract convergence analysis of the trust-region method was done by \cite{Manns2024}. This abstract
trust-region framework allows for streamlined proofs of classic convergence
results.

This strongly hints at the fact, that a trust-region method will be well suited
for \eqref{eq:prob}, however for any trust-region method an efficient solution
method for the arising subproblems is required.
For a problem with switching structure,
this can be done via
combinatorial algorithms, see \cite{Marko2023,Severitt2022}.
In \cite{Marko2023}, it was shown that Bellman's optimality principle applies for the solution of
subproblems in case of one-dimensional controls and it was extended in \cite{Marko2024} to vector-valued controls.

Mathematically, this work can be viewed as an application of \cite{Manns2024} by adapting the
proximal gradient algorithm to this framework and combining it with the switching point
optimization.
We will show that the trust-region subproblems can be efficiently solved
by an adaptation of the mentioned combinatorial algorithms. In
this process, we will provide a minor generalization of \cite{Manns2024}.
In particular, we show that the algorithm already applies without a metric
structure.
The convergence proof follows the ideas of \cite{Manns2024}.

For a detailed overview of proximal gradient methods,
we refer to the book \cite{Beck2017}
and the references therein.
Recent contributions in particular concerning trust-region globalizations
in the context of nonsmooth problems
are given in, e.g., \cite{Baraldi2022,AravkinBaraldiOrban2022}.

In \cref{sec:ATRM} the abstract trust-region method (\cref{alg:ATRM}) is introduced, together with its
necessary framework. As it will be shown, the analysis for the actual solution
algorithm \cref{alg:TRM_outer_solver_2} can be split into the part relating to
a proximal gradient step and a switching point optimization. To this end, it is
beneficial to prove the convergence of the proximal gradient method interpreted
as a trust-region method using \cref{alg:ATRM} in \cref{sec:convex_problem}.
First, we prove convergence of a criticality measure.
In the case
that $F$ has a compact Fréchet derivative, strong convergence of the control is proven.
The main \cref{sec:full_problem} is dedicated to the analysis of \eqref{eq:prob}.
We start by fixing the setting and showing the existence of solutions
in \cref{subsec:setting_assumptions}.
In \cref{subsec:stationarity}, we derive optimality conditions,
see \cref{thm:OC_stationarity}.
These stationarity conditions are used to define a criticality measure in
\cref{def:OC_crit_one_switch},
which involves a proximal gradient operator on the set $\set{\sgn(u)=1}$
and terms corresponding to the switching points.
In \cref{subsec:TR_application},
we apply the trust-region approach from \cref{sec:ATRM}
to \eqref{eq:prob}.
To this end, the oracle function $T$
(which yields the suggestion for the next iterate)
is given in \cref{def:OC_TRM}.
Afterwards, we verify \cref{asmp:ATRM}.
Interestingly, all proofs basically split into
a part which addresses the continuous contributions on $\set{\sgn(u) = 1}$
(where we can reuse results from \cref{sec:convex_problem})
and a part which addresses the switching nature of the control
(which can be handled by adapting the ideas from \cite{Manns2024}).
The final convergence result is \cref{thm:OC_iterate_convergence}.
\cref{sec:full_problem} is completed by the presentation of
the inner solver in \cref{subsec:solution_subproblem}, combining a proximal gradient step with
a switching-point solver based on Bellman's optimality principle.
We present some numerical results in \cref{section:numerics}.

\section{Abstract Trust-Region Method}
\label{sec:ATRM}
This section presents the result from \cite[Section 4]{Manns2024}, such that it serves
as an introduction to the abstract trust-region method with slight notational
changes.
While the problem is stated in the setting
of metric spaces
in \cite[Section 4]{Manns2024},
the convergence proof
does not rely on this structure.
Consequently, we state the method using an abstract set $X$.

% \subsection{Setting and Assumptions}
Let $X$ be a set and let $J:X\rightarrow \R$ denote a function that is bounded from
below.
We want to solve the abstract optimization problem
\begin{equation*}
	\label{eq:AOP}
	\tag{OP}
	\min_{x\in X}J(x)
	.
\end{equation*}
The basis of our method is a function
\begin{equation*}
	\tag{oracle}
	\label{eq:oracle}
	T:X\times(0,\infty)\rightarrow X
	,
\end{equation*}
which is called ``oracle'', and a ``prediction function''
(or ``predicted reduction'')
\begin{equation*}
	\tag{$\pred$}
	\label{eq:pred}
	\pred:X\times (0,\infty) \rightarrow [0,\infty).
\end{equation*}
Given an iterate $x \in X$ and a trust-region radius $\Delta > 0$,
$T(x, \Delta)$ is a suggestion for the next iterate and
$\pred(x, \Delta)$ is a prediction for the decrease in the objective $J$.
Typically, one has a model function $m_{x, \Delta} \colon X \to (-\infty, \infty]$
and $T(x, \Delta)$ is an approximate minimizer or Cauchy point of the model
while $\pred(x, \Delta) = m_{x, \Delta}(x) - m_{x, \Delta}(T(x, \Delta))$,
see \cref{sec:convex_problem,sec:full_problem}.
Note that a model function is not used in the abstract convergence proof.
Similarly, it is not enforced that $T(x, \Delta)$ stays within distance $\Delta$ to the point $x$
since $X$ is only an abstract set and might not carry a metric.

We define the ``actual reduction'' via
\begin{equation}
	\tag{$\ared$}
	\label{eq:ared}
	\ared:X\times(0,\infty) \to \R, \qquad \ared(x,\Delta) \coloneqq J(x) - J(T(x,\Delta)).
\end{equation}
Further, we are given a
``criticality measure''
\begin{equation}
	\tag{C}
	\label{eq:C}
	C:X \to [0,\infty).
\end{equation}
The criticality measure $C$ represents the encoding of an abstract
stationarity property, where $C(x)=0$ if and only if said property holds true.

The following set of assumptions connects the
above functions, such that the iterates of a trust-region method satisfy
the convergence $\lim_{n\rightarrow\infty}C(x_n) = 0$,
see \cref{thm:atrm_convergence} below.

\begin{assumption}
	\label{asmp:ATRM}
	Let the optimization problem \eqref{eq:AOP} be given. Assume that there exists a function
	$\Dmax_a:X\rightarrow[0,\infty]$, constants $\Dmax_b, \Dmax_c>0$, a nondecreasing function
	$f:[0,\infty)\rightarrow[0,\infty)$, and constants
	$c \ge 0$, $L,\delta > 0$, $s>1$, and $\eta \in (0,1)$, such that
	for all $x \in X$ with $C(x)>0$ the following properties hold:
	\begin{description}
		\item [monotonicity:]
			\begin{equation}
				\tag{I}
				\label{eq:monotonicity}
				\pred(x,\Delta_1) \ge \pred(x,\Delta_2)\qquad \forall\Delta_1 \ge \Delta_2 > 0
			\end{equation}
		\item[sufficient decrease:]
			\begin{equation}
				\tag{II}
				\label{eq:sufficient_decrease}
				\pred(x,\Delta) \ge f(C(x)) \Delta - c
				\Delta^{s}\qquad\forall
				\Delta \in (0, \Dmax_a(x))
			\end{equation}
		\item[accuracy:]
			\begin{equation}
				\tag{III}
				\label{eq:accuracy}
				\ared(x,\Delta) - \eta\,\pred(x,\Delta) \ge (1-\eta) \bigg(f(C(x)) \Delta
				- c \Delta^s\bigg)\qquad \forall
				\Delta \in (0, \Dmax_a(x))
			\end{equation}
			\begin{equation}
				\tag{III.E}
				\label{eq:accuracy_E}
				\begin{aligned}
					\ared(x,\Delta) &\ge \eta\,\pred(x,\Delta) \\
					\pred(x,\Delta) &\ge \delta  
				\end{aligned} \qquad\forall \Delta\in [\Dmax_a(x), \Dmax_b]
			\end{equation}
		\item[stability:] $\forall \Delta \in (0,\Dmax_c]$:     
			\begin{align}
				\abs{C(T(x,\Delta))-C(x)}&\le L\Delta
				\tag{IV}\label{eq:stability}\\
				\text{or}\qquad\qquad\pred(x,\Delta) &\ge \delta \tag{IV.E}\label{eq:stability_E}
			\end{align}
	\end{description}     
\end{assumption}
With all necessary assumptions in place, \cref{alg:ATRM} can be stated.
\begin{algorithm2e}
	\KwData{$\eta\in (0,1)$,
		$0 < \gamma_1 < 1 \le \gamma_2$,
		$0 < \Delta_0
		< \Delta_{\max} \in (0,\infty]$,
	$x_0 \in X$}
	$n \coloneqq 0$\\
	\While{$C(x_n)>0$}
	{    \eIf{$\ared(x_n,\Delta_n) \ge \eta\
		\pred(x_n,\Delta_n)$}
		{
			$x_{n+1} \coloneqq T(x_n,\Delta_n)$\\
			$\Delta_{n+1} \coloneqq \min\{\gamma_2\Delta_n,\Delta_{\max}\}$
		}
		{
			$x_{n+1} \coloneqq x_n$\\
			$\Delta_{n+1} \coloneqq \gamma_1\Delta_n$
		}
		$n \coloneqq n + 1$
	}
	\caption{Abstract trust-region method}
	\label{alg:ATRM}
\end{algorithm2e}
Under our abstract assumptions, we obtain convergence of the criticality measure.
\begin{theorem}  
	Let the optimization problem \eqref{eq:AOP} be given. If \cref{asmp:ATRM} is valid, then
	\cref{alg:ATRM} either terminates at an iterate $x_n$ with $C(x_n)=0$ or produces
	an infinite sequence of iterates $(x_n)$, such that
	$\lim_{n\rightarrow\infty}C(x_n) = 0$.
	\label{thm:atrm_convergence}
\end{theorem}
The proof of the theorem follows the ideas of \cite[Section~4]{Manns2024}
and can be found in \cref{sec:appendix}.
\begin{remark} 
	In the context of \cref{alg:ATRM}, an iteration $n$ is called ``accepted'' or ``successful''
	if $\ared(x_n,\Delta_n) \ge \eta\,\pred(x_n,\Delta_n)$ and
	subsequently $x_{n+1} = T(x_n,\Delta_n)$. Otherwise iteration $n$ is called
	``unsuccessful''.

	Regarding \cref{asmp:ATRM}\eqref{eq:sufficient_decrease} and \eqref{eq:accuracy}, it is important
	to point out
	that a negative term on the right hand side is allowed, and that this term is
	of higher order than $\Delta$. Whether there are finitely many of such terms with
	different exponents $s_j>1$ ultimately is of no importance. In the worst
	case one can always limit the convergence proof to estimate every
	negative term using $\Delta^{s_i} \le \Delta^{\min_{i}s_i}$ for $\Delta \in
	[0,1]$.

	As it becomes apparent from the proof, the trust-region radius is required to
	be scaled with $\gamma_1<1$ upon an unsuccessful iteration and is
	nondecreasing in case of a successful one, i.e.,
	$\Delta_{n+1} \in [\Delta_n, \min\set{\gamma_2 \Delta_n,
	\Delta_{\max}}]$ could be used. This
	corresponds to the classic trust-region method where the decision on
	increasing the trust-region radius $\Delta$ is made on the basis of how well prediction and actual
	reduction match.
	\label{rem:ATRM_accuracy}
\end{remark}

\section{The proximal gradient algorithm as a trust-region method}
\label{sec:convex_problem}
The main work in the convergence proof for \cref{alg:ATRM} applied to the problem \eqref{eq:OC}
is the verification of \cref{asmp:ATRM} for a composite problem.
Roughly speaking,
the oracle function for \eqref{eq:OC}
includes a proximal gradient step on the set where the control is already active.
Consequently,
the aim of the current section is twofold.
First, it provides
a self-contained convergence proof for a proximal gradient algorithm
by interpreting it as an abstract trust-region method. 
Second,
the estimates of this section will be used in \cref{sec:full_problem}
to analyze the trust-region method for \eqref{eq:OC}.
For this section, we will be investigating problems of class
\begin{equation*}
	\label{eq:OC_conv}
	\min_{v\in U} J(v) \coloneqq F(v) + G(v). 
	\tag{$OC_\text{conv}$}
\end{equation*}
Here, $U$ is a Hilbert space
and $F \colon U \to \R$, $G \colon U \to (-\infty,\infty]$
are given functions.
Throughout this section, we require the following standing assumption.
\begin{assumption}[Standing assumption for \cref{sec:convex_problem}]
	\label{asm:standing_prox_grad}
	The function $G$ is assumed to be convex, proper and lower semicontinuous
	and $F$ is assumed to be Fréchet differentiable with Lipschitz continuous
	derivative with constant $L_{\nabla F}$.
	Further, we require that $F$ and $G$ are Lipschitz continuous on
	$\dom(G) \coloneqq \set{v \in U \given G(v) < \infty}$ with constants $L_F,L_G$.
	Finally, $F + G$ is assumed to be bounded from below.
\end{assumption}
The final assumption on the Lipschitz continuity of $F$ and $G$
on the domain of $G$ is quite strong
and actually not needed for the convergence of a proximal gradient method.
However, it is needed for the straightforward application of the abstract trust-region method
from the previous section, see, in particular, \cref{lem:OC_conv_iterate_in_trust_region}
below.

In \cref{sec:OC_conv_prelims} it will be discussed briefly, why the function 
\begin{equation}
	\tag{$TRM_{\text{conv}}$}
	\label{eq:OC_conv_TRM}
	m_{v,\Delta}(w) \coloneqq \du {\nabla F(v),w - v}\al_{U}
	+ G(w) - G(v) + \frac{1}{2\Delta}\norm{w - v}^2_{U}
	,
\end{equation}
whose minimizer is the proximal gradient step, is
also a reasonable choice as a model function for \cref{alg:ATRM}.
Herein, $v$ is the current iterate, the trust-region radius $\Delta > 0$ serves as a regularization parameter
and the minimizer of the model function will act as the oracle, i.e., the proposal for the next iterate.
Convergence with respect to the criticality measure will be proven in
\cref{subsec:convergence_prox}, see
\cref{thm:OC_conv_convergence}. In \cref{thm:OC_conv_iterate_convergence} it will be
shown how the inclusion of a regularization term leads to strong
convergence of the iterates under the additional assumption of $\nabla F$ being
a compact operator.

\subsection{Preliminaries}
\label{sec:OC_conv_prelims}
\begin{definition}
	\label{def:prox}
	Given $r>0$, the proximal point $\hat v \coloneqq \prox_{rG}(v)$ of $v \in U$
	w.r.t.\ $r G$ is defined as the unique minimizer $\bar v$ of the strongly convex optimization problem
	\begin{equation*}
		\min_{\hat v\in U} \parens*{ G(\hat v) + \frac{1}{2r}\norm{\hat v - v}_U^2 }.
	\end{equation*}
\end{definition}
It is well known that the proximal point $\hat v = \prox_{r G}(v)$
can be characterized via the optimality condition $v - \hat v\in r\partial
G(\hat v)$.
From the literature we recall two crucial properties that
will be used in the proofs.
\begin{corollary}
	\label{cor:prox_properties}
	For all $r>0$ the proximal point mapping is nonexpansive, i.e.,
	\begin{equation}
		\label{eq:prox_nonexpansive}
		\norm{\prox_{rG}(v) - \prox_{rG}(w)}_{U}\le \norm{v-w}_{U}
		\qquad\forall v,w \in U.
	\end{equation}
	Further, given $v,d \in U$ and defining $\Phi \colon (0,\infty) \to [0,\infty)$
	via
	\begin{equation}
		\Phi(r) \coloneqq \norm{\prox_{rG}(v + rd) - v}_U,
		\label{eq:prox_monotone}
	\end{equation}
	the mapping $r\mapsto \Phi(r)$ is nondecreasing, while $r\mapsto
	\frac{\Phi(r)}{r}$ is nonincreasing.
\end{corollary}
\begin{proof}
	\cref{eq:prox_nonexpansive} is proven in \cite[Corollary
	23.11(i)]{Combettes2017}.
	The monotonicity properties related to
	$\Phi$ are proven in the finite-dimensional setting in \cite[Theorem~10.9]{Beck2017},
	and the proof directly carries over to the Hilbert space setting.
\end{proof}

It is well known that the first-order stationarity conditions
at a point $\bar v \in U$
for \eqref{eq:OC_conv}
are given by
$-\nabla F(\bar v) \in \partial G(\bar v)$,
due to the Fréchet differentiability of $F$.
The next, classical result shows that this can be expressed
via a criticality measure involving the proximal map of $G$.
\begin{corollary}
	\label{lem:OC_conv_stationarity}
	For fixed $r > 0$,
	we introduce the criticality measure $C \colon U \to [0,\infty)$ via
	\begin{equation}
		\label{eq:OC_conv_crit}
		C(v) \coloneqq \frac{1}{2r^2} \norm{v - \prox_{rG}(v - r\nabla F(v))}^2_{U}.
	\end{equation}
	Then, a point
	$\bar v \in U$ is stationary if and only if $C(\bar v)=0.$
\end{corollary}
\begin{proof}
	We write the first-order stationarity condition as
	$(\bar v - r \nabla F(\bar v)) - \bar v \in \partial (r G)(\bar v)$.
	Using the characterization of the proximal point in \cref{def:prox},
	this is equivalent to
	$\bar v = \prox_{rG}(\bar v - r\nabla F(\bar v))$, proving the claim.
\end{proof}

Since the (possibly) nonconvex function $F$ is linearized in \eqref{eq:OC_conv_TRM}, the problem
\eqref{eq:OC_conv_TRM} is
strongly convex and computing the proximal point mapping directly yields the
unique minimizer. 
\begin{definition}
	\label{def:OC_conv_oracle}
	The oracle function $T:U\times (0,\infty)\rightarrow U$ for \eqref{eq:OC_conv} is
	defined as 
	\begin{equation*}
		T(v,\Delta) \coloneqq \prox_{\Delta G}(v - \Delta \nabla F(v))
		.
	\end{equation*}
\end{definition}
Since $\Delta$ is an established symbol for the trust-region radius, the
proximal gradient step used for the oracle function reads slightly
unconventional, but will not lead to any ambiguities later.
In a classical trust-region method, a local model is optimized over the trust region,
i.e., the distance between the current iterate
and the proposed next iterate is limited by the trust-region radius $\Delta$.
In the proximal gradient method, the proximity is
only weakly enforced by the quadratic term in the model.
Nonetheless, using the Lipschitz continuity from
\cref{asm:standing_prox_grad} we can prove that
iterates remain in a scaled trust region.
\begin{lemma}
	\label{lem:OC_conv_iterate_in_trust_region}
	There exists a $\trrprox > 0$, such that the oracle
	from \cref{def:OC_conv_oracle}
	satisfies
	\begin{equation*}
		\norm{T(v,\Delta) - v}_{U}\le \trrprox\Delta
		\qquad\forall v \in \dom(G), \Delta \in (0,\infty)
		.
	\end{equation*} 
	In particular, we can choose
	$\trrprox= 2(L_F+L_G)>0$, where $L_F, L_G$ denote the Lipschitz constants
	of $F$ and $G$, both on the domain of $G$.
\end{lemma}
\begin{proof}
	Since $\hat v\coloneqq T(v,\Delta)$ is the minimizer of
	the model function
	$m_{v, \Delta}$ from
	\eqref{eq:OC_conv_TRM},
	we obtain
	\begin{equation*}
		0
		=
		m_{v, \Delta}(v)
		\ge
		m_{v, \Delta}(\hat v )
		=
		\du {\nabla F(v),\hat v- v}\al_{U}
		+ G(\hat v) - G(v) + \frac{1}{2 \Delta}\norm{\hat v- v}^2_{U}.
	\end{equation*}
	Rearranging the terms and using \cref{asm:standing_prox_grad},
	we find
	\begin{equation}
		\label{eq:estimate_prox}
		\frac{1}{2 \Delta}\norm{\hat v- v}^2_{U}
		\le
		\du {\nabla F(v),v - \hat v}\al_{U}
		+ G(v) - G(\hat v)
		\le
		(L_F + L_G) \norm{ \hat v - v }_U
		.
	\end{equation}
	Here, we used
	\begin{equation*}
		\du {\nabla F(v),v - \hat v}\al_{U}
		=
		-\lim_{t \searrow 0} \frac{ F(v + t (\hat v - v)) - F(v) }{t}
		\le
		L_F \norm{\hat v - v}_U
		.
	\end{equation*}
	This yields the claim.
\end{proof}

We briefly mention that the assumed Lipschitz continuity of $G$ is
necessary for the claim of \cref{lem:OC_conv_iterate_in_trust_region}.
As a simple example, we could consider
$U = \R$, $F \equiv 0$ and
\begin{equation*}
	G(x)
	\coloneqq
	\begin{cases}
		-\sqrt{x} & \text{if } x \ge 0, \\
		\infty & \text{else}.
	\end{cases}
\end{equation*}
Then,
\begin{align*}
	T(0, \Delta)
	&=\argmin_{\hat v\in \R} \parens*{ 0\cdot(\hat v - 0) + G(\hat v) - G(0)
	+ \frac{1}{2\Delta}(\hat v - 0)^2 }\\
	&=\argmin_{\hat v \in [0,\infty)} \parens*{-\sqrt{\hat v} + \frac{1}{2\Delta}\hat
	v^2}
	=\left(\frac{\Delta}{2}\right)^{\frac{2}{3}}
	.
\end{align*}
This shows that for any value of $\trrprox > 0$,
the assertion of
\cref{lem:OC_conv_iterate_in_trust_region} is violated for $\Delta
< \frac{1}{4\trrprox^3}$. As
the issue occurs for all trust-region radii below a certain threshold, it
can not be excluded in the algorithm.
\subsection{Convergence}
\label{subsec:convergence_prox}
In this section the individual parts of \cref{asmp:ATRM} are verified
in separate lemmata.
The convergence of the critical measure will
be concluded from \cref{thm:atrm_convergence} in \cref{thm:OC_conv_convergence}.
Strong convergence for a subsequence of iterates is proven under additional
assumptions in
\cref{thm:OC_conv_iterates_converge}.

In view of \cref{lem:OC_conv_iterate_in_trust_region},
the space $X$ appearing in the abstract setting of \cref{sec:ATRM}
is defined as
\begin{equation*}
	X \coloneqq \dom(G).
\end{equation*}
Note that the oracle $T$ maps into $X = \dom(G)$
by definition of the proximal mapping.

\begin{definition}[Prediction function]
	\label{def:OC_conv_prediction_function}
	We define the prediction function $\pred:U\times (0,\infty) \to [0,\infty)$ using the
	model function \eqref{eq:OC_conv_TRM} via
	\begin{equation*}
		\pred(v,\Delta) \coloneqq m_{v,\Delta}(v) - m_{v,\Delta}(T(v,\Delta))
		= -m_{v,\Delta}(T(v,\Delta)).
	\end{equation*}
\end{definition}
\begin{lemma}[Monotonicity]
	\label{lem:OC_conv_monotonicity}
	The prediction function $\pred(v,\Delta)$ is nondecreasing with respect to
	$\Delta$.
\end{lemma}
\begin{proof}
	As $\Delta$ increases, the model function $m_{v, \Delta}$ decreases.
	This yields the claim.
\end{proof}
\begin{lemma}[Sufficient decrease] 
	\label{lem:OC_conv_predictability}
	For the value $r > 0$ that was used
	in the definition  \eqref{eq:OC_conv_crit} of the criticality measure, the estimate
	\begin{equation*}
		\pred(v,\Delta) = -m_{v,\Delta}(T(v,\Delta))  \ge C(v)\Delta 
		\qquad\forall v \in U, \Delta \in (0,r]
	\end{equation*}
	holds and therefore \cref{asmp:ATRM}\eqref{eq:sufficient_decrease} is valid with
	$\Dmax_a \equiv r$ and $c = 0$.
\end{lemma}
\begin{proof}
	For brevity, we set $\hat v \coloneqq T(v, \Delta)$. Due to the convexity of $G$,
	the model function $m_{v,\Delta}$
	is strongly convex with constant $\Delta^{-1}$.
	Consequently, the quadratic growth condition
	\begin{equation*}
		m_{v,\Delta}(w)
		\ge
		m_{v,\Delta}(\hat v)
		+
		\frac{1}{2\Delta} \norm{w - \hat v}_U^2
		\qquad\forall w \in U
	\end{equation*}
	holds at its minimizer $\hat v \coloneqq T(v, \Delta) = \prox_{\Delta
	G}(v-\Delta\nabla F(v))$.
	In particular, inserting $w = v$ implies
	\begin{equation*}
		\pred(v, \Delta)
		= 
		m_{v,\Delta}(v)
		-
		m_{v,\Delta}(\hat v)
		\ge
		\frac{1}{2\Delta} \norm{v - \hat v}_U^2
		.
	\end{equation*}
	It remains to
	relate the last term with the criticality measure
	\eqref{eq:OC_conv_crit}. To
	this end, we employ the function $r\mapsto \Phi(r) = \norm{
	v - \prox_{rG}(v - r\nabla F(v))}$ from \cref{cor:prox_properties} and the
	fact that $r\mapsto \frac{\Phi(r)}{r}$ is nonincreasing.
	We get
	\begin{align*}
		\frac{1}{2\Delta}\norm{v - \hat v}^2_U &=
		\frac{1}{2\Delta}\norm{v - \prox_{\Delta G}(v - \Delta\nabla F(v))}^2_U
		=\frac{1}{2\Delta}\Phi(\Delta)^2
		=\frac{\Delta}{2}\left( \frac{\Phi(\Delta)}{\Delta} \right)^2\\
		&\ge \frac{\Delta}{2} \left( \frac{\Phi(r)}{r} \right)^2
		=\frac{\Delta}{2 r^2}\norm{\prox_{rG}(v - r\nabla F(v))-v}^2_{U} = C(v) \Delta
		,
	\end{align*}
	closing the proof.
\end{proof}
\begin{lemma}[Accuracy]
	\label{lem:OC_conv_accuracy}
	For an arbitrary $\eta \in (0,1]$,
	we have
	\begin{equation*}
		\ared(v,\Delta) - \eta\,\pred(v,\Delta) \ge (1-\eta) C(v) \Delta
		- \frac12 L_{\nabla F}\trrprox^2\Delta^2
		\qquad\forall v \in \dom(G), \Delta \in (0,r]
		,
	\end{equation*}
	where $r > 0$ is from \eqref{eq:OC_conv_crit} and
	$L_{\nabla F}$ denotes the Lipschitz constant of $\nabla F$.
	Thus,
	\cref{asmp:ATRM}\eqref{eq:accuracy} holds for $\Dmax_a \equiv r$, $s=2$, and
	$c = \frac12 L_{\nabla F}\trrprox^2$.
\end{lemma}
\begin{proof}
	We define $\hat v\coloneqq T(v,\Delta)$. Since $\nabla F$ is assumed to be Lipschitz continuous with constant $L_{\nabla F}$,
	the fundamental theorem of calculus implies
	\begin{equation*}
		F(\hat v)
		\le
		F(v) + \du{\nabla F(v), \hat v - v}\al_U
		+ \frac{L_{\nabla F}}{2} \norm{ \hat v - v }_U^2
		.
	\end{equation*}
	Consequently, 
	\bgroup\allowdisplaybreaks
	\begin{align*}
		\ared(v, \Delta)
		&=
		J(v) - J(\hat v)
		=
		F(v) - F(\hat v) + G(v) - G(\hat v)
		\\
		&\ge
		\du{\nabla F(v), v - \hat v}\al_U
		+ G(v) - G(\hat v)
		- \frac{L_{\nabla F}}{2} \norm{ \hat v - v }_U^2
		\\
		&=
		\pred(v, \Delta) + \frac12 \parens*{\frac1\Delta - L_{\nabla F}}
		\norm{\hat v - v}_U^2\\
		&\ge \eta \pred(v,\Delta) + (1-\eta)C(v)\Delta
		- \frac{1}{2} L_{\nabla F} \trrprox^2 \Delta^2
		,
	\end{align*}
	\egroup
	where we used \cref{lem:OC_conv_iterate_in_trust_region,lem:OC_conv_predictability}
	in the final estimate.
\end{proof}
\begin{remark}
	\label{rem:accuracy}
	The estimate in the proof of \cref{lem:OC_conv_accuracy}
	even yields
	\begin{equation*}
		\Delta \le L_{\nabla F}^{-1}
		\quad\implies\quad
		\frac{1}{\Delta}
		- L_{\nabla F} \ge 0
		\quad\implies\quad
		\ared(v,\Delta) - \eta
		\pred(v,\Delta) \ge 0
		.
	\end{equation*}
	This means that iterations with a sufficiently small trust-region radius
	are always successful. Such a result is expected under the given
	assumptions, cf.\ \cite[Lemma~10.4]{Beck2017}.
\end{remark}
\begin{lemma}[Stability]
	\label{lem:OC_conv_stability}
	The criticality measure $C(v)$ from
	\eqref{eq:OC_conv_crit}
	satisfies
	\begin{equation*}
		\abs{C(T(v,\Delta))-C(v)} \le
		L_C \Delta
		\qquad
		\forall v \in \dom(G), \Delta \in (0,\infty),
	\end{equation*}
	where $L_C\coloneqq \frac{\trrprox^2}{r}(2+rL_{\nabla F})$ with $L_{\nabla
	F}$ being the Lipschitz constant of $\nabla F$ on $\dom(G)$
	and $\trrprox$ as defined in \cref{lem:OC_conv_iterate_in_trust_region}.
	Hence, \cref{asmp:ATRM}\eqref{eq:stability} holds for arbitrary $\Dmax_c > 0$.
\end{lemma}
\begin{proof}
	Again, we set $\hat v \coloneqq T(v, \Delta)$.
	We first recall that
	\cref{lem:OC_conv_iterate_in_trust_region} implies
	\begin{equation*}
		\norm{v - \prox_{rG}(v - r\nabla F(v))}_U
		=
		\norm{v - T(v, r)}_U
		\le
		\trrprox r,
	\end{equation*}
	and similarly for $v$ replaced by $\hat v$.
	This lemma also yields $\norm{v - \hat v}_U \le \trrprox \Delta$.
	Further,
	we utilize
	the nonexpansivness of the proximal mapping as stated in \cref{cor:prox_properties}
	to estimate
	\begin{align*}
		2r^2\abs{C(\hat v) - C(v)}
		&= \abs[\big]{\norm{\hat v - \prox_{rG}(\hat v - r\nabla
			F(\hat v))}^2_U - \norm{v - \prox_{rG}(v - r\nabla
		F(v))}^2_U}\\
		&= (\norm{\hat v - \prox_{rG}(\hat v - r\nabla
			F(\hat v))}_U 
			+ \norm{v - \prox_{rG}(v - r\nabla
		F(v))}_U)\\
		&\qquad \cdot \abs[\big]{ \norm{\hat v - \prox_{rG}(\hat v - r\nabla
			F(\hat v))}_U - \norm{v - \prox_{rG}(v - r\nabla
		F(v))}_U}\\
		&\le 2\trrprox r\norm{v - \hat v + \left( \prox_{rG}(\hat v - r\nabla
				F(\hat v)) - \prox_{rG}(v - r\nabla
		F(v)) \right)}_U\\
		&\le 2 \trrprox r \left( \trrprox \Delta + \norm{\prox_{rG}(\hat v - r\nabla
				F(\hat v)) - \prox_{rG}(v - r\nabla
		F(v))}_U\right)\\
		&\le 2 \trrprox r\left(\trrprox \Delta + \norm{\hat v - r\nabla F(\hat v) - (v -r\nabla
		F(v))\right) }_U\\
		&\le 2 \trrprox r\left(2 \trrprox \Delta + r\norm{\nabla F(\hat v) - \nabla
		F(v)}_U\right)
		\le 2 \trrprox^2 r(2+rL_{\nabla F}) \Delta. 
	\end{align*}
\end{proof}
\begin{theorem}
	The abstract trust-region method \cref{alg:ATRM} applied to \eqref{eq:OC}
	with the proximal gradient step as oracle $T(v,\Delta)$ and the prediction
	function $\pred(v,\Delta)$ from
	\cref{def:OC_conv_prediction_function} either terminates at a stationary point
	$\bar v$ or produces an infinite sequence of iterates $(v_n)$, such that
	\begin{equation*}
		\norm*{
			v_n
			-
			\prox_{\hat r G}\left(v_n - \hat r \nabla F(v_n) \right)
		}_U
		\to 0
	\end{equation*}
	for all $\hat r > 0$.
	\label{thm:OC_conv_convergence}
\end{theorem}
\begin{proof}
	We will apply \cref{thm:atrm_convergence}
	with $X = \dom(G)$.
	Monotonicity was proven in \cref{lem:OC_conv_monotonicity}, predictability in \cref{lem:OC_conv_predictability}, accuracy in \cref{lem:OC_conv_accuracy} and stability in \cref{lem:OC_conv_stability}. 
	This verifies \cref{asmp:ATRM} with
	$\Dmax_a(v) = \Dmax_c = r$,
	$\Dmax_b = r/2 < r$.
	Note that none of the exception clauses were needed, i.e.,
	\eqref{eq:accuracy_E} is void due to $\Dmax_b < \Dmax_a$
	and
	\eqref{eq:stability_E} is not needed since
	\eqref{eq:stability} always holds.
	Therefore, convergence of the criticality measure $C$ for the iterates of the trust-region
	algorithm \cref{alg:ATRM} follows from \cref{thm:atrm_convergence}.
	This shows the convergence for $\hat r = r$.
	In particular, the function $\Phi_n(\hat r) \coloneqq \norm{u_n - \prox_{\hat r G}(u_n - \hat r \nabla F(u_n))}_U$
	satisfies $\Phi_n(r) \to 0$.
	Using the monotonicity properties from \cref{cor:prox_properties},
	this implies
	\begin{equation*}
		\Phi_n(\hat r) \le \Phi_n(r) \to 0 \quad\forall \hat r < r
		\qquad\text{and}
		\qquad
		\frac{\Phi_n(\hat r)}{\hat r} \le \frac{\Phi_n(r)}{r} \to 0 \quad\forall \hat r > r
		.
	\end{equation*}
	This finishes the proof.
\end{proof}
\begin{remark}
	In \cref{rem:ATRM_conv_with_Delta_bound} it will be argued that convergence with
	respect to the criticality measure already
	follows from \cref{asmp:ATRM} \eqref{eq:monotonicity}--\eqref{eq:accuracy} if
	a global lower bound of the trust-region radius is known. With
	\cref{rem:accuracy} convergence of \cref{alg:ATRM} for \eqref{eq:OC} with the
	proximal gradient oracle already follows from
	\cref{lem:OC_conv_monotonicity}, \cref{lem:OC_conv_predictability} and
	\cref{lem:OC_conv_accuracy}. Notably, global Lipschitz continuity of $F$ and $G$ is
	not needed.
	Regarding regularity of $F$, we still need Lipschitz continuity of $\nabla F$. This recovers the
	standard convergence result for the proximal gradient method,
	cf.\ \cite[Theorem~10.15(b)]{Beck2017}.
	\label{rem:OC_TRM_conv_Delta_bound}
\end{remark}

The result in this form does not yet yield convergence of the iterates.
This can be shown under
additional assumptions on the function $F$.
\begin{theorem}
	\label{thm:OC_conv_iterate_convergence}
	We assume that $G \colon U \to (-\infty,\infty]$ is convex, proper and lower semicontinuous.
	Let $\hat F \colon U \to \R$ be Fréchet differentiable
	such that
	$\nabla\hat F \colon U \to U$ is compact and
	Lipschitz continuous.
	We further assume that $\hat F$ and $G$ are Lipschitz continuous on $\dom(G)$.
	We define $F$ via
	\label{thm:OC_conv_iterate_convergence_F}
	\begin{equation*}
		F(v) \coloneqq \hat F(v) + \frac{\sigma}{2}\norm{v}^2_U
		\qquad\forall v \in U,
	\end{equation*}
	where
	$\sigma > 0$.
	Finally, we assume that $F + G$ is bounded from below.
	Then, the sequence of iterates
	contains a strongly convergent subsequence and every accumulation point is stationary.
	\label{thm:OC_conv_iterates_converge}
\end{theorem}
\begin{proof}
	We have to check that \cref{asm:standing_prox_grad}
	is satisfied for $F$ and $G$.
	The only problem is that $v \mapsto \norm{v}_U^2$ is
	only Lipschitz continuous on bounded sets.
	For this, we use that $F(v_n) + G(v_n)$ is decreasing,
	in particular,
	\begin{equation*}
		\norm{v_n}_U
		\le
		\sqrt{2\cdot\frac{F(v_n) + G(v_n)- C}{\sigma}}
		\le
		\sqrt{2\cdot\frac{F(v_0) + G(v_0)- C}{\sigma}}
		=:
		R
		,
	\end{equation*}
	where $C$ is a lower bound for $F + G$.
	This means that we can redefine $G$ to take the value $+\infty$
	outside the ball of radius $R$
	without changing the iterates of \cref{alg:ATRM}.
	Consequently, $F$ becomes Lipschitz continuous on $\dom(G)$.

	From the boundedness of the iterates $(v_n)$ we
	get weak
	convergence of a subsequence. However, it is not clear that the weak
	limit $\bar v$ satisfies $C(\bar v) =0$. We will first prove that every weakly convergent
	subsequence is even strongly convergent and conclude that the limit point $\bar v$ is stationary.

	Consider (without relabeling) a weakly convergent subsequence $(v_n)\rightharpoonup
	\bar v$.
	As shorthand we introduce $\hat v_n$ via a proximal gradient step with
	step length $\frac{1}{\sigma}$, i.e.,
	\begin{equation*}
		\hat v_n \coloneqq \prox_{\frac{G}{\sigma}}\left(v_n - \frac{\nabla
		F(v_n)}{\sigma}\right).
	\end{equation*}
	From \cref{thm:OC_conv_convergence},
	we get
	$\norm{\hat v_n- v_n} \rightarrow 0$.
	This implies $\hat v_n\rightharpoonup \bar v$.
	The proximal point mapping is characterized by the optimality condition for the corresponding
	optimization problem,
	i.e.,
	\begin{equation*}
		v_n - \hat v_n - \frac{1}{\sigma}\nabla F(v_n) \in \frac1\sigma \partial G(\hat v_n)
		.
	\end{equation*}
	Using the definition of $F$, this yields
	\begin{equation*}
		- \frac{1}{\sigma}\nabla \hat
		F(v_n) \in \partial
		\left(\frac{1}{\sigma}G+\frac{1}{2}\norm{\cdot}^2\right)(\hat
		v_n)
		.
	\end{equation*}
	Due to the compactness of $\nabla\hat F$, the left-hand side converges to $-\frac1\sigma \nabla \hat F(\bar v)$.
	Since the function $\frac{1}{\sigma}G+\frac{1}{2}\norm{\cdot}^2$
	is strongly convex,
	its subdifferential has a Lipschitz continuous inverse,
	see \cite[Theorem~18.15 (vii)$\Rightarrow$(i)]{Combettes2017}.
	This implies the strong convergence $\hat v_n \to \bar v$
	and, consequently, $v_n \to \bar v$.
	By passing to the limit, we also get
	\begin{equation*}
		-\frac1\sigma\nabla F(\bar v)
		=
		- \bar v -\frac{1}{\sigma}\nabla \hat F(\bar v)\in \frac{1}{\sigma}\partial G(\bar v)
		.
	\end{equation*}
	This shows that $\bar v$ is a stationary point.
\end{proof}
\section{Control problems with continuous-or-off controls and TV regularization}
%%fakesubsection: Intro
\label{sec:full_problem}
In this section, we will apply the trust-region algorithm
to the problem
\begin{equation*}
	\label{eq:OC}
	\tag{OC}
	\min_{u\in \Uad} J(u) \coloneqq F(u) + \sum_{i=1}^{N}\left(\int_{t_0}^{T}   g_i(u_i(t))\dt
	+ \TV(\sgn(u_i))\right)
	.
\end{equation*}
The notation and assumptions are introduced in \cref{subsec:setting_assumptions}
and we also prove existence of solutions therein.
In \cref{subsec:stationarity}, we derive stationarity conditions
and define the criticality measure.
The application of \cref{sec:ATRM}
is done in \cref{subsec:TR_application}.
Finally,
we discuss the fast solution of the trust-region subproblems
in \cref{subsec:solution_subproblem}.

\subsection{Problem setting, assumptions and existence of solutions}
\label{subsec:setting_assumptions}
In problem \eqref{eq:OC},
the control variable is vector-valued and non-negative,
i.e., $u \colon [t_0, T] \to [0,\infty)^N$.
By $L^p(t_0, T; \R^N)$ we denote the usual Lebesgue space
of measurable functions $u \colon [t_0, T] \to \R^N$ satisfying
\begin{equation*}
	\norm{u}_{L^p(t_0, T; \R^N)}^p
	\coloneqq
	\int_{t_0}^T
	\norm{u(t)}_p^p \d t
	=
	\int_{t_0}^T
	\sum_{i = 1}^N \abs{u_i(t)}^p
	\d t
	<
	\infty,
\end{equation*}
with the usual modification for $p = \infty$.

One of the main features of \eqref{eq:OC} is the total variation term
which acts on $\sgn(u_i)$, which is defined pointwise in time,
i.e., $\sgn(u_i)(t) = 1$ if $u_i(t) > 0$ and $\sgn(u_i)(t) = 0$ if $u_i(t) = 0$.
We extend a function $x \in \BV(t_0,T; \set{0,1})$ by $0$ to
the interval $[t_0 - \varepsilon,T+\varepsilon]$ (for some arbitrary $\varepsilon > 0$)
and define
\begin{equation*}
	\TV(x) \coloneqq
	\sup\set*{ \int_{t_0-\varepsilon}^{T+\varepsilon}x(t)
		\varphi'(t)\dt\given \varphi\in C_c^1(t_0-\varepsilon,T+\varepsilon; [-1,1])
	}.
\end{equation*}
In contrast to the usual definition,
our definition of $\TV(x)$ also includes possible jumps at $t_0$ and $T$.
If we consider \eqref{eq:OC}, the term $\TV(\sgn(u_i))$ counts how often the control $u_i$ is switched on or off,
including the boundary points $t_0$ and $T$.
Note that $\TV(\sgn(u_i))$ is always an even number.
In the vector-valued case $x \in \BV(t_0, T; \set{0,1}^N)$,
we set $\TV(x) \coloneqq \sum_{i = 1}^N \TV(x_i)$.

The admissible set is given by
\begin{equation}
	\Uad
	\coloneqq
	\set*{
		u \in L^\infty(t_0, T; \R^N)
		\given
		\begin{aligned}
			&u_i(t) \in \set{0} \cup [a_i, b_i] \text{ for a.a. $t \in [t_0, T]$}
			\text{ and }
			\\
			&\TV(\sgn(u_i)) < \infty \text{ for all $i = 1,\ldots,N$}
		\end{aligned}
	},
	\label{def:notation}
\end{equation}
where $0 < a_i \le b_i$.
For convenience, we also define the sets
\begin{align*}
	\swp &\coloneqq \bigtimes_{i = 1}^N [a_i, b_i] , &
	\swpz&\coloneqq \bigtimes_{i = 1}^N (\{0\}\cup[a_i,b_i]),
\end{align*}
which encode
the possible range of the controls.

The assumptions on the objective terms $F$ and $g_i$
are collected next.
\begin{assumption}
	\label{asmp:OC}\phantom{i}
	\begin{enumerate}
		\item \label{asmp:F}
			We assume that $F \colon L^1(t_0, T; \R^N) \to \R$ is bounded from below and Fréchet
			differentiable
			such that the Fréchet derivative
			$\nabla F \colon L^1(t_0, T; \R^N) \to L^\infty(t_0, T; \R^N)$
			is Lipschitz continuous with constant $L_{\nabla F} < \infty$.
			Further, we assume that for all $u \in L^1(t_0,T; \R^N)$,
			the function $\nabla F(u) \colon [t_0, T] \to \R^N$
			is Lipschitz continuous with constant $L_{\nabla F(u)} < \infty$
			and that these constants are uniformly bounded.
			For simplicity of the presentation,
			we denote the bound again by $L_{\nabla F}$, i.e.,
			we assume $L_{\nabla F(u)} \le L_{\nabla F}$
			for all $u \in L^1(t_0,T; \R^N)$.

		\item \label{asmp:g} The function $g_i: \R  \rightarrow (-\infty,\infty]$ has the value $g_i(0) = 0$,
			while $g|_{[a_i,b_i]}$ is assumed to be strongly convex with constant $m
			> 0$ and Lipschitz continuous with constant $L_g$
			on $[a_i, b_i]$. Finally, $g_i|_{\R\setminus (\{0\}\cup [a_i,b_i])}
			= \infty$. 
		\item
			The control bounds satisfy $0<a\le a_i \le b_i \le b <\infty$
			for constants $a$ and $b$
	\end{enumerate}
\end{assumption}
\begin{remark}
	\label{rem:F_Lipschitz}
	Clearly,
	$\norm{u}_{L^1(t_0,T;\R^n)} \le b (T-t_0) N$
	for all $u \in \Uad$.
	The assumed Lipschitz continuity of $\nabla F$ with respect to the control directly implies
	\begin{equation}
		\label{eq:bound_nabla_F}
		\begin{aligned}
			\norm{\nabla F(u)}_{L^\infty(t_0,T;\R^N)}
			&\le
			\norm{\nabla F(u) - \nabla F(0)}_{L^\infty(t_0,T;\R^N)}
			+
			\norm{\nabla F(0)}_{L^\infty(t_0,T;\R^N)}
			\\
			&\le
			C_{\nabla F}
			\coloneqq
			L_{\nabla F} b (T-t_0) N
			+
			\norm{\nabla F(0)}_{L^\infty(t_0,T;\R^N)}
			< \infty
		\end{aligned}
	\end{equation}
	for all $u\in \Uad$.
	Subsequently the first-order Taylor expansion with
	integral remainder yields 
	\begin{equation}
		\label{eq:F_Lipschitz}
		\abs{F(u) - F(\hat u)} \le C_{\nabla F}\norm{u-\hat u}_{L^1(t_0,T;\R^N)}
		\qquad \forall u,\hat u \in \Uad.
	\end{equation}
\end{remark}

Note that the assumption $g_i(0) = 0$
is not restrictive and can be achieved by a simple shift.
For brevity,
we define
the functional $G \colon L^2(t_0, T; \R^N) \to (-\infty,\infty]$
via
\begin{equation*}
	G(u)
	\coloneqq
	\sum_{i=1}^{N}\int_{t_0}^{T} g_i(u_i(t))\dt
	.
\end{equation*}
In particular, we have that $G(u) < \infty$
implies
$u \in \swpz$ a.e.\ on $[t_0, T]$.
Note that with our notation, the objective in \eqref{eq:OC}
can now be written as
\begin{equation*}
	F(u) + G(u) + \TV(\sgn(u))
	.
\end{equation*}
Finally, we need variants of $g_i$ and $G$, which restrict their arguments to be non-zero throughout,
i.e.,
we define $g_{i,+} \colon \R \to (-\infty, \infty]$
via
\begin{equation*}
	g_{i,+}(x) \coloneqq
	\begin{cases}
		g_i(x) & \text{if } x \in [a_i, b_i],\\
		\infty & \text{else}
	\end{cases}
\end{equation*}
and
$G_+ \colon L^2(t_0, T; \R^N) \to (-\infty,\infty]$
via
\begin{equation*}
	G_+(u)
	\coloneqq
	\sum_{i=1}^{N}\int_{t_0}^{T} g_{i,+}(u_i(t))\dt
	=
	\begin{cases}
		G(u)&\text{if } u(t)\in \swp \text{ for a.a.\ } t \in (t_0, T), \\
		\infty&\text{else}.
	\end{cases}
\end{equation*}
Note that $g_i$ and $G$ are nonconvex while $g_{i,+}$ and $G_+$
are convex functions,
due to the convexity of $g_i$ on $[a_i, b_i]$.

In order to prove the existence of solutions,
we provide a lemma.
\begin{lemma}
	\label{lem:compatibility_limits}
	Let $(u_k) \subset \Uad$ be a sequence of admissible controls
	such that
	$u_k \weakly \bar u$ in $L^2(t_0,T; \R^N)$
	and
	$\sgn(u_k) \weaklystar \bar\alpha$ in $\BV(t_0,T; \R^N)$,
	where weak-$\star$ convergence in $\BV(t_0,T; \R^N)$ is defined as in \cite[Definition~3.11]{Ambrosio2000}.
	Then, $\sgn(\bar u) = \bar\alpha$
	and $\bar u \in \Uad$.
\end{lemma}
\begin{proof}
	For simplicity of notation,
	we only argue in the scalar case $N = 1$.

	From \cite[Corollary~3.49]{Ambrosio2000},
	we get $\sgn(u_k) \to \bar\alpha$ in $L^1(t_0,T)$.
	We choose
	a subsequence (without relabeling)
	such that
	$\norm{\sgn(u_{k})-\bar\alpha}_{L^1(t_0,T)}\le 2^{-k}$.
	Since $\sgn(u_{k})$ and $\bar \alpha$ have values in $\{0,1\}$ almost everywhere,
	the Lebesgue measure of the set on which they differ
	equals the $L^1(t_0,T)$ norm of their difference,
	i.e.,
	\begin{align*}
		\abs*{\set{\sgn(u_{k})\neq
		\bar \alpha}}
		% =\abs*{\set{\sgn(u_{k})=1} \mathbin{\triangle} \set{\bar \alpha =1}}
		=  \norm{\sgn(u_{k}) - \bar \alpha}_{L^1(t_0,T)}\le 2^{-k}.
	\end{align*}
	On the set
	$\bigcap_{l = k}^{\infty} \set{\sgn(u_{l}) = 1}$
	we have $b \ge u_{l} \ge a > 0$ for all $l \ge k$ and, consequently, $b \ge \bar u \ge a > 0$ a.e.
	This shows $\sgn(\bar u) = 1 = \bar\alpha$ and $\bar u \in [a,b]$ on this set.
	Similarly, one can argue on the set
	$\bigcap_{l = k}^{\infty} \set{\sgn(u_{l}) = 0}$
	and one obtains
	$\sgn(\bar u) = 0 = \bar\alpha$ and
	$\bar u = 0$ a.e.\ on this set.
	Thus,
	we get
	$\sgn(\bar u) = \bar\alpha$ and
	$\bar u \in \set{0} \cup [a,b]$
	a.e.\ on the set
	$\bigcap_{l = k}^{\infty} \set{\sgn(u_{l}) = \bar\alpha}$
	for all $k \in \N$.
	For the Lebesgue measure of the complement, we get
	\begin{equation*}
		\abs*{
			(t_0, T)
			\setminus
			\bigcap_{l = k}^{\infty} \set{\sgn(u_{l}) = \bar\alpha}
		}
		\le
		\sum_{l = k}^\infty
		\abs{\set{\sgn(u_{l}) \ne \bar\alpha} }
		\le
		2^{-k + 1}
		.
	\end{equation*}
	This shows
	$\sgn(\bar u) = \bar\alpha$ and
	$\bar u \in \set{0} \cup [a,b]$
	a.e.\ on $(t_0, T)$.
	In order to get $\bar u \in \Uad$,
	we additionally need $\TV(\sgn(\bar u)) < \infty$,
	but this follows from
	$\TV(\sgn(\bar u)) = \TV(\bar\alpha) \le \liminf_{k \to \infty}
	\TV(\sgn(u_{k}))$,
	see
	\cite[Remark~3.5]{Ambrosio2000}.
\end{proof}
\begin{theorem}
	\label{thm:existence}
	We assume that $F$ is sequentially weakly lower semicontinuous on $L^2(t_0, T; \R^N)$.
	Then, problem \eqref{eq:OC} has a solution.
\end{theorem}
\begin{proof}
	Apply the direct method of calculus of variations.
	A minimizing sequence $(u_k)$ is bounded in $L^2(t_0,T;\R^N)$,
	since $\Uad$ is bounded in this space.
	Further, $(\sgn(u_k))$ is bounded in $\BV(t_0,T; \set{0,1}^N)$
	due to the structure of the objective,
	since $F$ and $G$ are bounded from below.
	Consequently, we can extract a subsequence (without relabeling)
	such that
	$u_k \weakly \bar u$ in $L^2(t_0,T; \R^N)$
	and
	$\sgn(u_k) \weaklystar \bar\alpha$ in $\BV(t_0,T; \R^N)$,
	see \cite[Theorem~3.23]{Ambrosio2000}.
	From \cref{lem:compatibility_limits}, we get
	$\sgn(\bar u) = \bar\alpha$ and $\bar u \in \Uad$.
	It remains to pass to the limit in the objective.
	By assumption, we have
	\begin{equation*}
		\liminf_{k \to \infty} F(u_k)
		\ge
		F(\bar u),
	\end{equation*}
	and
	\begin{equation*}
		\liminf_{k \to \infty} \TV(\sgn(u_k))
		\ge
		\TV(\bar\alpha)
		=
		\TV(\sgn(\bar u))
	\end{equation*}
	follows from the properties of $\TV$, see \cite[Remark~3.5]{Ambrosio2000}.
	It remains to discuss $G$.
	To this end, we define the function $\tilde g_i \colon \R \to \R$ via
	\begin{equation*}
		\tilde g_i( z )
		=
		\begin{cases}
			g_i(z) & \text{if } z \in [a_i, b_i], \\
			g_i(a_i) - (z - a_i) L_g & \text{if } z \le a_i, \\
			g_i(b_i) + (z - b_i) L_g & \text{if } z \ge b_i.
		\end{cases}
	\end{equation*}
	Since $L_g$ is a Lipschitz constant of $g_i$ on $[a_i, b_i]$,
	one can check that $\tilde g_i$ is a convex function.
	For every $i \in \set{1, \ldots, N}$, the sequence $\tilde g_i(u_{k,i})$
	is bounded in $L^2(t_0, T)$, thus, w.l.o.g., we have $\tilde g_i(u_{k,i}) \weakly v_i$ in $L^2(t_0,T)$.
	We have $(u_{k,i}, \tilde g_i(u_{k,i})) \in \epi(\tilde g_i)$ a.e.\ on $[t_0, T]$.
	Since the epigraph $\epi(\tilde g_i) \subset \R^2$ is closed and convex, the weak convergence implies
	$(\bar u_i, v_i) \in \epi(\tilde g_i)$ a.e.\ on $[t_0,T]$, i.e., $v_i \ge \tilde g_i(\bar u_i)$.
	Since $\sgn(u_{k,i}) \to \sgn(\bar u_i)$ in $L^2(t_0,T)$, see \cite[Corollary~3.49]{Ambrosio2000},
	we get
	\begin{equation*}
		g_i(u_{k,i})
		=
		\tilde g_i(u_{k,i}) \sgn( u_{k,i} )
		\weakly
		v_i \sgn( \bar u_i )
		\ge
		\tilde g_i(\bar u_i) \sgn( \bar u_i )
		=
		g_i(\bar u_i)
	\end{equation*}
	in $L^1(t_0, T)$.
	This shows
	\begin{equation*}
		\lim_{k \to \infty} G_i(u_{k,i})
		=
		\int_{t_0}^T v_i \sgn( \bar u_i ) \d t
		\ge
		G_i(\bar u_i).
	\end{equation*}
	By combining the above estimates, we get
	\begin{equation*}
		\liminf_{k \to \infty} J(u_k)
		\ge
		J(\bar u).
	\end{equation*}
	Since $(u_k)$ is a (subsequence of a) minimizing sequence, this shows that $\bar u$ is a minimizer.
\end{proof}

\subsection{Stationarity conditions and criticality measure}
\label{subsec:stationarity}
In this section, we will first
investigate local optimality conditions
for \eqref{eq:OC} by combining theory from convex optimal control with results from \cite[Proposition~4.4]{Leyffer2022}.
These stationarity conditions will be used to define the criticality measure associated with problem \eqref{eq:OC}.

\begin{definition} [Minimal representation] 
	\label{def:minimal_representation}
	Let $u\in L^\infty(t_0,T; \set{0} \cup [a,b])$ be given such that
	$\TV(\sgn(u))<\infty$.
	We set $K \coloneqq \TV(\sgn(u)) / 2$.
	A vector $\hat t\in \R^{2 K}$, such that
	\begin{equation*}
		t_0 \le \hat t_1 < \hat t_2 < \cdots
		< \hat t_{2 K} \le T
	\end{equation*}
	and such that
	\begin{equation}
		\label{eq:representation}
		\sgn(u) = \sum_{j=1}^{K}\chi_{[\hat t_{2 j-1},\hat t_{2 j}]}
	\end{equation}
	holds a.e.\ on $(t_0, T)$ is called minimal representation.
\end{definition}
\begin{proposition}
	\label{prop:minimal_representation}
	The minimal representation from \cref{def:minimal_representation} is
	uniquely determined.
\end{proposition}
\begin{proof}
	This follows by adapting the arguments in
	\cite[Proposition~4.4]{Leyffer2022}.
\end{proof}
In what follows, we always redefine $u$
on a set of measure zero such that \eqref{eq:representation}
holds everywhere on $[t_0, T]$.
This implies
\begin{equation*}
	\{\sgn(u_i) = 1\} \coloneqq
	\set{ t \in [t_0,T] \mid \sgn(u_i(t)) = 1 }
	=
	\bigcup_{j=1}^{K_j}[t_{2i-1},t_{2i}]
	.
\end{equation*}
We refer to $\hat t_u \coloneqq \hat t$ as the minimal
representation of $\sgn(u)$.

In case of a vector-valued function $u \in L^\infty(t_0, T; \swpz)$
with $\TV(\sgn(u)) < \infty$,
we say that
$\hat t_u \coloneqq (\hat t_{u_1}, \ldots, \hat t_{u_N})$ is
the minimal representation of $\sgn(u)$.

We remind the reader, that the jump set of a function $\alpha$ of bounded variation
is commonly referred to as $J_\alpha$.
Since we are interested in the switching points
of functions $u \in L^\infty(t_0,T)$ with $\TV(\sgn(u)) < \infty$, which are the jump points of $\sgn(u)$,
with some abuse of notation, we denote by $J_u$ the jump set of $\sgn(u)$.
Note that the jump set of $\sgn(u)$ satisfies
\begin{equation*}
	J_u = \set{ \hat t_1, \ldots, \hat t_{2 K} },
\end{equation*}
where $\hat t = \hat t_u \in \R^{2 K}$
is the minimal representation of $\sgn(u)$.
This directly yields
\begin{equation*}
	J_{u}\coloneqq \set{s\in[t_0,T]\given\sgn(u)(s-)\neq
	\sgn(u)(s+)},
\end{equation*}
with the conventions
\begin{equation*}
	\sgn(u)(t_0-) \coloneqq \sgn(u)(T+) \coloneqq 0.
\end{equation*}
Note that we used one-sided limits of $\sgn(u)$,
since the function $u$ itself might not possess one-sided limits.
We mention that this characterization of $J_u$ implies that $t_0 \in J_{u}$
if and only if $\sgn(u(t_0+)) = 1$,
and similarly for $T$.

Since we are interested in defining an appropriate criticality
measure, we investigate the necessary local optimality conditions of first order.
\begin{theorem}
	\label{thm:OC_stationarity}
	Let $\bar u$ be locally optimal for \eqref{eq:OC} in the sense of $L^2(t_0,T;\R^N)$.
	Then
	for all $i = 1,\dots,N$,
	the function $\bar u_i$ (after modification on a set of measure zero)
	is continuous on the closed set $\set{\sgn(\bar u_i) = 1}$
	and the following conditions hold:  
	\begin{subequations}
		\label{eq:OC_stationarity}
		\begin{align}
			\label{eq:OC_stationarity:1}
			-(\nabla F(\bar u)(t))_i&\in \partial g_{i,+}(\bar u_i(t))&& \forall t\in \set{
			\sgn(\bar u_i) =1},\\
			\label{eq:OC_stationarity:2}
			-(\nabla F(\bar u)(t))_i &=
			\frac{ g_i( \bar u_i(t)) }{\bar u_i(t) }&&\forall t\in J_{\bar
			u_i}\cap (t_0,T),\\
			\label{eq:OC_stationarity:3}
			-(\nabla F(\bar u)(t_0))_i &\ge
			\frac{ g_i( \bar u_i(t_0)) }{\bar u_i(t_0) }&&\text{if}\ t_0 \in
			J_{\bar u_i}, \\
			\label{eq:OC_stationarity:4}
			-(\nabla F(\bar u)(T))_i &\ge
			\frac{ g_i( \bar u_i(T)) }{\bar u_i(T) }&&\text{if}\ T \in
			J_{\bar u_i}
			.
		\end{align}
	\end{subequations}
\end{theorem}
\begin{proof}
	We prove \eqref{eq:OC_stationarity:1} and
	\eqref{eq:OC_stationarity:2}--\eqref{eq:OC_stationarity:4} separately.

	We start by fixing $i \in \set{1,\ldots,N}$
	and vary only the component $\bar u_i$ of the solution.
	For this, we introduce the abbreviation
	\begin{equation*}
		F_i(u_i)
		\coloneqq
		F(\bar u_1, \ldots, \bar u_{i-1}, u_i, \bar u_{i+1}, \ldots, \bar u_N)
		.
	\end{equation*}
	First,
	we vary $\bar u_i$
	only on the set $\set{\sgn(\bar u_i) = 1}$.
	That is, we consider
	\begin{align*}
		\min_{u_i \in L^2(\set{\sgn(\bar u_i) = 1})}
		\quad
		F_i(\chi_{\set{\sgn(\bar u_i) = 1}} u_i)
		+
		\int_{\set{\sgn(\bar u_i) = 1}} g_{i,+}(u_i(t))\dt
		.
	\end{align*}
	Clearly, $\bar u_i$ is a local minimizer of this auxiliary problem.
	The convexity of $g_{i,+}$ implies that the integral
	in the objective is a convex function of $u_i$.
	Therefore, the necessary optimality condition directly follows from the standard theory of convex optimization, i.e.,
	we get
	\begin{equation*}
		- \left[\nabla F(\bar u)(t) \right]_i  \in \partial g_{i,+}(\bar u_i(t))
		\qquad \text{for a.a.\ }  t\in \set{\sgn(\bar u_i) = 1}.
	\end{equation*}
	By assumption, $g_{i,+}$ is strongly convex
	and this implies that $(\partial g_{i,+})^{-1} = \partial g_{i,+}\conjugate$
	is single-valued and Lipschitz continuous.
	Consequently,
	\begin{equation*}
		\bar u_i(t)
		=
		\partial g_{i,+}^{-1}\parens*{ - \left[\nabla F(\bar u)(t) \right]_i }
		\qquad \text{for a.a.\ }  t\in \set{\sgn(\bar u_i) = 1}.
	\end{equation*}
	Together with the
	assumed continuity of $\nabla F(\bar u) \colon [t_0, T] \to \R^N$,
	this implies continuity of the right-hand side w.r.t.\ $t$
	and, therefore, $\bar u_i$ (after modification on a set of measure zero)
	is continuous on
	the closed set $\set{\sgn(\bar u_i) = 1}$. 
	This verifies \eqref{eq:OC_stationarity:1}.

	For the other conditions we only vary $\sgn(\bar u_i)$.
	To this end, we denote by $\bar v_i \in C([t_0,T]; [a_i, b_i])$
	an arbitrary continuous extension of $\bar u_i |_{\set{\sgn(\bar u_i) = 1}}$.
	Further, for $\alpha\in L^1(t_0,T)$
	we introduce
	\begin{equation*}
		\tilde F(\alpha) \coloneqq
		F_i( \alpha \bar v_i)   + \int_{t_0}^{T}\alpha(t)g(\bar v_i(t))\dt
		.
	\end{equation*}
	Now, it is clear that $\bar\alpha = \sgn(\bar u_i)$
	is a local minimizer of
	\begin{equation*}
		\min_{\alpha \in \BV(t_0,T;\{0,1\})}
		\quad
		\tilde F(\alpha) + \TV(\alpha).
	\end{equation*}
	Since $g(\bar v_i(t))$ is continuous with respect to $t$ by
	\cref{asmp:OC}\ref{asmp:g} and $F$ already has the assumed regularity from
	\cref{asmp:OC}\ref{asmp:F}, the function $\tilde F$
	is Fréchet differentiable from $L^1(t_0,T)$ to $\R$
	and its derivative is given by
	\begin{equation*}
		\nabla \tilde F(\alpha)
		=
		\nabla F_i(\alpha \bar v_i) \bar v_i + g(\bar v_i)
		\in C([t_0,T])
		.
	\end{equation*}
	This regularity is enough to apply the optimality conditions from
	\cite[Lemma~4.10]{Leyffer2022}
	(see also the note after \cite[Theorem~3.6]{Marko2023}),
	to obtain
	\begin{equation*}
		\nabla \tilde F(\sgn(\bar u_i))(t)  = 0\qquad \forall t\in J_{\bar u_i}
		\cap (t_0, T)
	\end{equation*}
	for each point from the jumping set, except at
	the boundary
	(recall, that in these references, the usual definition of the total variation is used).
	The continuity of $\bar u_i$ on $\set{\sgn(\bar u_i) = 1}$ implies that 
	$\bar v_i(t) = \bar u_i(t) \ge a_i > 0$
	at all switching points $t \in J_{\bar u_i}$.
	Consequently, \eqref{eq:OC_stationarity:2} follows.

	If a jump occurs at a boundary point,
	then this switching point can only be moved in one direction
	and
	a simple argument yields
	\begin{equation*}
		\nabla \tilde F(\sgn(\bar u_i))(t_0) \le 0 \quad \text{if $t_0 \in J_{\bar u_i}$}
		\qquad\text{and}\qquad
		\nabla \tilde F(\sgn(\bar u_i))(T) \le 0 \quad \text{if $T \in J_{\bar u_i}$} .
	\end{equation*}
	Together with the above formula for $\nabla \tilde F (\sgn(\bar u_i))$,
	the stationarity conditions \eqref{eq:OC_stationarity} follow.
\end{proof} 
In order to shed some light on the connection of \eqref{eq:OC_stationarity:1} and \eqref{eq:OC_stationarity:2},
we provide a small lemma.
\begin{lemma}
	\label{lem:tangets}
	For all $i = 1, \ldots, N$,
	\begin{itemize}
		\item
			the function $v \mapsto g_{i}(v)/v$ has a unique minimizer $u^\star_i$ on $[a_i, b_i]$,
		\item
			the inclusion $g_{i}(v)/v \in \partial g_{i,+}(v)$ has a unique solution $v^\star_i \in [a_i, b_i]$,
	\end{itemize}
	and both coincide, i.e.,
	$u^\star_i = v^\star_i$.
\end{lemma}
\begin{proof}
	For simplicity, we drop the index $i$ in the proof.
	We have the chain of equivalences
	\begin{align*}
		v^\star \text{ minimizes } v \mapsto g(v)/v
		\qquad
		&\Longleftrightarrow\qquad
		\frac{g(v)}{v} \ge \frac{g(v^\star)}{v^\star}
		&& \forall v \in [a, b]
		\\
		&\Longleftrightarrow\qquad
		g(v) \ge g(v^\star) + \frac{g(v^\star)}{v^\star} ( v - v^\star )
		&& \forall v \in [a, b]
		\\
		&\Longleftrightarrow\qquad
		\frac{g(v^\star)}{v^\star} \in \partial g_+(v^\star)
		.
	\end{align*}
	Note that $g$ and $g_+$ coincide on $[a,b]$.
	This shows that minimizers of $v \mapsto g(v)/v$
	and solutions of $g(v)/v \in \partial g_+(v)$
	coincide.
	The strong convexity of $g$ implies
	\begin{equation*}
		g(v) > g(v^\star) + \frac{g(v^\star)}{v^\star} ( v - v^\star )
		\qquad\forall v \in [a, b] \setminus \set{v^\star}
	\end{equation*}
	and, consequently,
	$g(v)/v > g(v^\star)/v^\star$ for all $v \in [a, b] \setminus \set{v^\star}$.
	This provides the desired uniqueness.
\end{proof}
\begin{remark}
	Let us assume that $\bar u$ satisfies the conditions
	\eqref{eq:OC_stationarity:1}
	and
	\eqref{eq:OC_stationarity:2}.
	For a switching point $s \in J_{\bar u_i} \cap (t_0, T)$,
	this implies
	\begin{equation*}
		\frac{g_i(\bar u_i(s))}{\bar u_i(s)}
		\in
		\partial g_{i,+}(\bar u_i(s))
		.
	\end{equation*}
	Due to \cref{lem:tangets},
	the value $\bar u_i(s)$ is uniquely 
	determined
	and it equals
	\begin{equation*}
		u_i^\star
		\coloneqq
		\argmin_{v \in [a_i, b_i]} \frac{g_i(v)}{v}
		,
	\end{equation*}
	i.e.,
	$\bar u_i(s) = u_i^\star$.
	This is somewhat surprising. It characterizes
	$\bar u_i(s)$ as the
	point at which a tangent of $g_i$
	goes through the origin.
	This means
	that independently of the choice for $F$, switches can only occur from and to
	the predetermined value $u_i^\star$, solely depending on $g_i$.

	At the boundary switching points $s \in J_{\bar u_i} \cap \set{t_0, T}$,
	one can argue similarly.
	This shows that \eqref{eq:OC_stationarity:3} and \eqref{eq:OC_stationarity:4}
	imply $\bar u_i(s) \ge u_i^\star$.
	\label{rem:switching_point_optimality}
\end{remark}
We already know from \cref{sec:convex_problem} that a criticality measure
for a convex problem can be defined via the proximal point mapping.
The remaining part of the criticality measure is built from the additional
stationarity conditions at the switching points. For the sake of clarity
it is beneficial to introduce the criticality measure in two parts. 
For simplicity
of notation we introduce the following shorthand.
\begin{definition}
	Let $u\in \Uad$ be given.
	For $p\in [1,\infty)$ and $f\in L^p(t_0,T;\R^N)$ we define
	\begin{equation*}
		\norm{f}^p_{L^p(\set{\sgn(u) = 1})} \coloneqq
		\sum_{i=1}^{N}\int_{\{\sgn(u_i) = 1\}} \abs{f_i(t)}^p \dt 
		\in
		[0,\infty)
		.
	\end{equation*}
	Similarly, for $u, w \in \Uad$
	and
	$f\in L^p(t_0,T;\R^N)$ we define
	\begin{equation*}
		\norm{f}^p_{L^p(\set{\sgn(u w) = 1})} \coloneqq
		\sum_{i=1}^{N}\int_{\{\sgn(u_i) = 1\} \cap \{\sgn(w_i) = 1\}} \abs{f_i(t)}^p \dt 
		\in
		[0,\infty)
		.
	\end{equation*}
	\label{def:alpha_norm}
\end{definition}
\begin{definition}
	\label{def:OC_crit_one_switch}
	The criticality measure for a component $i$ of a feasible control $u\in\Uad$ at
	a single switching point $s\in J_{u_i}$ is defined via:
	\begin{equation*}
		\begin{aligned}
			v_{u,i}(s)
			&\coloneqq
			\argmin_{z \in [a_i, b_i]} [\nabla F(u)(s)]_i z + g_i(z)
			,
			\\
			V_{u,i}(s)
			&\coloneqq
			\min_{z \in [a_i, b_i]} [\nabla F(u)(s)]_i z + g_i(z)
			,
			\\
			C_{\text{switch}}(u,i,s) &\coloneqq
			\max\left\{V_{u,i}(s),- V_{u,i}(s)\frac{\min\{s - t_0,T-s\}
			}{T-t_0}\right\}
			.
		\end{aligned}
	\end{equation*}
	For a fixed value $r > 0$, we define the criticality measure for the active part of the
	control and combine both to obtain the actual criticality measure:
	\begin{align*}
		C_{\text{prox}}(u) &\coloneqq\frac{1}{2r^2} \norm{u - \prox_{r G_+}(u
		- r\nabla F(u))}^2_{L^2(\{\sgn(u) = 1\})}
		,
		\\[1em]
		C_{\text{switch}}(u) &\coloneqq \max_{i= 1,\ldots,N}\max_{s\in
		J_{u_i}} C_{\text{switch}}(u,i,s)
		,
		\\[1em]
		C(u) &\coloneqq \max\{C_{\text{prox}}(u),C_{\text{switch}}(u)\}.
	\end{align*}

\end{definition}
In case $N = 1$
(which will be assumed in some proofs for simplicity of the presentation),
we write
$v_u(s)$, $V_u(s)$ and $C_{\text{switch}}(u,s)$
instead of
$v_{u,1}(s)$, $V_{u,1}(s)$ and $C_{\text{switch}}(u,1,s)$,
respectively.

\begin{corollary}
	\label{cor:OC_crit}
	For all $u \in \Uad$, we have $C(u) \ge 0$.
	A feasible point $\bar u\in \Uad$ is stationary
	in the sense of \eqref{eq:OC_stationarity}
	if and only if
	$C(\bar u) = 0$.
\end{corollary}
\begin{proof}
	The term $\frac{\min\{s - t_0,T-s\}}{T-t_0}$
	appearing in $C_{\text{switch}}(u,i,s)$ is always non-negative,
	consequently,
	$C_{\text{switch}}(u,i,s) \ge 0$.
	Together with $C_{\text{prox}}(u) \ge 0$,
	we get $C(u) \ge 0$.

	From the definition of $G_+$,
	it follows that $\prox_{G_+}$ can be evaluated
	componentwise and pointwise.
	This implies that $C_{\text{prox}}(\bar u) = 0$
	is equivalent to \eqref{eq:OC_stationarity:1}.

	For a switching point
	$s \in J_{\bar u_i} \cap (t_0,T)$
	we get 
	\begin{align*}
		C_{\text{switch}}(u,i,s) = 0
		&\iff V_{u,i}(s) = 0 \\
		&\iff [\nabla F(\bar u)(s)]_i v_{u,i}(s) +  g_i(v_{u,i}(s)) = 0\\
		&\iff -[\nabla F(\bar u)(s)]_i = \frac{g_i(v_{u,i}(s))}{v_{u,i}(s)}
		.
	\end{align*}
	Thus, the condition
	$C_{\text{switch}}(u,i,s) = 0$ is equivalent
	to \eqref{eq:OC_stationarity:2} for all $s \in J_{\bar u_i} \cap (t_0, T)$.
	Similarly, for all $s \in J_{\bar u_i} \cap \{t_0,T\}$,
	$C_{\text{switch}}(u,i,s) = 0$
	if and only if $V_{u,i}(s)\le 0$ which is equivalent to
	\eqref{eq:OC_stationarity:3}, \eqref{eq:OC_stationarity:4}.
\end{proof}
\begin{remark}
	The asymmetric definition of $C_{\text{switch}}(u,i,s)$ may appear surprising
	at first, but is a consequence of including switches at the boundary points $t_0$ and
	$T$ for the computation of $\TV(\sgn(u))$. Scaling $-V_{u,i}(s)$ with
	the distance to the boundary prevents a discontinuity in the criticality measure, when
	moving the first switch to $t_0$ or the last switch to $T$,
	see the proof of \cref{lem:OC_stability} below.
	\label{rem:OC_crit}
\end{remark}

\begin{remark}
	\label{rem:Vui_intuition}
	We give a motivation for the definition of $V_{u,i}(s)$.
	For an interior switching point $s \in J_{u_i} \cap (t_0, T)$,
	the optimality condition $C(u) = 0$ from \cref{cor:OC_crit}
	implies $V_{u,i}(s) = 0$.
	If $V_{u,i}(s) \ne 0$,
	the sign of $V_{u,i}(s)$ indicates
	in which direction the switching point $s$ should be moved
	in order to improve the objective value.
	For another control $w \in \Uad$ in the neighborhood of $u$,
	a partial linearization of the objective leads to
	\begin{equation*}
		\int_{t_0}^T \nabla F(u)(t)^\top (w(t) - u(t)) + g(w(t)) - g(u(t)) \d t
		.
	\end{equation*}
	In case $V_{u,i}(s) > 0$,
	we have $[\nabla F(u)(s)]_i u_i(t) + g_i(u_i(t)) > 0$
	for all $t$ with $u_i(t) > 0$.
	By continuity,
	we get $[\nabla F(u)(t)]_i u_i(t) + g_i(u_i(t)) > 0$
	for all $t$ with $u_i(t) > 0$ close enough to $s$.
	Thus, we can improve the objective functional by switching off $u_i$
	in the neighborhood of the switching point $s$.
	In the other case $V_{u,i}(s) < 0$,
	we similarly get
	$[\nabla F(u)(t)]_i v_{u,i}(s) + g_i(v_{u,i}(s)) < 0 = [\nabla F(u)(t)]_i u_i(t) + g_i(u_i(t))$
	for all $t$ with $u_i(t) = 0$ close enough to $s$.
	Consequently,
	the objective can be improved by
	activating the control $u_i$ in this neighborhood of the switch $s$
	and setting it to the value $v_{u,i}(s)$.
	This will be made precise in the proof of \cref{lem:OC_sufficient_decrease}.
\end{remark}

\subsection{Application of the trust-region method}
\label{subsec:TR_application}

The goal of this section is to apply the abstract trust-region method from \cref{sec:ATRM}
to solve \eqref{eq:OC}.
We define the set $X$ from \cref{asmp:ATRM} to be the admissible set, i.e.,
\begin{equation*}
	X \coloneqq \Uad.
\end{equation*}
\cite{Leyffer2022,Manns2024} showed convergence of a  trust-region method for optimal control of ordinary
differential equations with integer-valued controls of bounded variation.
This would be applicable to \eqref{eq:OC} if we would have $a_i = b_i$.
In order to generalize to our situation $a_i < b_i$,
we
incorporate the proximal gradient method into the trust-region method for
switching point optimization and verify \cref{asmp:ATRM} to prove convergence
for the trust-region method applied to \eqref{eq:OC} in \cref{thm:OC_iterate_convergence}.
Now, we define the model function, the oracle and the prediction function.
\begin{definition}
	\label{def:OC_TRM}
	Consider a point $u\in \Uad$ and $\Delta > 0$.
	The model function $m_{u, \Delta} \colon \Uad \to \R$ is defined as 
	\begin{equation}
		\label{eq:OC_TRM}
		\begin{aligned}
			m_{u,\Delta}(w) &\coloneqq \du {\nabla F(u),w-u}\al_{L^2(t_0,T;\R^N)} + G(w)
			- G(u) \\
			&\quad+ \TV(\sgn(w)) - \TV(\sgn(u))
			+ \frac{1}{2\Delta}\norm{w-u}^2_{L^2(\set{\sgn(uw) = 1})}
			.
		\end{aligned}
	\end{equation}
	The oracle function $T \colon \Uad \times \R_+ \to \Uad$
	is defined
	by minimizing the model function over a trust region,
	i.e.,
	\begin{equation*}
		T(u,\Delta)
		\in
		\argmin\set*{
			m_{u,\Delta}(w)
			\given
			w\in \Uad,
			\norm{\sgn(u)-\sgn(w)}_{L^1(t_0,T;\R^N)}\le \Delta
		}
		.
	\end{equation*}
	The prediction function is defined as
	\begin{equation}
		\label{eq:OC_pred}
		\pred(u,\Delta)
		\coloneqq
		m_{u,\Delta}(u) - m_{u,\Delta}(T(u,\Delta))
		=
		- m_{u,\Delta}(T(u,\Delta))
		.
	\end{equation}
\end{definition}
Conceptually, the oracle function will select an optimal switching pattern $\sgn(w)$ for
the model problem inside the trust region, while
a proximal gradient step
with step length $\Delta$ is performed on the set $\set{\sgn(uw) = 1}$.
On the newly active set $\set{\sgn(w) = 1} \setminus \set{\sgn(u) = 1}$,
an optimal continuation is chosen.
Deriving an explicit formula for the oracle function will be done in \cref{subsec:solution_subproblem}.

Since the model function $m_{u,\Delta}$ is not convex,
the minimizer might not be unique
and
for the oracle, we can use any of its minimizers.
Note that the existence of a minimizer of the model function inside the trust region
can be shown along the lines of the proof of \cref{thm:existence} by using the following lemma.
\begin{lemma}
	\label{lem:something_is_lsc}
	Let $u \in \Uad$ and a sequence $(w_k) \subset \Uad$
	be given such that
	$w_k \weakly \bar w$ in $L^2(t_0,T; \R^N)$
	and
	$\sgn(w_k) \to \sgn(\bar w)$ in $L^1(t_0,T; \R^N)$.
	Then,
	\begin{equation*}
		\norm{ u - \bar w }_{L^2(\set{\sgn(u \bar w) = 1})}^2
		\le
		\liminf_{k \to \infty}
		\norm{ u - w_k }_{L^2(\set{\sgn(u w_k) = 1})}^2
		.
	\end{equation*}
\end{lemma}
\begin{proof}
	For simplicity, we argue in the case $N = 1$.
	From the weak lower semicontinuity of the norm in $L^2$, we get
	\begin{equation*}
		\norm{ u - \bar w }_{L^2(\set{\sgn(u \bar w) = 1})}^2
		\le
		\liminf_{k \to \infty}
		\norm{ u - w_k }_{L^2(\set{\sgn(u \bar w) = 1})}^2
		.
	\end{equation*}
	Further,
	\begin{align*}
		\MoveEqLeft
		\abs*{
			\norm{ u - w_k }_{L^2(\set{\sgn(u \bar w) = 1})}^2
			-
			\norm{ u - w_k }_{L^2(\set{\sgn(u w_k) = 1})}^2
		}
		\\&
		=
		\norm{ u - w_k }_{L^2(\set{\sgn(u \bar w) = 1} \mathbin{\triangle} \set{\sgn(u w_k) = 1})}^2
		\\&
		\le
		\norm{ u - w_k }_{L^2(\set{\sgn(\bar w) = 1} \mathbin{\triangle} \set{\sgn(w_k) = 1})}^2
		\\&
		\le
		\norm{ u - w_k }_{L^\infty(t_0,T)}^2
		\abs{\set{\sgn(\bar w) = 1} \mathbin{\triangle} \set{\sgn(w_k) = 1}}
		\\&
		=
		\norm{ u - w_k }_{L^\infty(t_0,T)}^2
		\norm{\sgn(\bar w) - \sgn(w_k)}_{L^1(t_0,T)}
		\to
		0
		,
	\end{align*}
	where $\triangle$ denotes the symmetric difference.
	In the last equality, we used that $\sgn(\bar w)$ and $\sgn(w_k)$ only take values in $\set{0,1}$.
	This shows the claim.
\end{proof}

Now, we have defined the ingredients
for applying the abstract \cref{alg:ATRM}
to \eqref{eq:OC}.
It remains to verify \cref{asmp:ATRM}.
\begin{lemma}[Monotonicity]
	\label{lem:OC_monotonicity}
	The prediction function \eqref{eq:OC_pred} is nondecreasing,
	i.e., \cref{asmp:ATRM}\eqref{eq:monotonicity} holds.
\end{lemma}
\begin{proof}
	As $\Delta$ increases, the model function $m_{u, \Delta}$ decreases
	and the trust region becomes larger.
	This yields the claim.
\end{proof}
For the verification of the remaining parts of \cref{asmp:ATRM},
we combine the arguments from \cite[Section~5]{Manns2024} with the results from \cref{subsec:convergence_prox}.
This means that we have to define $\Dmax_a(u)$ accordingly.
\begin{definition}
	\label{def:minimal_gap}
	Given $u\in \Uad$ and the minimal representation $\hat t_u$ for
	$\sgn(u)$ from \cref{prop:minimal_representation}, we define $\Dmax_a(u)\in (0,r]$ via the smallest gap
	between two sign changes
	\begin{equation*}
		\Dmax_a(u) \coloneqq
		\min\set*{ r,
			\min_{i = 1,\ldots,N} \min_{j = 1,\ldots,2K_i
			- 1} ((\hat t_{u_i})_{j+1} - (\hat t_{u_i})_{j})
		}
		,
	\end{equation*}
	where $K_i \coloneqq \TV(\sgn(u_i))/2$ and $r$ is as in the definition of the criticality measure,
	see \cref{def:OC_crit_one_switch}.
\end{definition}
\begin{lemma}[Sufficient decrease]
	\label{lem:OC_sufficient_decrease}
	The prediction function \eqref{eq:OC_pred} satisfies
	\cref{asmp:ATRM}\eqref{eq:sufficient_decrease} with $\Dmax_a(u)$ as defined
	in  \cref{def:minimal_gap}, $f\equiv\Id$,
	$c = \frac{1}{2}L_{\nabla F} b$ and $s = 2$.
\end{lemma}
\begin{proof}
	Let $u \in \Uad$ and $\Delta \in (0, \Dmax_a(u))$ be given.
	Since the predicted decrease $\pred(u, \Delta)$ is defined via the solution 
	of the model function $m_{u,\Delta}$,
	we obtain the estimate
	\begin{equation*}
		\pred(u, \Delta)
		\ge
		-m_{u,\Delta}(w)
	\end{equation*}
	for all $w \in \Uad$ with $\norm{\sgn(u) - \sgn(w)}_{L^1(t_0,T;\R^N)} \le \Delta$.
	We will construct several feasible points $w$ for comparison.
	First, we consider
	\begin{equation*}
		\wprx\coloneqq \prox_{\Delta G_+}(u - \Delta \nabla F(u))
	\end{equation*}
	and set $w_i \coloneqq w_{\prox,i} \rchi_{\set{\sgn(u_i)=1}}$, $i = 1,\ldots,N$.
	It is clear that $\sgn(w) = \sgn(u)$.
	Consequently,
	this case is completely analogous to \cref{lem:OC_conv_predictability}.
	Since $\Delta < \Dmax_a(u) \le r$,
	where $r > 0$ was used in \cref{def:OC_crit_one_switch} for the definition
	of the criticality measure,
	this lemma implies
	\begin{equation}
		\label{eq:OC_predictability_part_1}
		\pred(u,\Delta)
		\ge
		-m_{u,\Delta}(w)
		\ge \frac{1}{2r^2} \norm{u - \prox_{rG}(u - r\nabla
		F(u))}^2_{L^2(\set{\sgn(u) = 1})} \Delta
		= C_{\text{prox}}(u) \Delta
		.
	\end{equation}

	Next, we consider a switching point $s \in J_{u_i} \cap [t_0, T)$.
	We construct $w$ by
	shifting this switch by
	at most $\Delta < \Dmax_a(u)$.
	Since the distance between any  two switches is at least
	$\Dmax_a(u)$, we have $\TV(\sgn(u)) = \TV(\sgn(w))$.
	W.l.o.g., we consider the case $N=1$
	and, further, $\sgn(u(s-)) = 0$ and $\sgn(u(s+)) = 1$,
	i.e., the control is switched on at $s$.
	The other cases can be handled similarly.
	From \cref{asmp:OC}\ref{asmp:F}
	we get that $\nabla F(u)$ is Lipschitz continuous
	with constant $L_{\nabla F(u)} \le L_{\nabla F}$
	and this implies
	\begin{align}
		\nonumber
		\abs*{ \int_{s}^{s+h}(\nabla F(u)(s) - \nabla F(u)(t)) w(t) \d t }
		&\le
		\norm{w}_{L^\infty(s, s+h)}
		\int_{s}^{s+h} \abs{ \nabla F(u)(s) - \nabla F(u)(t)} \d t
		\\&
		\label{eq:grad_F_pointwise}
		\le
		L_{\nabla F(u)}
		b
		\int_s^{s+h} \abs{s - t} \d t
		\le
		\frac12 h^2 L_{\nabla F} b
	\end{align}
	for feasible $w \in \Uad$
	and $h > 0$ small enough such that $s + h \in (t_0,T)$.

	We make a distinction by cases with respect to the sign of
	$V_u(s)$ from \cref{def:OC_crit_one_switch}. As mentioned in
	\cref{rem:Vui_intuition}, the sign of $V_{u}(s)$ encodes
	in which direction the switching point should be moved.

	In case $V_u(s) > 0$,
	we first note that
	$s < T - \Delta$, since the control must be switched off at least once in $(s, T]$
	and any two switching points are more than $\Delta$ apart.
	Thus, we can define
	$w \coloneqq
	u - u\rchi_{[s,s+\Delta]}$, i.e., the control is switched on a little bit later.
	We immediately get $\TV(\sgn(w)) = \TV(\sgn(u))$
	and
	$\norm{w-u}^2_{\set{\sgn(uw) = 1}} = 0$.
	Together with $w = 0$ on $(s, s+\Delta)$,
	we have 
	\begin{align*}
		\pred(u,\Delta) &\ge - m_{u,\Delta}(w)
		= \du {-\nabla F(u),w-u}\al_{L^2(t_0,T)} + G(u) - G(w)
		\\
		&=\int_{s}^{s+\Delta}-\nabla F(u)(t)(w(t)-u(t))\dt
		+ \int_{s}^{s+\Delta} g(u(t)) - g(w(t))\dt \\
		&\stackrel{\mathclap{\eqref{eq:grad_F_pointwise}}}{\ge}\int_{s}^{s+\Delta}\nabla F(u)(s)
		u(t) + g(u(t))\dt - \frac 12 \Delta^2 L_{\nabla
		F}b
		\\
		&\ge \Delta V_u(s) - \frac 12 \Delta^2 L_{\nabla F} b
		= \Delta C_{\text{switch}}(u,1,s) - \frac 12 \Delta^2 L_{\nabla F} b
		,
	\end{align*}
	where the last estimate follows from the definition of $V_u(s)$
	and the last equality is the definition of $C_{\text{switch}}(u,1,s)$.

	In case $V_u(s) < 0$, we have
	\begin{equation*}
		0
		<
		C_{\text{switch}}(u, 1, s)
		=
		- V_{u}(s)\frac{\min\{s - t_0,T-s\}
		}{T-t_0}
	\end{equation*}
	and this implies $s > t_0$.
	Now, we want to move the switching point to the left.
	If $s$ is near the initial time $t_0$, this might limit our possibility to move the switch.
	Therefore, we set $h \coloneqq \min\set{\Delta, s - t_0}$.
	We use $v_{u}(s)$ from \cref{def:OC_crit_one_switch} (i.e., the minimizer of $w \mapsto [\nabla
	F(u)(s)]w + g(w)$) 
	to define $w \coloneqq u + v_u(s)\rchi_{[s-h,s]}$,
	i.e.,
	the control is switched on slightly earlier with the fixed value $v_u(s)$.
	We proceed similar to
	the above estimate.
	Using $u = 0$ on $(s-h, s)$,
	we get
	\begin{align*}
		\pred(u,\Delta) &\ge\int_{s-h}^{s}-\nabla
		F(u)(t)(w(t)-u(t)) \dt + 
		\int_{s-h}^{s}- g(w(t))+ g(u(t)) \dt\\
		&\stackrel{\mathclap{\eqref{eq:grad_F_pointwise}}}{\ge}\int_{s-h}^{s}-\nabla
		F(u)(s)v_u(s)
		- g(v_u(s))\dt - \frac 12 h^2 L_{\nabla F} \norm{v_u(s)}_{L^\infty(s-h,s)}\\
		&= h(-V_u(s))- \frac 12 h^2 L_{\nabla F} b
		.
	\end{align*}
	We continue with
	\begin{align*}
		h(-V_u(s))
		&=
		-V_u(s) \min\set{\Delta, s - t_0}
		=
		-V_u(s) \frac{(s - t_0) \Delta}{\max\set{s - t_0, \Delta}}
		\\&
		\ge
		-V_u(s) \frac{s - t_0}{T - t_0} \Delta
		\ge
		\Delta C_{\text{switch}}(u,1,s),
	\end{align*}
	where the first estimate uses $\Delta \le T - t_0$, $s - t_0 \le T - t_0$,
	and the second estimate follows from the definition of $C_{\text{switch}}(u,1,s)$.
	Combining the last two equations,
	we again get
	\begin{equation*}
		\pred(u, \Delta)
		\ge
		\Delta C_{\text{switch}}(u,1,s) - \frac 12 h^2 L_{\nabla F} b.
	\end{equation*}
	Finally, if $V_u(s) = 0$, the same inequality holds trivially,
	since $\pred(u, \Delta) \ge 0$ and $C_{\text{switch}}(u,1,s) = 0$.

	Since the criticality measure is defined as the maximum of the proximal part
	and the term for the
	switching points,
	we obtain
	\begin{align*}
		\pred(u,\Delta)&\ge
		\max\left\{C_{\text{prox}}(u)\Delta,C_{\text{switch}}(u)\Delta
		-\frac12 L_{\nabla F} \Delta^2 b\right\}\\
		&\ge C(u)\Delta
		- \frac12 L_{\nabla F} \Delta^2 b
		\qquad\forall \Delta \in (0, \Dmax_a(u))
		.
	\end{align*}
\end{proof}

\begin{remark}
	The above proof reveals some interesting properties. First, we can guarantee the
	prediction estimate by analyzing the position of the switches and the improvement by
	the proximal gradient method on $\{\sgn(u) = 1\}$ independently. Further a switching point is either stationary in the sense of \eqref{eq:OC_stationarity:2}--\eqref{eq:OC_stationarity:3}
	or can always be moved to guarantee an improvement. While the proof
	is constructive, the constructed improvement always retains the switching behavior.
	However, a minimizer $T(u,\Delta)$ of the trust-region subproblem (\cref{def:OC_TRM})
	is not forced to follow the same switching pattern as $u$.
	Consequently, the trust-region algorithm may produce solutions
	with better function values compared to a pure
	switching-point optimization approach.
	\label{rem:OC_predictability}
\end{remark}
In the next lemma, we check that the oracle stays in a scaled trust region,
similar to \cref{lem:OC_conv_iterate_in_trust_region}.
\begin{lemma}
	\label{lem:OC_iterate_in_trust_region}
	For all $u \in \Uad$, we have
	\begin{equation*}
		\norm{T(u, \Delta) - u}_{L^1(t_0, T; \R^N)}
		\le
		\rho \Delta
		\qquad
		\forall \Delta \in (0,\infty),
	\end{equation*}
	with $\rho \coloneqq 2 (C_{\nabla F} + L_g) N(T - t_0) + b$.
\end{lemma}
\begin{proof}
	We set $\hat u = T(u, \Delta)$.
	We use the fact
	that the norm can be split into the part where both controls
	are active and a remainder.
	For $i \in \set{1, \ldots, N}$ and $t \in \set{\sgn(u_i \hat u_i) = 1}$,
	we have
	\begin{equation*}
		\hat u_i(t)
		=
		\prox_{\Delta g_i^+}(
			u_i(t)
			-
			\Delta [\nabla F(u)(t)]_i
		)
		.
	\end{equation*}
	Now,
	we argue as in \eqref{eq:estimate_prox}
	and obtain
	\begin{equation*}
		\abs{
			\hat u_i(t) - u_i(t)
		}
		\le
		2 \parens*{ \abs{[\nabla F(u)(t)]_i} + L_g } \Delta
		\le
		2 \parens*{ C_{\nabla F} + L_g } \Delta
		\qquad
		\forall t \in \set{\sgn(u_i \hat u_i) = 1}
		,
	\end{equation*}
	see \eqref{eq:bound_nabla_F}.
	On the other hand, if $t \not\in \set{\sgn(u_i \hat u_i) = 1}$,
	at least one of the controls is zero.
	Unless both are zero, we have
	$\abs{\sgn(\hat u_i(t)) - \sgn(u_i(t))} = 1$.
	This implies
	\begin{equation*}
		\abs{
			\hat u_i(t) - u_i(t)
		}
		\le
		b \abs{\sgn(\hat u_i(t)) - \sgn(u_i(t))}
		\qquad
		\forall t \not\in \set{\sgn(u_i \hat u_i) = 1}
		.
	\end{equation*}
	By integrating over these inequalities, we get
	\begin{align*}
		\norm{\hat u - u}_{L^1(t_0, T; \R^N)}
		&=
		\int_{t_0}^T
		\sum_{i = 1}^N
		\abs{
			\hat u_i(t) - u_i(t)
		}
		\d t
		\\&
		\le
		2 \parens*{ C_{\nabla F} + L_g } \Delta N (T - t_0)
		+
		b \norm{\sgn(\hat u) - \sgn(u)}_{L^1(t_0,T;\R^N)}
		\\&
		\le
		2 \parens*{ C_{\nabla F} + L_g } \Delta N (T - t_0)
		+
		b \Delta
		,
	\end{align*}
	where we used the trust-region constraint from the definition of the oracle,
	see \cref{def:OC_TRM}.
	This shows the claim.
\end{proof}
\begin{lemma}[Accuracy]
	For an arbitrary $\eta \in (0,1]$,
	\cref{asmp:ATRM}\eqref{eq:accuracy} and \eqref{eq:accuracy_E} 
	are satisfied with
	$C, r$ as in \cref{def:OC_crit_one_switch},
	$\Dmax_a$ from \cref{def:minimal_gap}, $\rho$ from
	\cref{lem:OC_iterate_in_trust_region}, $\delta = 1$, $s = 2$,
	\begin{equation*}
		c\coloneqq
		\frac{1}{2} L_{\nabla F} \parens*{ b +  \frac{\rho^2}{1 - \eta}}
		,\qquad\text{and}\qquad
		\Dmax_b \coloneqq \min\set*{\frac{r}{2}, \frac{1}{(C_{\nabla
			F}+L_g)b},
			\sqrt{
					\frac{2 (1 - \eta)}{L_{\nabla F} \rho^2}
			}
		}
		.
	\end{equation*}
	\label{lem:OC_accuracy}
\end{lemma}
\begin{proof}
	The proof can be carried out similarly to \cref{lem:OC_conv_accuracy}, but we
	have to be careful with the case of switching points being too close,
	i.e., $\Dmax_a(u) \le \Delta$.
	This leads to the corresponding exception \eqref{eq:accuracy_E},
	which will be verified analogously to \cite[Lemma~5.4]{Manns2024}.

	We start with an arbitrary $\Delta \in (0, r)$.
	First, we define $\hat u \coloneqq T(u,\Delta)$ and $\ared(u,\Delta)$ will be estimated using
	the first-order approximation
	\begin{equation*}
		F(u) - F(\hat u)
		\ge
		\dual{\nabla F(u)}{u - \hat u}_{L^2} - \frac12 L_{\nabla F} \norm{u - \hat u}_{L^1(0,T; \R^N)}^2
	\end{equation*}
	of $F$,
	which follows from the Lipschitz continuity of $\nabla F \colon L^1(t_0,T;\R^N) \to L^\infty(t_0,T;\R^N)$.
	Using this estimate yields
	\begin{align*}
		\ared(u, \Delta) &= J(u) - J(\hat u)\\
		&=F(u) + G(u) + \TV(\sgn(u)) - F(\hat u) - G(\hat u)
		- \TV(\sgn(\hat u))\\
		&\ge\bigg(- \du {\nabla F(u),\hat u - u}\al_{L^2} +G(u)
			- G(\hat u) \\
			&\qquad +(\TV(\sgn(u)) - \TV(\sgn(\hat u)))
		-\frac{1}{2\Delta}\norm{\hat u-u}^2_{\sgn(u \hat u) = 1}\bigg) \\
		&\quad
		- \frac12 L_{\nabla F} \norm{u - \hat u}_{L^1(0,T; \R^N)}^2
		+ \frac{1}{2\Delta}\norm{\hat u- u}^2_{L^2(\set{\sgn(u\hat u) =1})}\\ 
		&\ge \pred(u,\Delta)
		- \frac12 L_{\nabla F} \norm{u - \hat u}_{L^1(0,T; \R^N)}^2
		.
	\end{align*}
	Combining this with
	\cref{lem:OC_iterate_in_trust_region},
	we get
	\begin{equation}
		\label{eq:ared_sub_eta_pred}
		\ared(u,\Delta) - \eta\,\pred(u,\Delta)
		\ge
		(1-\eta) \pred(u, \Delta)
		-
		\frac12 L_{\nabla F} 
		\rho^2 \Delta^2
		\qquad\forall \Delta \in (0, r)
		.
	\end{equation}
	In case $\Delta \in (0, \Dmax_a(u))$, we can insert the predicted reduction from \cref{lem:OC_sufficient_decrease} and obtain
	\begin{equation*}
		\ared(u,\Delta) - \eta\,\pred(u,\Delta) \ge (1-\eta) \left(C(u) \Delta
		-  c\Delta^2\right)
	\end{equation*}
	with the constant $c$ as defined in the claim, proving \cref{asmp:ATRM}\eqref{eq:accuracy}
	for $\Delta \in (0,\Dmax_a(u))$.

	It remains to check \cref{asmp:ATRM}\eqref{eq:accuracy_E},
	i.e., the case $\Delta \in [\Dmax_a(u), \Dmax_b]$.
	Since $\Dmax_b < r$, this implies
	$\Dmax_a(u) < r$ and, consequently,
	there exist indices $i \in \set{1,\ldots,N}$
	and $j \in \set{1, \ldots, 2 K_i - 1}$
	such that
	$(\hat t_{u_i} )_{j+1} - (\hat t_{u_i} )_{j} = \Dmax_a(u) \le \Delta$.
	This means that
	in the component $u_i$
	we have a short time distance
	between two consecutive switches
	(either on-off-on or off-on-off).
	We will see that this situation can be improved
	by simply removing both switches.
	For brevity, we set $N=1$ and only consider the on-off-on situation.
	We further set
	$s \coloneqq (\hat t_{u} )_{j}$
	and
	$h \coloneqq \Dmax_a(u) = (\hat t_{u} )_{j+1} - (\hat t_{u} )_{j}$.
	We define $\hat u \in \Uad$
	via
	$\hat u = u + a \chi_{(s, s + h)}$.
	Due to
	$\norm{\sgn(\hat u) - \sgn(u)}_{L^1(t_0,T)} = h \le \Delta$,
	$\hat u$ belongs to the trust region.
	One has $\TV(\sgn(\hat u)) = \TV(\sgn( u)) - 2$ and the model function \eqref{eq:OC_TRM} can be
	estimated via
	\begin{equation*}
		m_{u,\Delta}(\hat u) \le \int_{s}^{s+h} \nabla F(u)(t)(\hat u(t) -  u(t))\dt + \int_{s}^{s+h}
		g(t,\hat u(t)) - g(t,u(t))\dt - 2
		.
	\end{equation*}
	With the global Lipschitz constant
	$L_g$ of $g$ and the Lipschitz constant $C_{\nabla F}$ of $F$
	from \cref{rem:F_Lipschitz}, one has
	\begin{equation*}
		\pred(u,\Delta)
		\ge
		-m_{u, \Delta}(\hat u)
		\ge 2 - C_{\nabla F}\norm{\hat u - u}_{L^1(t_0,T;\R^N)}
		- L_g\norm{\hat u-u}_{L^1(t_0,T;\R^N)}.
	\end{equation*}
	The trust-region constraint for \eqref{eq:OC_TRM} yields $\norm{\hat u - u}_{L^1(t_0,T)}\le b \Delta$
	and this implies
	\begin{equation*}
		\pred(u,\Delta) \ge 2 - (C_{\nabla F} + L_g) b\Delta \ge 1\qquad \forall\Delta
		\in [\Dmax_a(u), \Dmax_b]
		.
	\end{equation*}
	If we use this estimate for the predicted reduction we can revisit \eqref{eq:ared_sub_eta_pred}
	and obtain
	\begin{equation*}
		\ared(u,\Delta) - \eta\,\pred(u,\Delta)
		\ge
		(1-\eta)
		-
		\frac12 L_{\nabla F} \rho^2 \Delta^2
		\ge
		0
		\qquad
		\forall\Delta\le \Dmax_{b}.
	\end{equation*}
	This verifies
	\cref{asmp:ATRM}\eqref{eq:accuracy_E}. 
\end{proof}
\begin{lemma}[Stability]
	\label{lem:OC_stability}
	For the criticality measure $C(u)$ defined in \cref{def:OC_crit_one_switch} and the corresponding prediction function
	\eqref{eq:OC_pred}, there exist constants  $\Dmax_c > 0$ and $L > 0$, such that
	\begin{equation*}
		\abs{C(T(u,\Delta)) - C(u)}\le L \Delta
		\quad\lor\quad
		\pred(u,\Delta) \ge \frac{1}{2}
		\qquad
		\forall \Delta\le \Dmax_c,
	\end{equation*}
	i.e.,
	\cref{asmp:ATRM}\eqref{eq:stability} or \eqref{eq:stability_E}
	hold for all $\Delta \le \Dmax_c$ with $\Dmax_c \coloneqq \min\set{((C_{\nabla
	F} + L_g)\rho)^{-1},\Dmax_b}$ for $\rho$ as defined in
	\cref{lem:OC_iterate_in_trust_region}.
\end{lemma}
\begin{proof}
	Define $\hat u = T(u,\Delta)$. The proof is done by estimating the
	individual parts of the criticality measure
	and we start with the criticality measure for the
	proximal gradient step $C_{\prox}$.
	We partition it into the part where $u$ and $\hat u$ are active and
	a remainder term, which is bounded using the trust-region radius:
	\begin{align*}
		C_{\prox,u}(u) &\coloneqq \frac{1}{2r^2}\norm{u - \prox_{rG_+}(u-r\nabla
		F(u))}^2_{L^2(\set{\sgn(u\hat u)=1})}\\
		C_{\prox}(u) &= C_{\prox,u}(u) +     \frac{1}{2r^2}\norm{u - \prox_{rG_+}(u-r\nabla
		F(u))}^2_{L^2(\set{\sgn(u) = 1 \land \sgn(\hat u)=0})} \\
		&\le C_{\prox,u}(u) +  \frac{(b-a)^2}{2r^2}\Delta
		.
	\end{align*} 
	Note that $\hat u = \prox_{\Delta G_+}( u - \Delta \nabla F(u))$ on $\set{\sgn(u \hat u) = 1}$.
	Hence, we can directly use
	\cref{lem:OC_conv_stability} for $\abs{C_{\text{prox},u} - C_{\text{prox},\hat
	u}} \le L_C \Delta$ to recover Lipschitz continuity for the proximal gradient
	part, i.e.,
	\begin{equation*}
		\abs{C_{\text{prox}}(u) - C_{\text{prox}}(\hat u)} \le \left(L_C
		+ \frac{(b-a)^2}{2r^2}\right)\Delta
		.
	\end{equation*}

	Regarding $C_{\text{switch}}$, an apparent issue is that the position of switching points of $\hat u$ in relation to
	those of $u$ is not known at first. W.l.o.g., we limit the discussion to the
	case $N=1$ for brevity of the notation. We will argue similar as in \cite[Lemma
	5.5]{Manns2024} to relate the switching points of $u$ and $\hat u$.
	We define  $M\coloneqq (C_{\nabla F} + L_g) \rho $.
	Then,
	\begin{align*}
		\abs*{
			\pred(u,\Delta) - (\TV(\sgn (u)) - \TV(\sgn(\hat u)))
		}
		&=
		\du {\nabla F(u),u-\hat
		u}\al_{L^2(t_0,T)} + G(u) - G(\hat u)
		\\&
		\le
		M\Delta
		.
	\end{align*}
	Now the choice $\Dmax_c \coloneqq \min\{M^{-1},\Dmax_b\}$ guarantees, that
	for $\Delta \le \Dmax_c$ and $\TV(\sgn(\hat u)) < \TV(\sgn(u))$ the predicted reduction
	admits the lower bound:
	\begin{equation*}
		\pred(u,\Delta) \ge 1,
	\end{equation*}
	where we used that $\TV(\sgn(u))$ can only attain even values.
	Note that
	$\TV(\sgn(\hat u)) > \TV(\sgn(u))$ cannot occur for $\Delta \le \Dmax_c$, since the prediction is always
	positive.
	By the second
	part of \cref{lem:OC_accuracy}, this prediction is also guaranteed, if
	$\Delta\in[\Dmax_a(u),\Dmax_b]$.

	This leaves the case $\TV(\sgn(\hat u)) = \TV(\sgn(u))$ and $\Delta < \Dmax_a(u)$.
	We abbreviate $t \coloneqq \hat t_u$ and $s \coloneqq \hat t_{\hat u}$,
	which both belong to $\R^{2 K}$ with $K \coloneqq \TV(\sgn(u)) / 2$.
	Since the trust-region radius directly limits the measure
	of the set $\set{\sgn(u) \neq \sgn(\hat u)}$, for all $j = 1,\dots, 2K-1$ there exists a $\xi_j$, such that
	\begin{equation*}
		\xi_j\in (t_j,t_{j+1})
		\quad\land\quad
		\sgn(\hat u(\xi_j)) = \sgn(u(t_j+))
		= \sgn(u(t_{j+1} - )),
	\end{equation*}
	since otherwise $\norm{\sgn(u) - \sgn(\hat u)}_{L^1(t_j,t_{j+1})}
	= t_{j+1} - t_j > \Dmax_a(u) > \Delta.$
	For completeness of the
	argument, it is beneficial to extend the control with zero past $[t_0,T]$
	and consider additionally $\xi_0 < t_0, T < \xi_T$, such that we can use the
	same argument as above to conclude $s_j\in
	(\xi_{j-1},\xi_{j})$ for $j = 1,\dots, 2K$. With the above arguments it becomes clear, that the trust-region radius enforces $s_j\in [t_j
	-\Delta,t_j+\Delta]$,
	i.e., $u$ and $\hat u$ possess the same switching patterns
	and the switching points of $\hat u$ are (small) perturbations
	of the switching points of $u$.
	Consequently, the signs of $u$ and $\hat u$ 
	only differ in the neighborhoods of the switching points
	and we get
	\begin{equation*}
		\Delta \ge \norm{\sgn(u) - \sgn(\hat u)}_{L^1(t_0,T)} = \sum_{j=1}^{2K} \abs{t_j - s_j}
		.
	\end{equation*}

	Before we can move to the final estimate, we need to investigate the parts of
	$C_{\text{switch}}(u)$ with respect to Lipschitz continuity. Therefore, we recall the assumed Lipschitz continuities of $\nabla F$:
	\begin{align*}
		\norm{\nabla F(u) - \nabla F(\hat u)}_{L^\infty(t_0,T)} &\le L_{\nabla F} \norm{u-\hat
		u}_{L^1(t_0, T)},\\
		\abs{\nabla F(u)(s) - \nabla F(u)(t)} &\le
		L_{\nabla F(u)}\abs{s-t}
		\le L_{\nabla F}\abs{s-t}
		.
	\end{align*}
	Together with \cref{lem:OC_iterate_in_trust_region}, this implies
	\begin{align*}
		\abs{\nabla F(u)(t_j)- \nabla F(\hat u)(s_j)}
		&\le
		\abs{
			\nabla F(u)(t_j)
			-
			\nabla F(\hat u)(t_j)
		}
		+
		\abs{
			\nabla F(\hat u)(t_j)
			-
			\nabla F(\hat u)(s_j)
		}
		\\&
		\le
		L_{\nabla F} ( \norm{u - \hat u}_{L^1(t_0,T)} + \abs{t_j - s_j} )
		\le
		L_{\nabla F} ( \rho \Delta + \Delta)
		.
	\end{align*}
	Using the notation of \cref{def:OC_crit_one_switch},
	we get
	\begin{align*}
		V_u(t_j)
		&=
		\nabla F(u)(t_j) v_u(t_j) + g(v_u(t_j))
		\\&
		\le
		\nabla F(u)(t_j) v_{\hat u}(s_j) + g(v_{\hat u}(s_j))\\
		&=
		(\nabla F(u)(t_j) - \nabla F(\hat u)(s_j) ) v_{\hat u}(s_j) + V_{\hat u}(s_j)
		\\
		&
		\le
		L_{\nabla F} (\rho + 1) \Delta b +  V_{\hat u}(s_j)
	\end{align*}
	and vice cersa. Consequently,
	\begin{align*}
		\abs{
			V_u(t_j) - V_{\hat u}(s_j)
		}
	% &\le
	% \abs{\xi_{t,j} - \xi_{s,j}} b\\
	% &\le
	% L_{\nabla F} ( \norm{u - \hat u}_{L^1(t_0,T)} + \abs{t_j - s_j} ) b  \\
	% &
		\le L_{\nabla F} (\rho + 1) \Delta b
		.
	\end{align*}
	Using the definition of $C_{\text{switch}}$, this allows for the estimate
	\begin{align*}
		C_{\text{switch}}(\hat u,s_j) 
		&= \max\left\{V_{\hat u}(s_j),-
			V_{\hat u}(s_j)\frac{\min\{s_j - t_0,T-s_j\}
		}{T-t_0}\right\}\\
		&\le C_{\text{switch}}(u,t_j) +
		L_{\nabla F} (\rho + 1) \Delta b
	\end{align*}
	and vice versa switching the roles of $u$ and $\hat u$. From this point we can
	make the final estimate
	\begin{align*}
		\abs{  C_{\text{switch}}(u) - C_{\text{switch}}(\hat u)} &= 
		\abs*{
			\max_{j=1,\ldots,2K} 
			C_{\text{switch}}(u,t_j) - 
			\max_{j=1,\ldots,2K} 
		C_{\text{switch}}(\hat u,s_j)}\\
		&\le  
		\max_{j=1,\ldots,2K} 
		\abs*{ C_{\text{switch}}(u,t_j)   - C_{\text{switch}}(\hat u,s_j)}
		\\
		&\le
		L_{\nabla F} (\rho + 1) \Delta b
		.
	\end{align*}
	This shows that
	\cref{asmp:ATRM}\eqref{eq:stability} is satisfied.
\end{proof}
Now, we are in position to prove the convergence theorem.
\begin{theorem}
	\label{thm:OC_iterate_convergence}
	\cref{alg:ATRM} applied to \eqref{eq:OC} either terminates at a stationary
	point or produces a sequence of iterates for which $C(u_n)\rightarrow 0$.

	If in addition  
	$\nabla F \colon L^2(t_0, T; \R^N) \to L^2(t_0, T; \R^N)$ is a compact operator,
	$(u_n)$ contains a subsequence which converges strongly in $L^2(t_0, T; \R^N)$
	towards a stationary
	point $\bar u$ and $\sgn(u_n)\rightarrow \sgn(\bar u)$ in
	$L^1(t_0, T; \R^N)$.
\end{theorem}
\begin{proof}
	Convergence of the criticality measure immediately follows
	from \cref{thm:atrm_convergence}, since
	\cref{lem:OC_monotonicity,lem:OC_sufficient_decrease,lem:OC_accuracy,lem:OC_stability}
	verified the four parts of
	\cref{asmp:ATRM}.

	Similar to the these lemmas, the idea of the remaining part of the proof is to separate the
	arguments relating to the proximal gradient method and the part relating
	to the switching points. We will carry out the proof for $N=1$, i.e.
	$(u_k)\in L^\infty(t_0,T)$ for clarity of notation. It will become apparent from the individual
	arguments that they generalize to $N > 1$ without any changes.

	To this end we consider the sequence $(\sgn(u_k))$, which is bounded in
	$\BV(t_0,T;\{0,1\})$, since the objective includes the total variation
	regularization. By \cite[Theorem~3.23]{Ambrosio2000}, there exists a subsequence
	(without relabeling) that is weak-$\star$ convergent
	towards a limit point $\bar \alpha$. Further by \cite[Proposition
	3.12]{Ambrosio2000}, this subsequence $(\sgn(u_k))$ already convergences strongly in
	$L^1(t_0,T;\{0,1\})$.
	Thus, we can choose another
	subsequence (without relabeling) such that $\norm{\sgn(u_k)
	- \sgn(u_{k+1})}_{L^1(t_0,T)} \le 2^{-k}$.
	Here,
	$u_{k+1}$ is the element of the subsequence
	following $u_k$.

	For the same the subsequence, $(u_k)\subset \Uad \subset L^2(t_0,T)$ is bounded by
	the \cref{asmp:OC}\ref{asmp:g}, see also \eqref{def:notation}.
	Thus, there exists (without
	relabeling) a weakly convergent subsequence $u_k \rightharpoonup \bar u$.
	From \cref{lem:compatibility_limits}, we get
	$\sgn(\bar u)=\bar \alpha$ almost
	everywhere.

	Let $n \in \N$ be arbitrary.
	As a next step we want to split $\norm{u_k-\bar u}_{L^2(t_0,T)}$ into
	the part where the sign no longer changes, i.e.,
	which is supported on
	\begin{equation*}
		\mathcal{I}_n\coloneqq \bigcap_{l=n}^\infty \set{\sgn(u_l)=1}
		,
	\end{equation*}
	and a vanishing remainder.
	From a straightforward argument,
	we get for all indices $k \ge n$
	the inclusion
	\begin{align*}
		\MoveEqLeft
		\left( \set{\bar \alpha = 1} \cup \set{\sgn(u_k) = 1} \right)
		\setminus
		\mathcal{I}_n
		\subset
		\left( \bigcup_{l=n}^\infty \set{\sgn(u_l) = 1} \mathbin{\triangle}
		\set{\sgn(u_{l+1})=1} \right)
		.
	\end{align*}
	Since the above functions all have values in $\{0,1\}$ we get
	\begin{equation}
		\label{eq:TRM_strong_convergence_h3}
		\begin{aligned}
			\abs*{(\set{\bar \alpha = 1} \cup \set{\sgn(u_k)=1})\setminus
			\mathcal{I}_n} &\le\sum_{l=n}^\infty \norm{\sgn(u_l)
			- \sgn(u_{l+1})}_{L^1(t_0,T)}  \\   
			&\le \sum_{l=n}^\infty 2^{-l} \le 2^{-n+1}
			.
		\end{aligned}
	\end{equation}
	We split the domain into $\mathcal{I}_n$ and its complement,
	i.e.,
	for $k \ge n$ we get
	\begin{equation*}
		\begin{aligned}
			\norm{u_k - \bar u}^2_{L^2(t_0,T;\R^N)} &= \norm{u_k - \bar
			u}^2_{L^2(\mathcal{I}_n)}
			+\norm{u_k - \bar u}^2_{L^2(\set{\sgn(u_k)=1}\cup \set{\bar \alpha
		= 1})\setminus \mathcal{I}_n)}\\
		&\le \norm{u_k - \bar u}^2_{L^2(\mathcal{I}_n)} \\
		&\qquad + \norm{u_k - \bar
		u}_{L^\infty(t_0,T)}^2 \abs*{\set{\sgn(u_k)=1}\cup \set{\bar \alpha
= 1})\setminus \mathcal{I}_n)}
.
		\end{aligned}
	\end{equation*}
	Since the sequence $(u_k)$ only contains admissible points, we get
	$\norm{u_k - \bar u}^2_{L^\infty(t_0,T)}\le b^2$. With
	\eqref{eq:TRM_strong_convergence_h3} for any
	$\varepsilon \ge 0$ there exists a $n\in\N$, such that for all $k\ge n$
	we have
	\begin{equation}
		\label{eq:TRM_strong_convergence_h2}
		\norm{u_k - \bar u}_{L^2(t_0,T)}^2
		\le
		\norm{u_k - \bar u}_{L^2(\mathcal{I}_n)}^2 + \frac{\varepsilon}{2} 
		.
	\end{equation}
	In order to bound the first addend on the right-hand side,
	we utilize that
	$u_k$ is always positive on $\mathcal{I}_n$ for $k \ge n$
	and this means that
	the criticality measure $C(u_k)$
	contains information on
	a proximal gradient step on $\mathcal{I}_n$,
	see \cref{def:OC_crit_one_switch}.
	We define
	$\hat u_k \coloneqq u_k - \prox_{r G_+}(u_k - r\nabla F(u_k))$. From the
	convergence of the criticality measure we get  
	\begin{equation*}
		\norm{\hat u_k}_{L^2(\set{\sgn(u_k) = 1})}^2
		=
		\norm{u_k - \prox_{rG}(u_k - r\nabla F(u_k))}_{L^2(\set{\sgn(u_k) = 1})}^2
		=
		2r^2 C_{\prox}(u_k) \rightarrow
		0.
	\end{equation*}
	Consequently,
	\begin{equation}
		\label{eq:hat_u_conv}
		\norm{\hat u_k}_{L^2(\mathcal{I}_n)}^2
		\to
		0.
	\end{equation}
	Since $G_+$ is an integral functional,
	the proximal map can be evaluated pointwise.
	In particular, this yields
	\begin{equation*}
		- r\nabla F(u_k) + \hat u_k \in r\partial g_+(u_k - \hat u_k)\qquad \text{a.e. on }
		\mathcal{I}_n,
	\end{equation*}
	where $g_+$ is defined via
	\begin{equation*}
		g_+(v) \coloneqq
		\begin{cases}
			g(v)&\text{if } a \le v \le b \text{ a.e.\ on } (t_0, T), \\
			\infty&\text{else}.
		\end{cases}
	\end{equation*}
	Since $g_+$ is assumed to be strongly convex,
	see \cref{asmp:OC},
	its convex conjugate $g_+\conjugate$
	has a Lipschitz continuous subdifferential $\partial g_+\conjugate$.
	Further,
	from $u_k \rightharpoonup \bar u$ in $L^2(t_0,T)$
	and the assumed compactness of $\nabla F$,
	we get
	$\nabla F(u_k)\rightarrow \nabla F(\bar u)$ in $L^2(t_0, T)$.
	Combining these two facts with
	$\partial g_+\conjugate = (\partial g_+)^{-1}$
	and
	the strong convergence of $\hat u_k$ on $\mathcal{I}_n$
	due to \eqref{eq:hat_u_conv},
	we get
	\begin{equation}
		\label{eq:TRM_strong_convergence_h4}
		\bar u
		\leftharpoonup
		u_k
		=
		\hat u_k
		+
		\partial g_+\conjugate\left(\frac{\hat u_k}{r}  - \nabla
		F(u_k)\right)
		\rightarrow \partial g_+\conjugate( - \nabla F (\bar u))
		\qquad\text{in } L^2(\mathcal{I}_n).
	\end{equation}
	Hence, there exists $K \in \N$, $K \ge n$, such that
	\begin{equation*}
		\norm{u_k - \bar u}_{L^2(\mathcal{I}_n)}
		\le
		\frac{\varepsilon}{2}
		\qquad\forall k \ge K
	\end{equation*}
	and with \eqref{eq:TRM_strong_convergence_h2}  we get
	\begin{equation*}
		\norm{u_k - \bar u}_{L^2(t_0, T)}
		\le
		\varepsilon
		\qquad\forall k \ge K.
	\end{equation*}
	This shows
	$u_k \rightarrow \bar u$ in $L^2(t_0,T)$.
	It remains to show that $\bar u$ is stationary,
	i.e.,
	$C(\bar u) = 0$.

	Regarding $C_{\prox}$, we note that \eqref{eq:TRM_strong_convergence_h4} implies
	$\bar u = \prox_{r g_+}( \bar u - r \nabla F(\bar u))$ on $\mathcal{I}_n$
	for all $n \in \N$.
	From \eqref{eq:TRM_strong_convergence_h3}, we get
	$\abs{ \set{\bar\alpha = 1} \setminus \mathcal{I}_n} \le 2^{-n+1}$
	and, thus,  $\set{\bar\alpha = 1} = \set{\sgn(\bar u) = 1}$
	coincides with $\bigcup_{n = 1}^\infty \mathcal{I}_n$
	(up to a set of measure zero).
	Consequently,
	$\bar u = \prox_{r g_+}( \bar u - r \nabla F(\bar u))$ a.e.\ on $\set{\sgn(\bar u) = 1}$.
	Hence, $C_{\prox}(\bar u) = 0$.

	To show the stationarity of the switching points,
	we fix a switching point $s \in J_{\bar u}$ of $\bar u$.
	From the convergence $u_k \to \bar u$ in $L^1(t_0,T)$,
	we infer the existence of
	a switching point $s_k \in J_{u_k}$ of $u_k$ with
	$\abs{s - s_k} \to 0$.
	From $C(u_k) \to 0$,
	we get
	$C_{\text{switch}}(u_k, 1, s_k) \to 0$.
	The continuity properties of the gradient $\nabla F$ from \cref{asmp:OC}
	ensure
	$\nabla F(u_k)(s_k) \to \nabla F(\bar u)(s)$
	and, consequently,
	$V_{u_k,1}(s_k) \to V_{\bar u, 1}(s)$.
	Altogether, this implies $C_{\text{switch}}(\bar u, 1, s) = 0$.
	Since $s$ was an arbitrary switching point,
	we get $C_{\text{switch}}(\bar u) = 0$.
	Thus, $C(\bar u) = 0$ and $\bar u$ is stationary.
\end{proof}
\subsection{Fast solution of subproblems}
\label{subsec:solution_subproblem}
It remains to give a possibility for the computation of the oracle function. We use the approach from
\cite[Section 5]{Marko2023} to solve the corresponding trust-region subproblems efficiently
using the Bellman principle. 

We discretize the problem by choosing an equidistant discretization
with step length $\tau = (T - t_0)/N_T$ and the controls are discretized by piecewise constant functions.
Thus, the controls are identified with matrices $u \in \R^{N \times N_T}$,
i.e., the entry $u_{i,j}$ corresponds the value of component $i$ on the time interval $[t_{j-1}, t_{j})$.
Consequently, the gradient $\nabla F(u)$ can be represented by a matrix
in $\R^{N \times N_T}$,
i.e., $\nabla F(u)_{i,j}$ is the partial derivative of $F$ at $u$ w.r.t.\ $u_{i,j}$ divided by $\tau$.
In each iteration of the trust-region method,
we have to minimize
a discretization of the model problem \eqref{eq:OC_TRM}
(dropping the constant terms which do not depend on $w$),
i.e.,
\begin{align*}
	m_{u,\Delta}(w) &\coloneqq
	\tau \sum_{j=1}^{N_T}\sum_{i=1}^{N}
	\parens*{
		\nabla
		F(u)_{i,j}w_{i,j}
		+ g_i(w_{i,j})
		+ \frac{\sgn(u_{i,j})\sgn(w_{i,j})}{2\Delta}\abs{u_{i,j}
		- w_{i,j}}^2
	}
	\\
	&\qquad
	+
	\sum_{j=0}^{N_T} \sum_{i=1}^N\abs{\sgn(w_{i,j}) - \sgn(w_{i,j+1})}
	.
\end{align*}
Note that the third term in parenthesis is the discretization of
$(2\Delta)^{-1}\norm{u-w}^2_{L^2(\set{\sgn(uw)=1})}$.

In order to account for the boundary terms in the TV regularization,
we use the convention
$\sgn(w_{i,0}) = \sgn(w_{i,N_T+1}) = 0$.
The components of $w$ are coupled due to the trust-region constraint
$\sum_{i = 1}^N \sum_{j = 1}^{N_T} \abs{ \sgn(u_{i,j}) - \sgn(w_{i,j})} \le
\BB$. The unconventional symbol choice for this trust-region radius
comes from the fact that
this value can be understood as a ``budget'' for the amount of change
regarding the switching pattern in the application of Bellman's principle below.
Originally, we would have $\BB = \Delta / \tau$,
but we will see that it is beneficial to treat
the proximal parameter $\Delta$
and the trust-region radius $\BB$ appearing in the constraint differently.
The oracle $\hat u = T(u, \Delta, \BB)$ is a minimizer of
the model over the trust region.

The main idea will be to split the model problem and apply the celebrated Bellman
principle.
In what follows,
we use the abbreviation $w_j \coloneqq (w_{1,j}, \ldots, w_{N,j})$
for the variable at time step $j$.
To this end we define for $1\le l \le r \le N_T$
the family of partial model problems
via
\begin{align*}
	m_{u,\Delta}^{(l,r)}(w) &\coloneqq \tau \sum_{j=l}^{r}\sum_{i=1}^{N}
	\parens*{
		\nabla
		F(u)_{i,j}w_{i,j}
		+ g_i(w_{i,j})
		+ \frac{\sgn(u_{i,j})\sgn(w_{i,j})}{2\Delta}\abs{u_{i,j}
		- w_{i,j}}^2
	}
	\\
	&\qquad
	+
	\sum_{j=l}^{r-1} \norm{\sgn(w_{j}) - \sgn(w_{j+1})}_1
	.
\end{align*}
Note that $m_{u,\Delta}^{(1,N_T)}(w)$
does not contain the jump terms at the boundary,
i.e.,
\begin{equation*}
	m_{u,\Delta}(w)
	=
	m_{u,\Delta}^{(1,N_T)}(w)
	+
	\norm{\sgn(w_{1})}_1 + \norm{\sgn(w_{N_T})}_1
	.
\end{equation*}
By optimizing a part of the model function
only w.r.t.\ $w_1, \ldots, w_l$
gives rise to the definition of the value function
\begin{equation*}
	\Phi(l,\alpha,B) \coloneqq \min\set*{
		\begin{aligned}
			&m^{(1,l)}_{u,\Delta}(w)
			+ \norm{ \sgn(w_{1}) }_1
		\end{aligned}
		\given
		\begin{aligned}
			&w_1,\ldots,w_l\in \swpz,\\
			&\sgn(w_l) = \alpha, \\
		&\sum_{j=1}^{l}\norm{\sgn(u_j) - \sgn(w_j))}_1 = B
	\end{aligned}
}
\end{equation*}
for $l \in \set{1,\ldots,N_T}$, $\alpha \in \set{0,1}^N$, $B \in \set{0,\ldots, \BB}$.
Note that $\alpha$ encodes the terminal value of $\sgn(w_l)$
and this will be important for the recursion formula below.
Now, we can rewrite the value of the model problem at the minimizer $\hat u = T(u, \Delta, \BB)$ via
\begin{align*}
	m_{u,\Delta}(\hat u) &=
	\min\set*{
		\Phi(N_T,\alpha,B)
		+
		\norm{\alpha}_{1}
		\given
		\alpha\in \{0,1\}^N,\ B \in \set{0,\ldots,\BB}
	}
	,
\end{align*}
where $\norm{\alpha}_{1}$ is the cost for a jump at the terminal time $T$.
Now, we derive a recursion formula for the value function.
For $1 \le l < N_T$ we have
\bgroup\allowdisplaybreaks
\begin{align*}
	\MoveEqLeft \Phi(l+1,\alpha,B) \\
	&=\min\set*{
		\begin{aligned}
			&m^{(1,l)}_{u,\Delta}(w)
			+
			m^{(l+1,l+1)}_{u,\Delta}(w)\\
			&+\norm{\sgn(w_1)}_1
			\\&
			+ \norm{\sgn(w_l) - \alpha}_1  
		\end{aligned}
		\given
		\begin{aligned}
			&w_1,\ldots, w_l, w_{l+1}\in \swpz, \\ 
			&\sgn(w_{l+1}) = \alpha,\\
			&\sum_{j=1}^{l + 1}\norm{\sgn(u_j) - \sgn(w_j)}_1 = B
	\end{aligned}}
	\\[10pt]
	&=\min\set*{
		\begin{aligned}
			&m^{(1,l)}_{u,\Delta}(w)
			+
			m^{(l+1,l+1)}_{u,\Delta}(w)\\
			&+\norm{\sgn(w_1)}_1
			\\&
			+ \norm{\beta - \alpha}_1  
		\end{aligned}
		\given
		\begin{aligned}
			&w_1,\ldots, w_l, w_{l+1}\in \swpz, \\ 
			&\beta \in \set{0,1}^N,\ \tilde B \in \set{0,\ldots,B},\\
			&\sgn(w_{l+1}) = \alpha,\ \sgn(w_l) = \beta\\
			&\norm{\sgn(u_{l+1}) - \alpha}_1 = B - \tilde B \\
			&\sum_{j=1}^{l}\norm{\sgn(u_j) - \sgn(w_j)}_1 = \tilde B
	\end{aligned}}
	\\[10pt]
	&= \min\set*{
		\begin{aligned}
			&\Phi(l,\beta,\tilde B)      +
			m^{(l+1,l+1)}_{u,\Delta}(w)\\
			&+\norm{\beta - \alpha}_1
		\end{aligned}
		\given
		\begin{aligned}
			&w_{l+1}\in \swpz,\ \beta\in \{0,1\}^N,\ \tilde B \in \set{0,\ldots,B},\\
			&\sgn(w_{l+1}) = \alpha,\ \norm{\sgn(u_{l+1}) - \alpha}_1 = B - \tilde B
	\end{aligned}}
	.
\end{align*}
\egroup
In the first equality,
we have rewritten the objective,
in the second equality we introduced the additional optimization variables $\beta$ and $\tilde B$,
which allows us to carry out the optimization w.r.t.\ $w_1, \ldots, w_l$
in the third equality explicitly.
The remaining optimization w.r.t.\ $w_{l+1}$
can be carried out
explicitly and
the minimizer (which also depends on $\alpha_i$)
is given by
\begin{equation}
	\label{eq:formula_w}
	w_{i,l+1}^{\alpha_i}
	\coloneqq
	\begin{cases}
		0 &\text{if } \alpha_i = 0, \\
		\prox_{\Delta g_{i,+}}\parens*{ u_{i,j} - \Delta \nabla F(u)_{i,j}}
		&\text{if } \alpha_i = 1 \text{ and } \sgn(u_{i,j}) = 1, \\
		\argmin_{v \in [a_i, b_i]} \parens*{ \nabla F(u)_{i,j} v + g_{i,+}(v)}
		&\text{if } \alpha_i = 1 \text{ and } \sgn(u_{i,j}) = 0. \\
	\end{cases}
\end{equation}

The initial value of the value function
can be computed explicitly as
\begin{equation*}
	\Phi(1,\alpha,B)
	=
	\min\set*{
		\begin{aligned}
			&m^{(1,1)}_{u,\Delta}(w)
			+ \norm{ \alpha }_1
		\end{aligned}
		\given
		\begin{aligned}
			&w_1\in \swpz,\\
		&\sgn(w_1) = \alpha,\ \norm{\sgn(u_1) - \alpha)}_1 = B
	\end{aligned}
}
.
\end{equation*}
Again, this minimization w.r.t.\ $w_1$ can be evaluated as in \eqref{eq:formula_w}.

Finally, we mention how the minimizer $\hat u = T(u, \Delta)$ of the model problem
can be evaluated.
For this, we store in $\Psi(l+1, \alpha, B)$ one of the minimizing values of $\beta$
appearing in the recursion formula for $\Phi(l+1, \alpha, B)$.
Then, the optimal switching vector $\hat\alpha$ can be evaluated via
\begin{subequations}
	\label{eq:computation_of_solution}
	\begin{equation}
		(\hat\alpha_{N_T}, \hat B_{N_T}) = \argmin\set*{
			\Phi(N_T,\alpha,B)
			+
			\norm{\alpha}_{1}
			\given
			\alpha\in \{0,1\}^N,\ B \in \set{0,\ldots,\BB}
		}
	\end{equation}
	and
	\begin{align}
		\hat\alpha_{l} &= \Psi(l+1, \hat\alpha_{l+1}, \hat B_{l+1} )
		,
		&
		\hat B_l &= \hat B_{l+1} - \norm{\sgn(u_{l+1}) - \hat\alpha_{l+1}}_1
	\end{align}
	for $l = N_T-1, \ldots, 1$.
	Finally, the minimizer $\hat u$ is given by
	\begin{equation}
		\hat u_{i,j} = w_{i,j}^{\hat\alpha_{i,j}}
	\end{equation}
	by using \eqref{eq:formula_w}.
\end{subequations}

In the definition of the oracle $T(u, \Delta)$,
see \cref{def:OC_TRM},
the value $\Delta$ is used both as a trust-region radius
and as a parameter scaling the proximal term.
In the discretized setting,
this corresponds to the fixed choice
$\BB = \Delta/\tau$.
If we use this fixed scaling in \cref{alg:ATRM},
the algorithm barely makes any progress,
since the trust-region radius
and the proximal-point parameter should scale very differently.
Consequently, we separated both parameters.

The utilization of the value function and the separation of $\Delta$ and $\BB$ have a further advantage.
Once we have computed $\Phi$ and $\Psi$,
we can optimize the model $m_{u,\Delta}$
over a smaller trust region with radius $\hat\BB \le \BB$.
That is, the minimizer $\hat u$ of
\begin{equation*}
	\min\set*{
		m_{u, \Delta}(w)
		\given
		\sum_{j = 1}^{N_T} \norm{ \sgn(u_{j}) - \sgn(w_{j})}_1 \le \hat\BB
	}
\end{equation*}
can be obtained
by replacing $\BB$ by $\hat\BB$ in
\eqref{eq:computation_of_solution}.
Note that this is computationally cheap
and does not require the
expensive recomputation of $\Phi$ and $\Psi$.
The separation of $\Delta$ and $\BB$
and the possible fast reevaluation for smaller trust-region radii $\BB$
give rise to \cref{alg:TRM_outer_solver_2}
which is a variation of \cref{alg:ATRM}.
Note that the convergence of \cref{alg:TRM_outer_solver_2}
is not covered directly by the convergence analysis of \cref{thm:OC_iterate_convergence}.
The investigation of convergence guarantees of \cref{alg:TRM_outer_solver_2} is subject to future research.
\begin{algorithm2e}[htp]
	\caption{TRM for \eqref{eq:OC}}
	\label{alg:TRM_outer_solver_2}
	\KwData{%
		$\eta\in (0,1)$,
		$0 < \gamma_1 < 1 \le \gamma_2$,
		$0 < \Delta_0 < \Delta_{\max} \in (0,\infty]$,
		$\BB_0 \coloneqq \BB_{\max}\in \N$,
		$u_0 \in \Uad$
	}
	$n\coloneqq0$\\
	\While{$C(u_n)>0$}
	{    \eIf{$\ared(u_n,\Delta_n,\BB_n) \ge \eta\
		\pred(u_n,\Delta_n,\BB_n)$}
		{
			$u_{n+1} \coloneqq T(u_n, \Delta_n,\BB_n)$\\
			$\Delta_{n+1} \coloneqq \min\{\gamma_2\Delta_n,\Delta_{\max}\}$\\
			$\BB_{n+1} \coloneqq \min\{\ceil{\gamma_2 \BB_n} + 1, \BB_{\max}\}$
		}
		{
			$u_{n+1} \coloneqq u_n$\\
			\eIf{$\BB_{n} \coloneqq 0$}
			{
				$\Delta_{n+1} \coloneqq \gamma_1 \Delta_n$\\
				$\BB_{n+1} \coloneqq \BB_{\max}$
			}
			{
				$\BB_{n+1} \coloneqq \floor{\gamma_1 \BB_{n}}$\\
				$\Delta_{n+1} \coloneqq \Delta_n$
			}
		}
		$n \coloneqq n + 1$
	}
\end{algorithm2e}

\section{Numerical experiments} 
\label{section:numerics}
In this section, we solve (discretizations of) \eqref{eq:OC}
with \cref{alg:TRM_outer_solver_2} by using
the approach of \cref{subsec:solution_subproblem} to solve the trust-region subproblems.
We intend to highlight
characteristic properties, as well as shortcomings of the algorithm.
We will consider examples that come with a qualitative understanding of the
resulting control, especially its switching behavior.
In both examples, the parameters used in \cref{alg:TRM_outer_solver_2} are
\begin{align*}
	\gamma_1 &= \frac12
	,
	&
	\gamma_2 &= 2
	,
	&
	\Delta_0 &= 10^{-7}
	,
	&
	\Delta_{\max} &= 10
	,
	&
	\eta &= 10^{-3}
	.
\end{align*}
The termination criterion ``$C(u_n) = 0$'' in \cref{alg:TRM_outer_solver_2}
is replaced by
\begin{equation*}
	C_{\text{prox}}(u_n) \le 10^{-10}
	\qquad\text{and}\qquad
	\BB_n = 0
	.
\end{equation*}
The first condition ensures that the continuous parts of the control
are almost stationary.
The second condition ensures that the algorithm did not found
reasonable steps and the budget $\BB_n$ for moving the switching
points has been reduced to $0$.
Note that $C_{\text{switch}}(u_n)$ will not converge to zero
after discretization and this will be discussed below.

\subsection{Decay problem}

As a first example,
we consider an identification problem
in which the terms of the objective are defined via
\begin{equation}
	\label{eq:decay}
	\begin{aligned}
		F(u) &\coloneqq \frac{\sigma_y}{2}\int_{t_0}^{T}\left(\frac{y_u(t) - y_d(t)}{y_d(t)}\right)^2
		\dt + \frac{\sigma_T}{2}(y_u(T)-y_d(T))^2
		,\\
		G(u) &\coloneqq\int_{t_0}^{T} g(u(t))\dt,\\
		J(u) &\coloneqq F(u) + G(u) + \sigma_{\sgn}\TV(\sgn(u))
	\end{aligned}
\end{equation}
with parameters $t_0=0$, $T = 140$, $\sigma_y = 10$, $\sigma_T = 0.3$,
$\sigma_{\sgn}  = 1$,
and
\begin{equation*}
	g(x) =
	\begin{cases}
		0& \text{if }x= 0,\\
		0.7x^2 - 0.5x + 0.4&\text{if }x\in [0.3,1],\\
		\infty&\text{else},
	\end{cases}
\end{equation*}
i.e., $a = 0.3$ and $b = 1$.
Further, the state $y_u$ is the solution of the ordinary differential equation
\begin{equation*}
	\dot y_u(t) = -(0.025 + 0.05 u(t))y_u(t),
	\qquad
	y_u(t_0) = 1000.
\end{equation*}
Note that this differential equation describes an exponential decay of the state $y_u$
and the control can be used to increase the decay rate from $0.025$ to a value from
$[0.04, 0.075]$.
The choice of parameters comes from balancing the order of magnitude of
the different terms in the objective,
but have no further meaning.
We use \cref{alg:TRM_outer_solver_2}
to solve discretized versions of this problem.

\begin{figure}[!ht]
	\centering
	\begin{subfigure}{0.45\linewidth}
		\includegraphics[width=\linewidth]{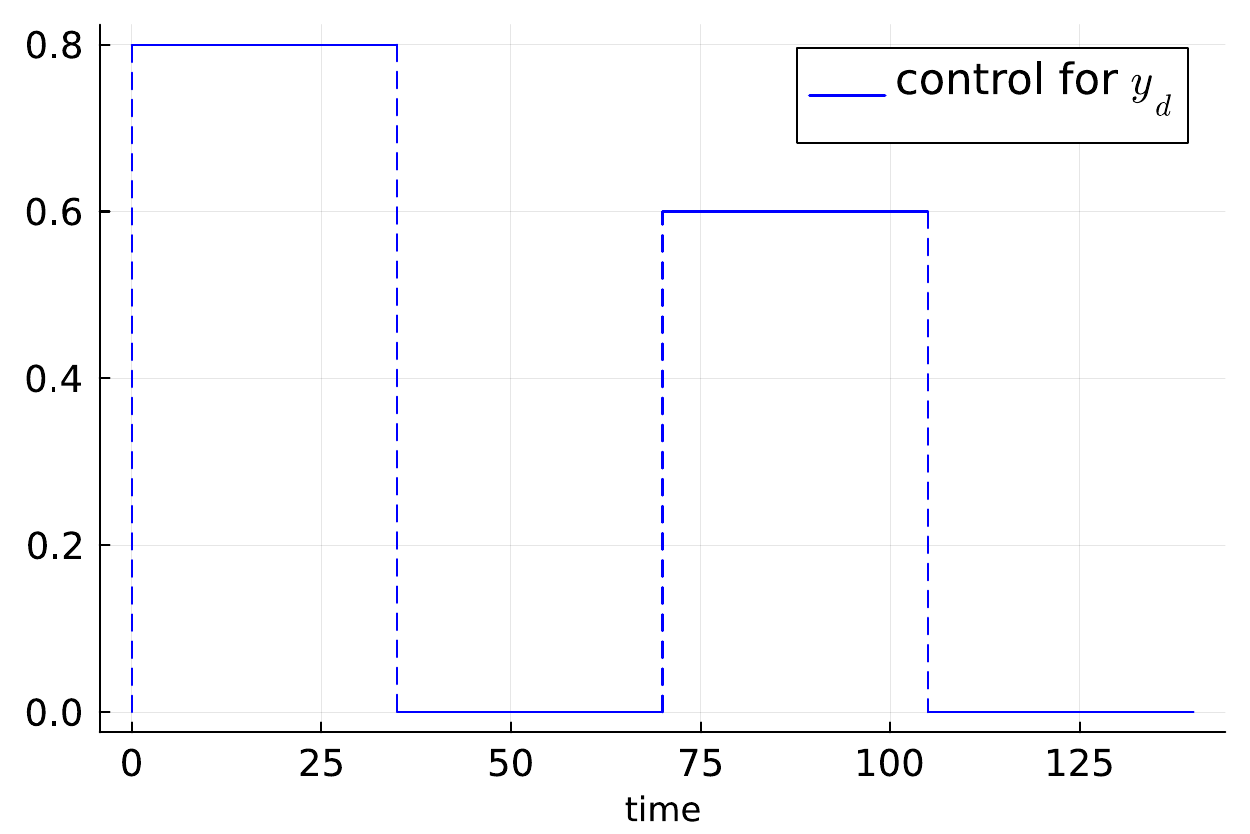}
		\caption{}
	\end{subfigure}
	\hfil
	\begin{subfigure}{0.45\linewidth}
		\includegraphics[width=\linewidth]{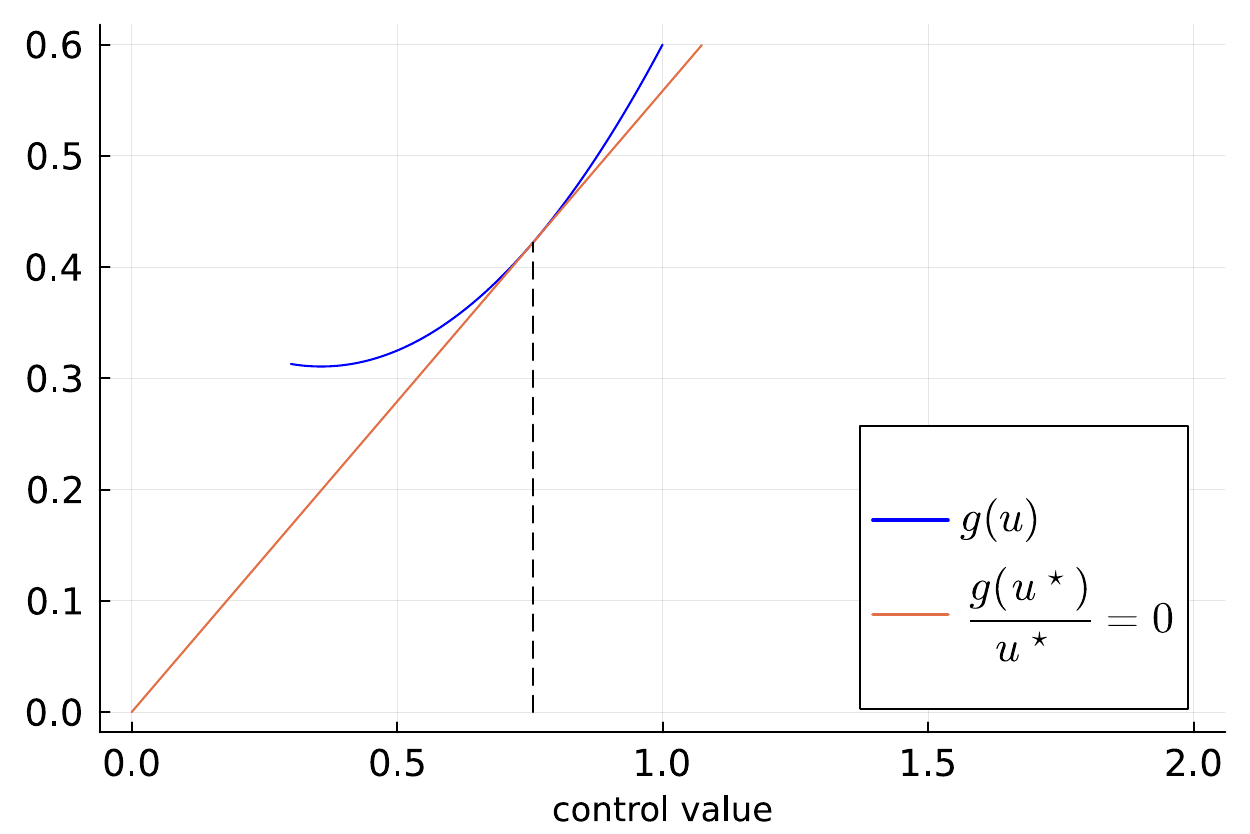}
		\caption{}
	\end{subfigure}

	\begin{subfigure}{0.45\linewidth}
		\includegraphics[width=\linewidth]{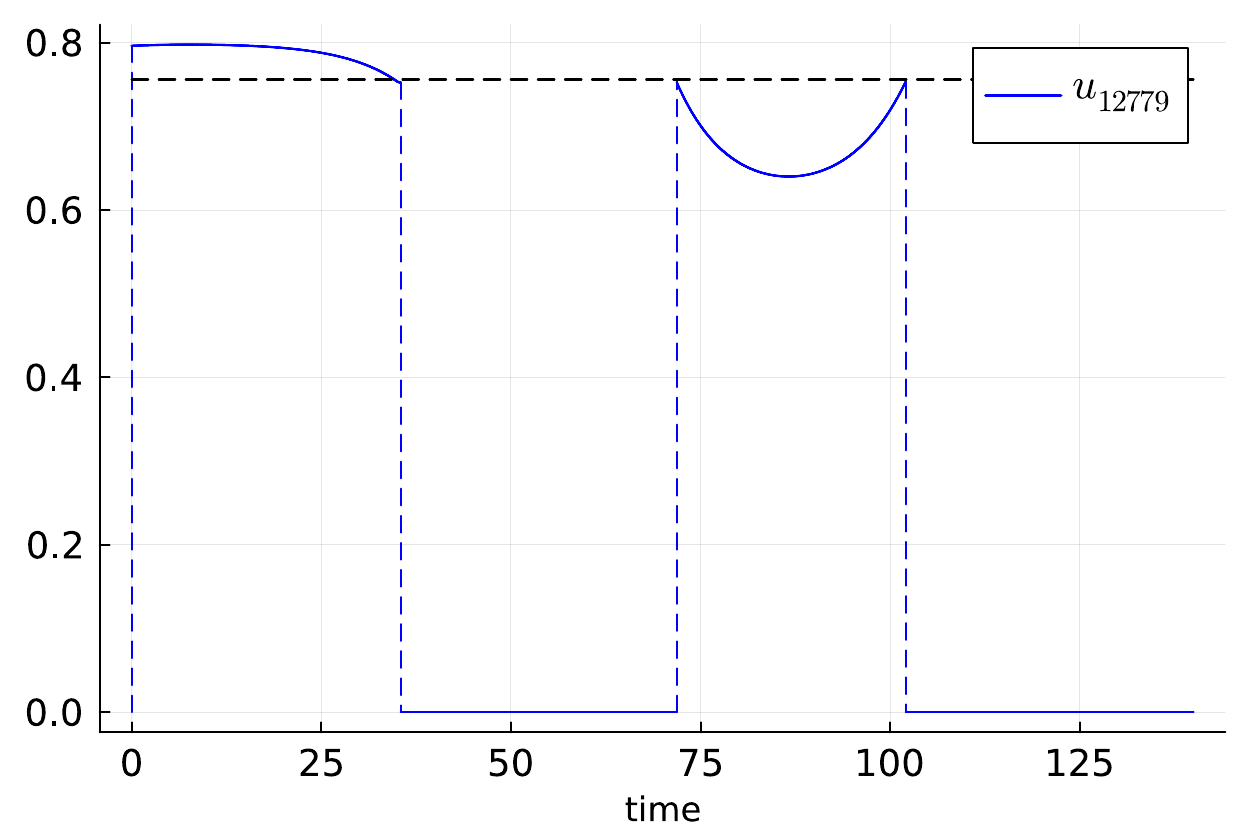}
		\caption{}
	\end{subfigure}
	\hfil
	\begin{subfigure}{0.45\linewidth}
		\includegraphics[width=\linewidth]{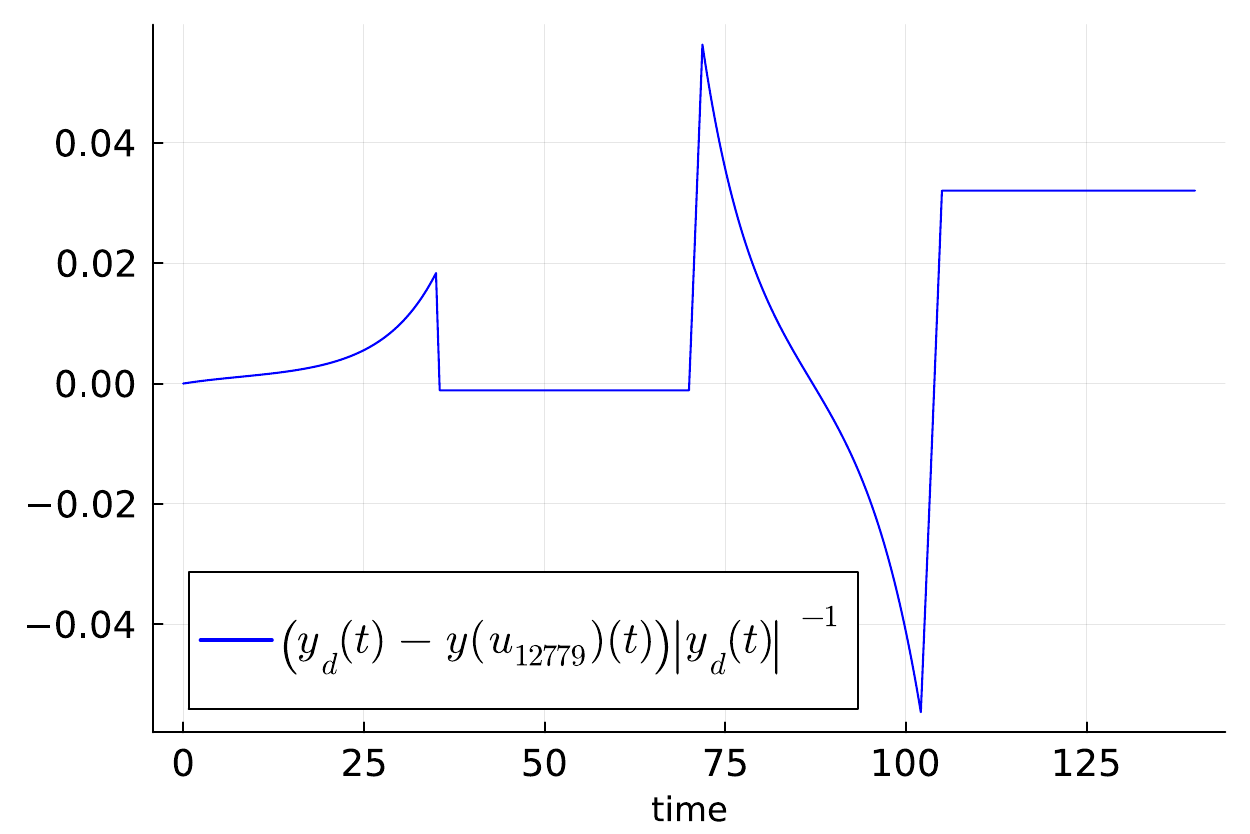}
		\caption{}
	\end{subfigure}

	\begin{subfigure}{0.45\linewidth}
		\includegraphics[width=\linewidth]{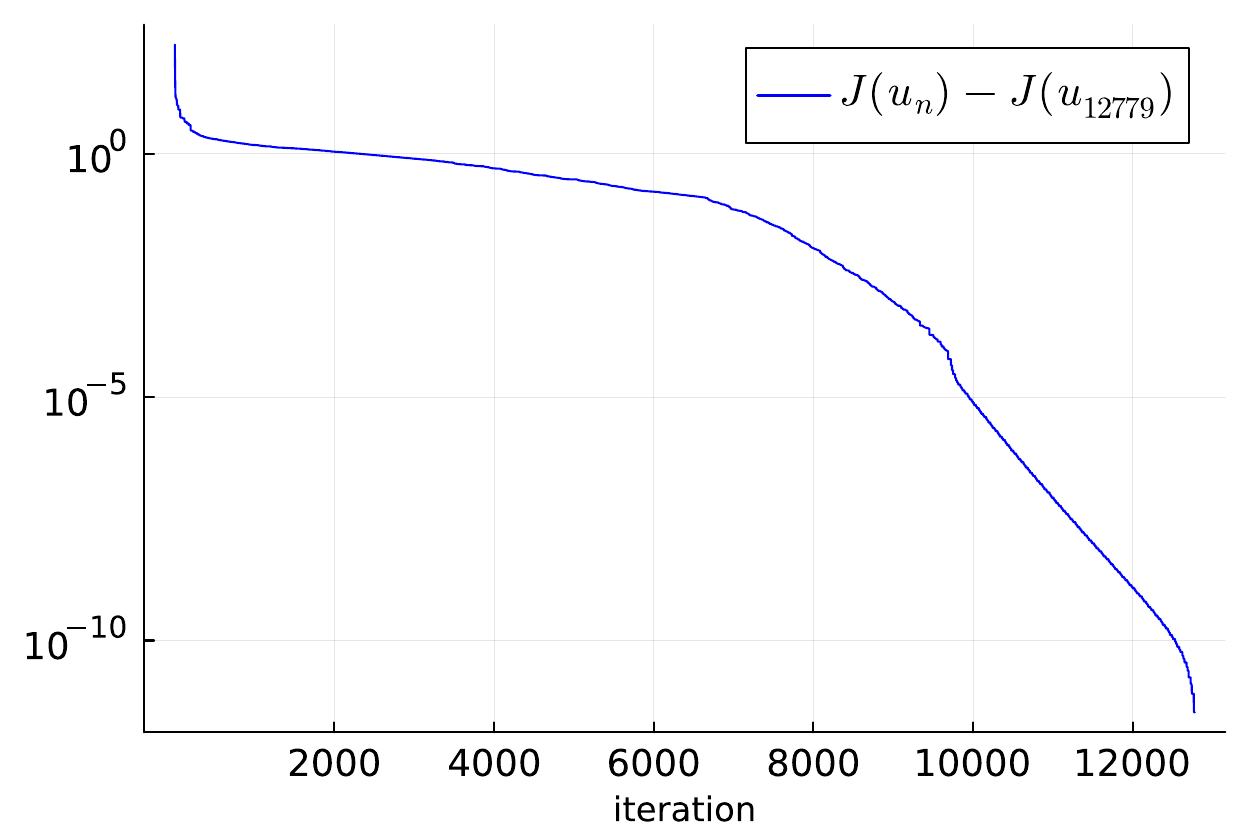}
		\caption{}
	\end{subfigure}
	\hfil
	\begin{subfigure}{0.45\linewidth}
		\includegraphics[width=\linewidth]{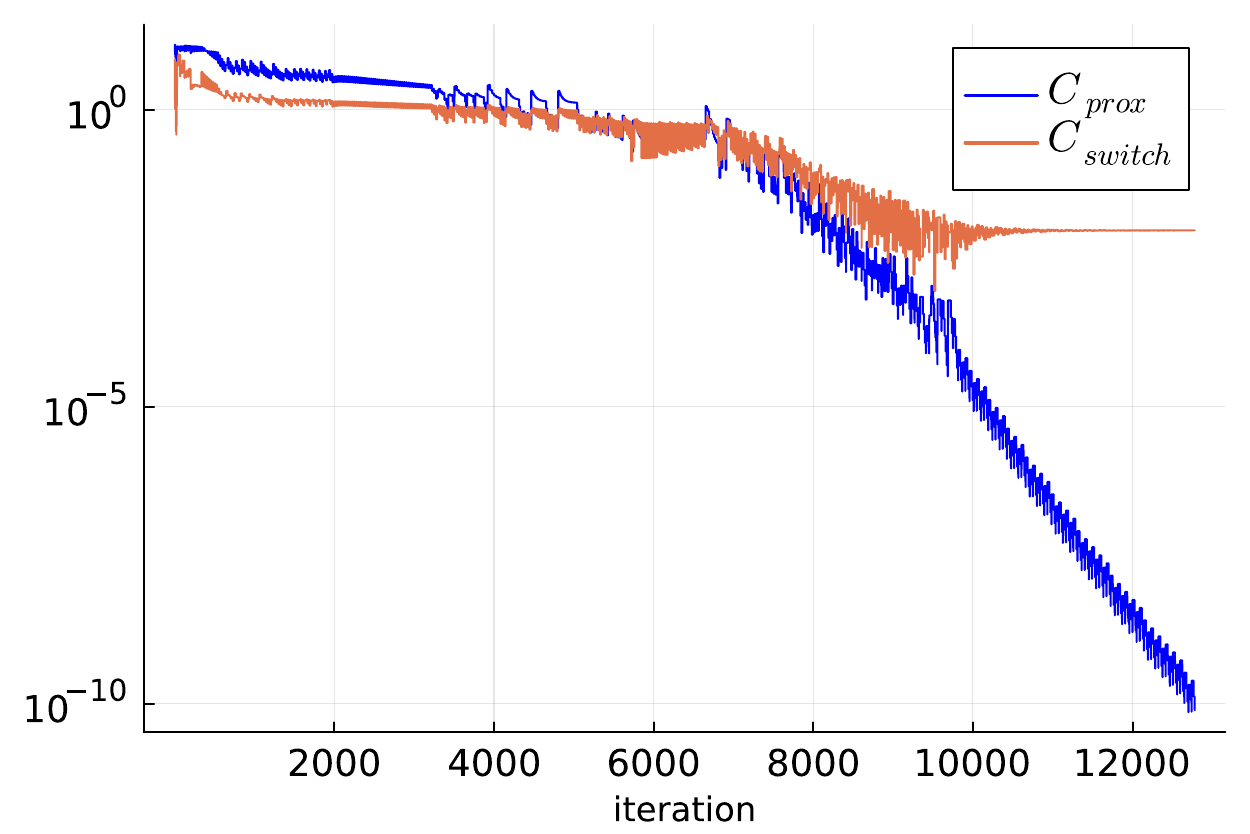}
		\caption{}
	\end{subfigure}
	\caption{Identification of the decay problem. }
	\label{fig:Decay}
\end{figure}
The target state $y_d$ is computed using the control shown in \cref{fig:Decay}(a).
Due to the structure of the state equation, $y_d$ is strictly positive.
The running cost of the active control is priced with the convex function $g$,
see
\cref{fig:Decay}(b). The tangent through the origin (shown in orange) gives the prognosis for the
control value $u^\star$ at which switches are allowed, see \cref{rem:switching_point_optimality}.
The problem is solved on a discretization with 8192 time steps.
\cref{fig:Decay}(c) depicts the computed solution.
The dotted vertical lines show the switches of the control.
The
dotted horizontal line visualizes the precomputed control value $u^\star$
(see \cref{rem:switching_point_optimality})
at which a jump may
occur. Note the exception at $t_0$, where the optimality
condition only yields that the control value is above the switching value.
\cref{fig:Decay}(d) shows the pointwise relative error between target state and computed
state. Note that the plot does not show jumps, but steep slopes between the points,
where the computed control and the control for the computation of $y_d$ have 
slightly shifted switching points.  \cref{fig:Decay}(e) visualizes the objective over the iterations of
the algorithm compared to the objective value obtained in the final iteration
and \cref{fig:Decay}(f) shows both parts of the criticality
measure. 
\begin{table}[!ht]
	\center
	\begin{tabular}{rrrrlll}
		\toprule
		$N_T$ & $\BB_{\max}$  & iterations & time in s            & $J$   & $C_{\prox}$          & $C_{\text{switch}}$\\
		\midrule
		32 &    8 &  4499 &    0.44 & 31.4678 & 9.80040e-11 & 0.785092 \\
		64 &    8 &  3496 &    0.57 & 32.1746 & 9.81943e-11 & 0.728040 \\
		128 &    8 &  3115 &    0.88 & 31.8243 & 3.91882e-11 & 0.549244 \\
		256 &   16 &  4179 &    2.67 & 31.5575 & 9.83474e-11 & 0.320774 \\
		512 &   32 &  4671 &    6.04 & 31.4768 & 9.03884e-11 & 0.122810 \\
		1024 &   64 &  5168 &   20.04 & 31.4857 & 8.43537e-11 & 0.062130 \\
		2048 &  128 &  6881 &   68.46 & 31.4955 & 8.73450e-11 & 0.025920 \\
		4096 &  256 &  9519 &  273.32 & 31.5029 & 6.94831e-11 & 0.014348 \\
		8192 &  512 & 12779 & 1157.26 & 31.5069 & 7.76904e-11 & 0.009325 \\
		16384 & 1024 & 15133 & 4568.48 & 31.5089 & 9.46679e-11 & 0.003030 \\
		\bottomrule
	\end{tabular}
	\caption{Computations with different discretizations for the decay problem}
	\label{tab:Decay_times}
\end{table}
In \cref{tab:Decay_times}
we present some numbers concerning the solution for different time
discretization levels $N_T$
and the choice $\BB_{\max} = \max\set{8, \floor{N_T/16}}$.

From this example it should be noted, that the predicted switching value is attained
in the solution and that the switching behavior mimics that of the control
which was used to compute the target state $y_d$. The
target state is identified well.
The criticality measure $C_{\text{prox}}$ is of magnitude
$10^{-10}$
(which means that the algorithm terminated successfully), while $C_{\text{switch}}$
stagnates at
a rather high level (depending on the time discretization),
which at first glance may not
be the expected convergence.
Note that we derived the optimality conditions
\cref{thm:OC_stationarity}
only for the non-discretized problem.
The proof of condition \eqref{eq:OC_stationarity:1}
(corresponding to $C_{\prox}$)
can be easily carried over
to the discretized setting,
but the proofs of
\eqref{eq:OC_stationarity:2}--\eqref{eq:OC_stationarity:4}
cannot be transferred since they rely on small shifts
of the switching points,
which is not possible after time discretization.
In \cref{tab:Decay_times}, we can see that the final value of
$C_{\text{switch}}$ decreases with $N_T$.

\begin{figure}[!ht]
	\centering
	\begin{subfigure}{0.3\linewidth}
		\includegraphics[width=\linewidth]{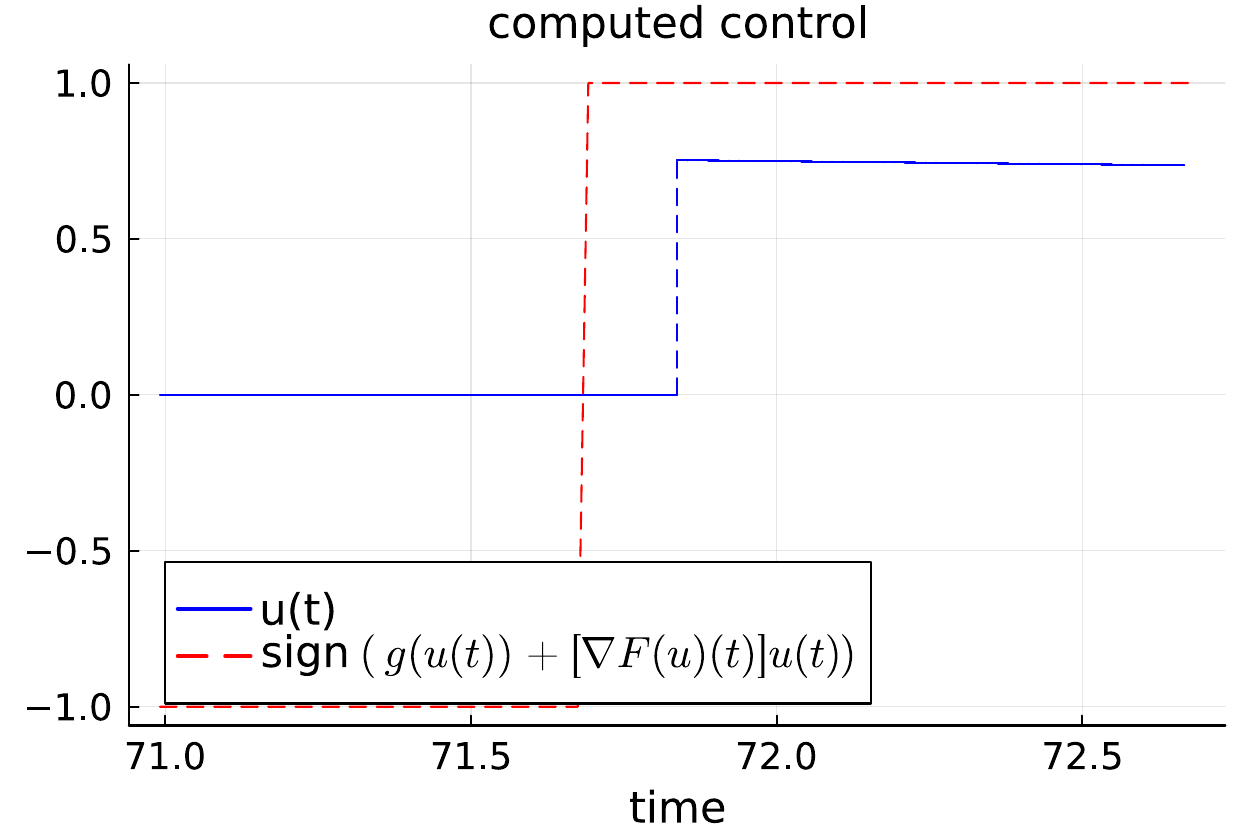}
		\caption{}
	\end{subfigure}
	\hfil
	\begin{subfigure}{0.3\linewidth}
		\includegraphics[width=\linewidth]{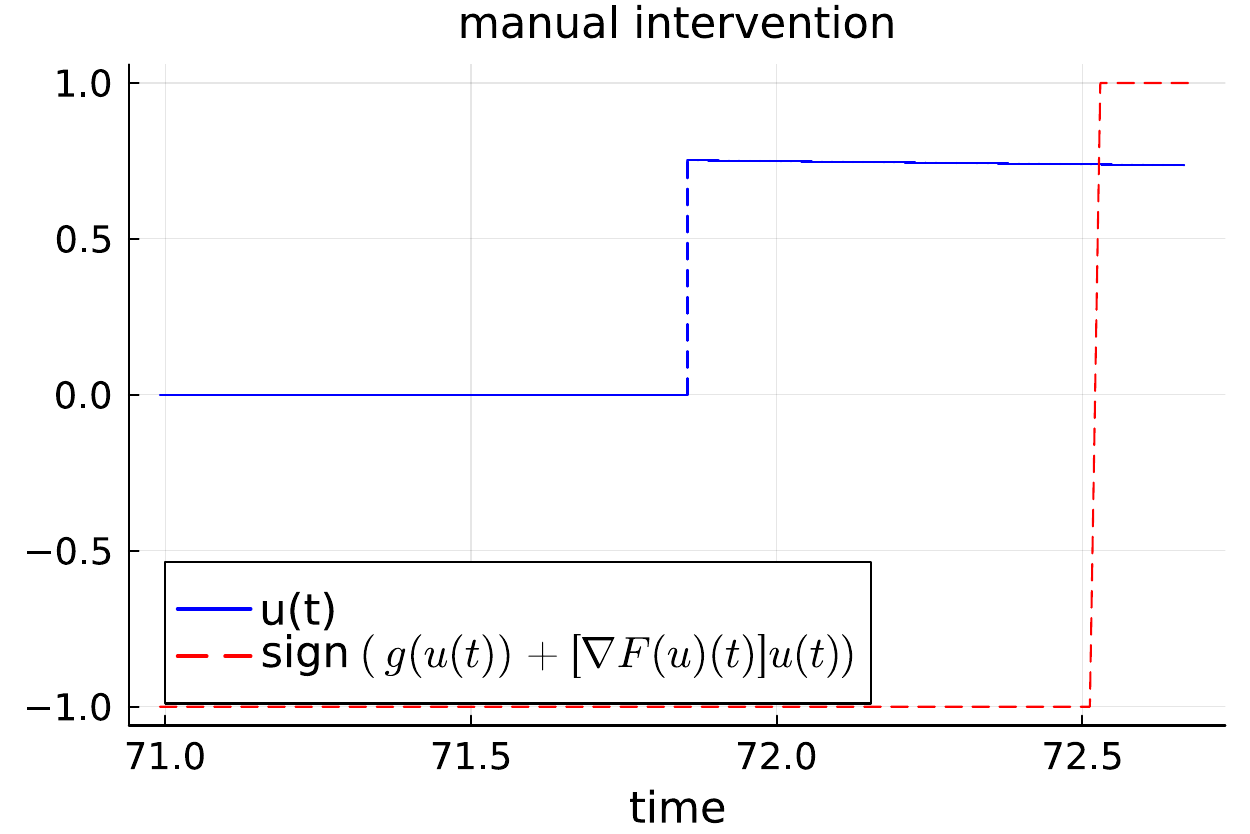}
		\caption{}
	\end{subfigure}
	\hfil
	\begin{subfigure}{0.3\linewidth}
		\includegraphics[width=\linewidth]{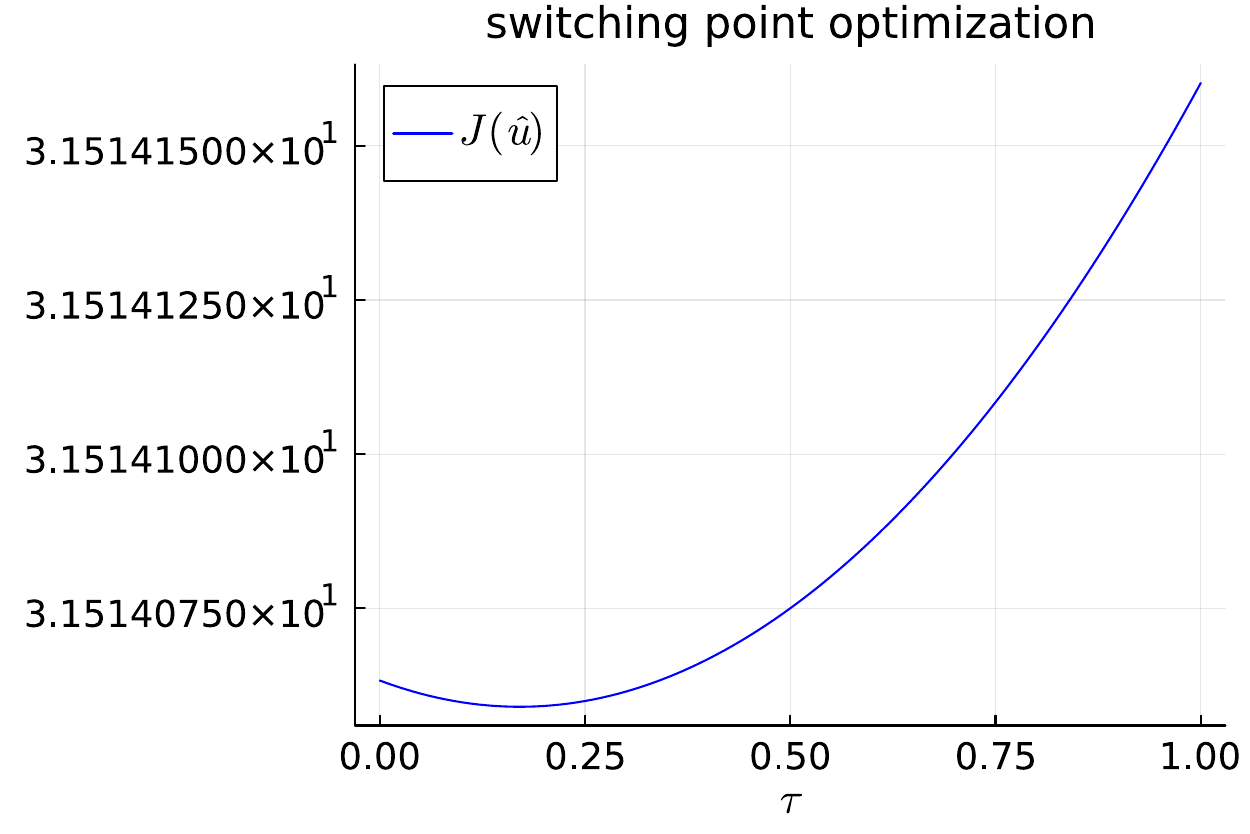}
		\caption{}
	\end{subfigure}
	\caption{Identification of the decay problem. }
	\label{fig:Decay_problems}
\end{figure}

Further visualization of this issue is given in \cref{fig:Decay_problems},
where we zoomed in around the switching point at $t \approx 71.75$.
The computed switching point can be seen in
\cref{fig:Decay_problems}(a) together with the signum of $g(u(t)) + [\nabla
F(u)(t)]u(t)$. The optimality condition suggests, that the sign changes of both
functions occur at the same time point. The value of $g(u(t)) + [\nabla
F(u)(t)]u(t)$ suggests moving the switch to a later
point in time. \cref{fig:Decay_problems}(b) shows the result of manually moving the
switch one time step $\tau$ to the right. We
observe that the manual correction shifts the sign change, such that the
algorithm in its next iteration would reverse our intervention. In this sense
the switching point is indeed locally optimal even if it does not perfectly align
with the exact position of the sign change of the control and the criticality measure retains a rather high value. 
Further, \cref{fig:Decay_problems}(c) shows a gradual movement of the switch between the original
grid points, where
$0$ represents no change and $1$ pushing it a full time step. This suggests
a discretization refinement around the switching points would benefit the
algorithm, but a global refinement of the discretization is not feasible looking
at the scaling behavior in \cref{tab:Decay_times}.

\subsection{SIR problem}
For the second example, we want to consider a more practical example with
multiple controls $u$ controlling a simple SIR model. The
underlying system of ordinary differential equations is given for a population with
a fixed number of $N$
individuals divided into ``susceptible'' ($S$), ``infected''
($I$) and ``recovered'' ($R$) individuals, where the infection is driven by
a contact term with time-dependent infection rate $\beta(t)$ (see below) and recovery occurs with rate
$\rho$.
This results in the system
\begin{equation*}
	\tag{SIR}
	\label{eq:SIR}
	\begin{aligned}
		\dot S(t) &= - \beta(t) S(t)I(t)/N, \\ 
		\dot I(t) &= \phantom{-{}}\beta(t) S(t)I(t)/N - \rho  I(t),\\
		\dot R(t) &= \phantom{- \beta(t) S(t)I(t)/N - {}} \rho I(t).
	\end{aligned}
\end{equation*}
When planning to contain the outbreak of a pandemic,
one possible
strategy is to achieve sufficient suppression in order to wait for a fundamental
change of circumstances, like a new treatment or vaccination.
A second approach is to obtain
a controlled infection while ``flattening the curve''.
We try to model the latter via the objective
\begin{equation*}
	F(I,S) \coloneqq \frac{\sigma_I}{2}\norm{I}^2_{L^2}
	+ \frac{\sigma_S}{2}S(T)^2 
	,
\end{equation*}
where the first term penalizes the severity of the outbreak and the second term
penalizes the remaining number of susceptibles at the end of the considered
time horizon.
The control modifies the infection parameter, such that $\beta(t) = \beta_0
(1 - u_1(t) - u_2(t))$ and the controls are priced with functions representing
a cheap measure

\begin{equation*}
	g_{1}(x) = \begin{cases}
		0 &\text{if } x = 0,\\
		30x^2 + x + 1&\text{if } x \in [0.1, 0.6],\\
		\infty & \text{else},
	\end{cases}\\
\end{equation*}
and an expensive measure
\begin{equation*}
	g_{2}(x) = \begin{cases}
		0 &\text{if } x = 0,\\
		30000x^2 + x + 3000&\text{if } x \in [0.1, 0.4],\\
		\infty& \text{else}.
	\end{cases}
\end{equation*}
The stark contrast is intentionally chosen to highlight shortcomings of the
algorithm and at the same time allows for some intuition on the structure
of the solution. The full objective is given by
\begin{equation*}
	J(u) \coloneqq  F(I(u),S(u)) + \sum_{i=1}^{2}\left(\int_{t_0}^{T}   g_i(u_i(t))\dt
	\right) + \sigma_{\sgn}\TV(\sgn(u)).
\end{equation*}
The chosen parameters are
\begin{equation*}
	\sigma_I = 2, \;
	\sigma_S = 800, \;
	\sigma_{\sgn} = 10000, \;
	\rho = 0.1, \;
	\beta_0 = 0.6, \;
	t_0 = 0, \;
	T = 140
	.
\end{equation*}
\begin{figure}[!ht]
	\centering
	\begin{subfigure}{0.45\linewidth}
		\includegraphics[width=\linewidth]{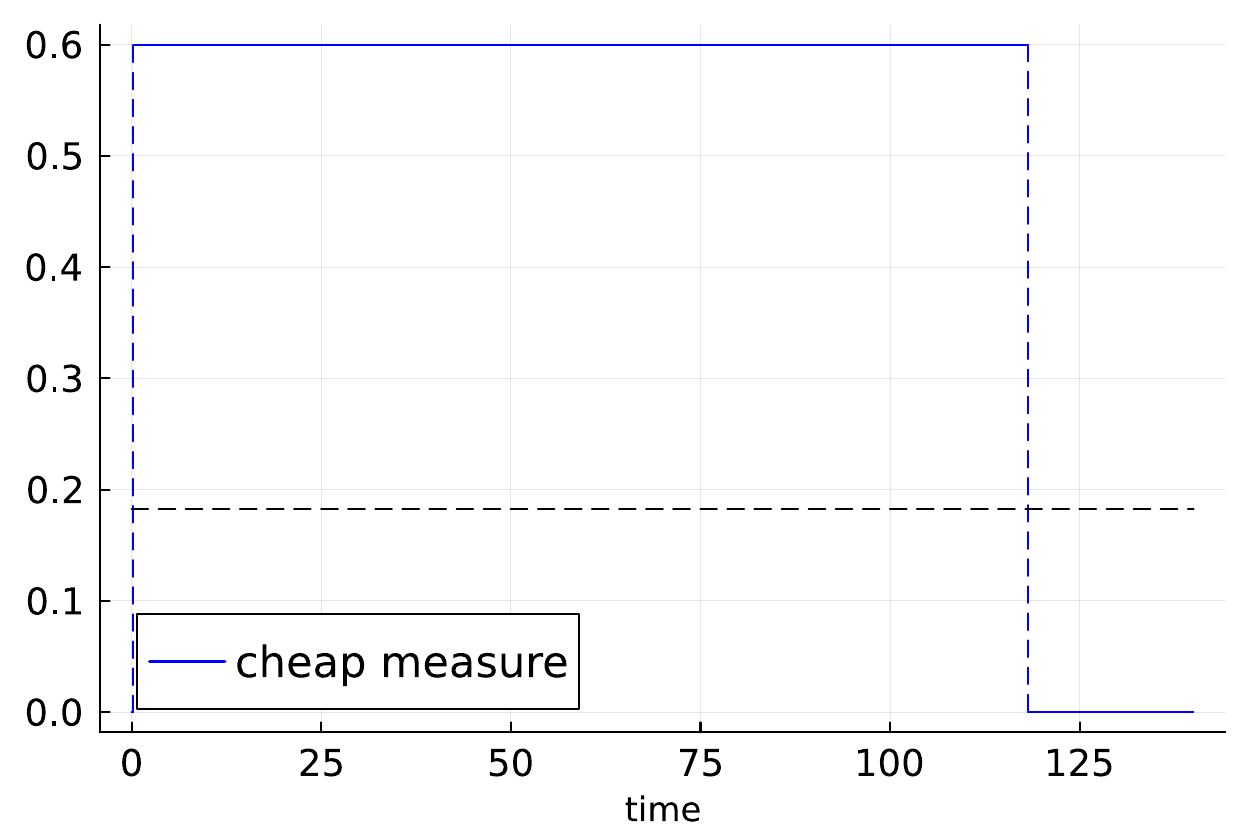}
		\caption{}
	\end{subfigure}
	\hfil
	\begin{subfigure}{0.45\linewidth}
		\includegraphics[width=\linewidth]{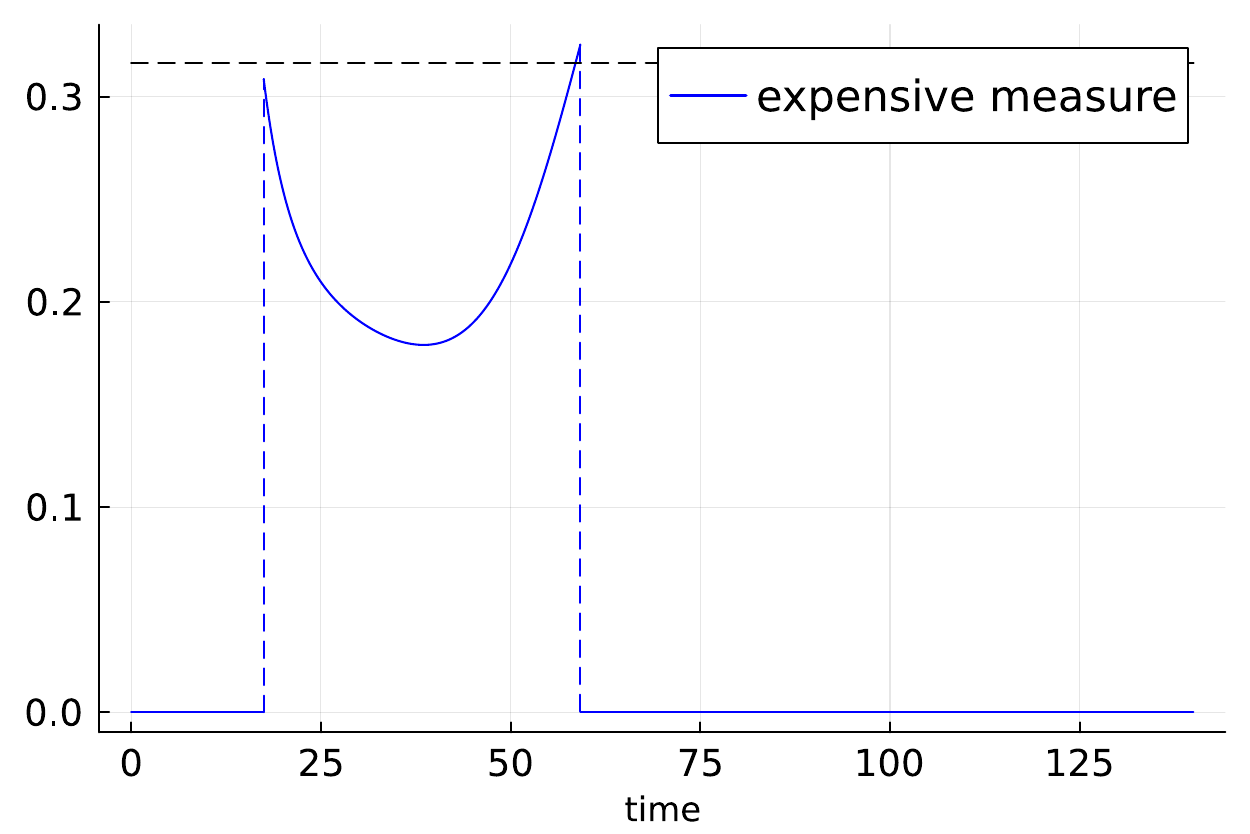}
		\caption{}
	\end{subfigure}

	\begin{subfigure}{0.45\linewidth}
		\includegraphics[width=\linewidth]{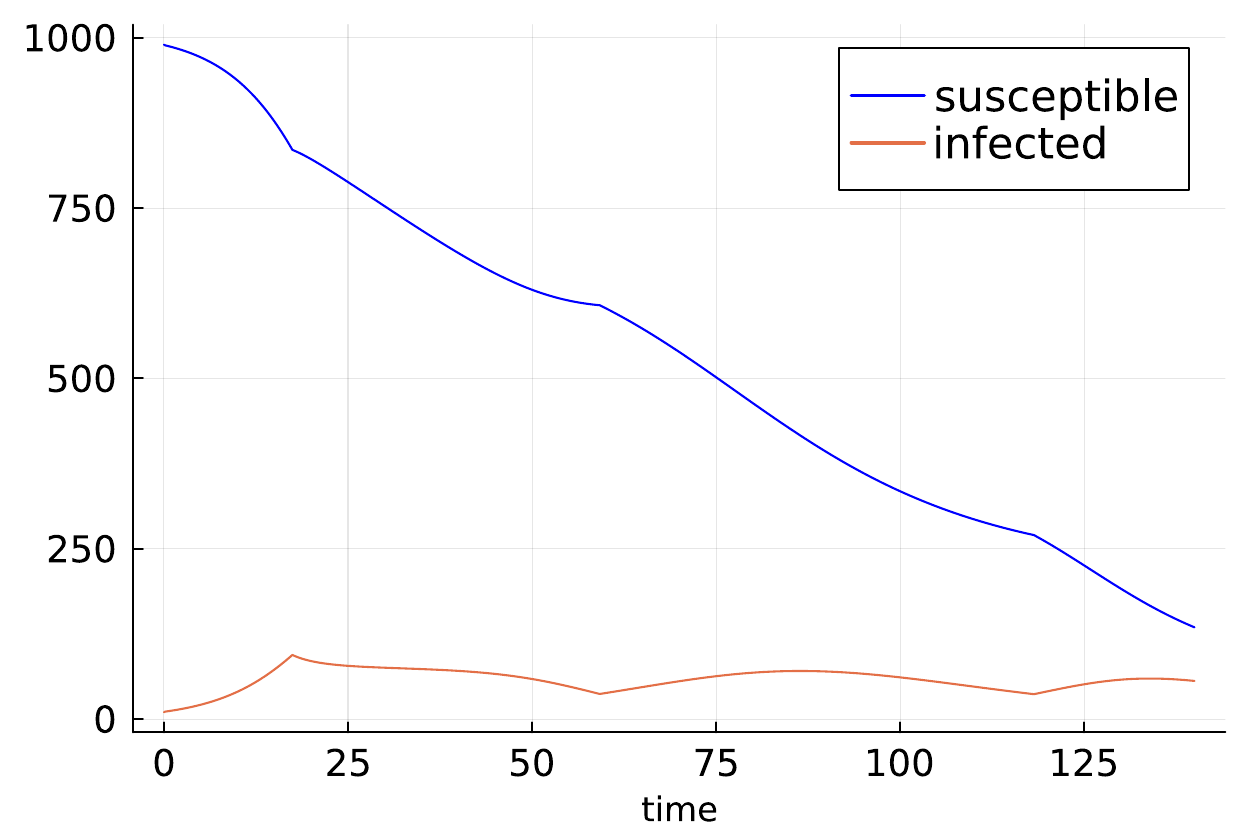}
		\caption{}
	\end{subfigure}
	\hfil
	\begin{subfigure}{0.45\linewidth}
		\includegraphics[width=\linewidth]{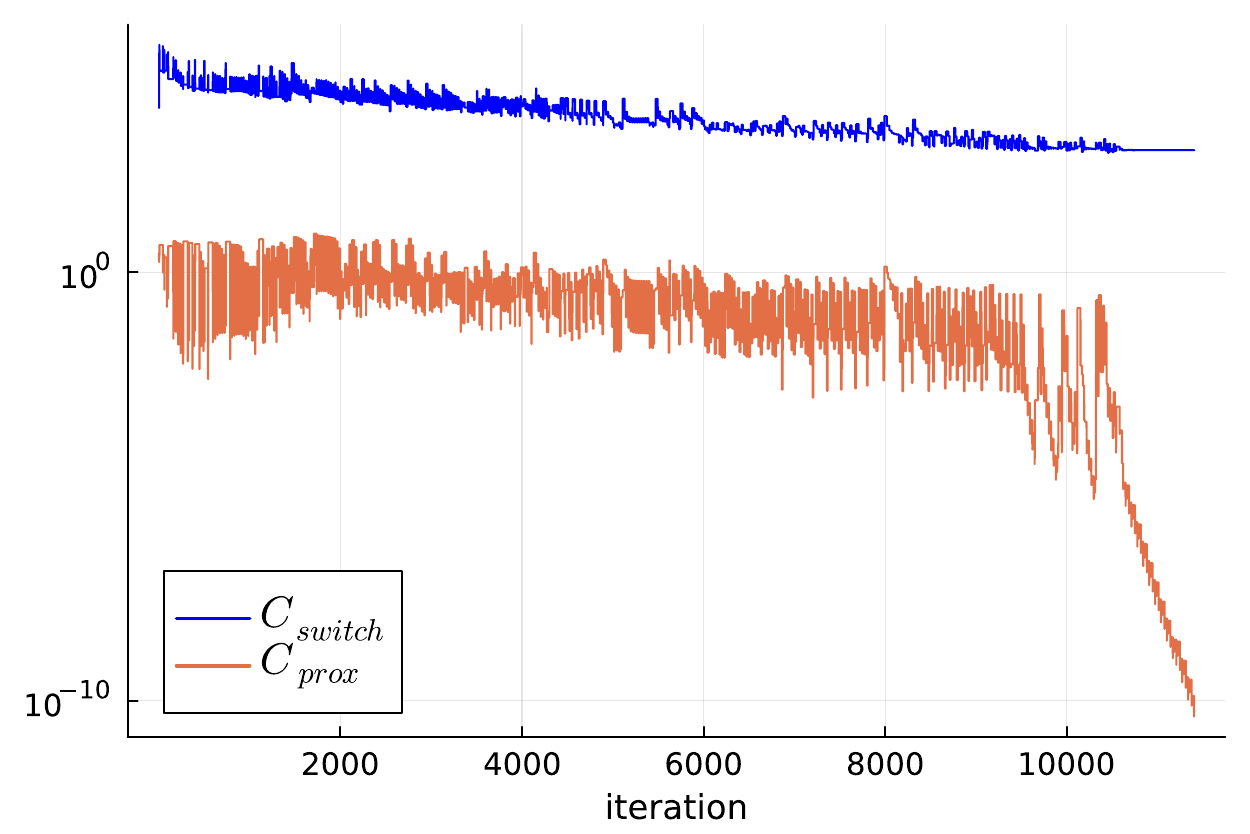}
		\caption{}
	\end{subfigure}

	\caption{Control of an SIR dynamic}
	\label{fig:SIR}
\end{figure}%
The solution for a discretization with 8192 time steps is shown in \cref{fig:SIR}.
From a qualitative standpoint, the solution is reasonable. The suggested control
manages a controlled infection with a flattened curve, while the cheap measure is
used to its fullest capacity and the expensive measure is balanced with its cost,
while the amount of switches are minimized, due to the high switching cost.
On closer inspection, we see an amplified version of the detailed problem
regarding the convergence of $C_{\text{switch}}$. For the cheap control the
switching value (visualized by the dotted line) seems to be ignored, while the expensive
measure does not exactly reach its switching value either.
This suggests that the current time discretization is too coarse to resolve these features
of the optimal control.
\begin{table}[!ht]
	\center
	\begin{tabular}{rrrrllr}
		\toprule
		$N_T$ & $\BB_{\max}$  & iterations & time in s            & $J$   & $C_{\prox}$          & $C_{\text{switch}}$\\
		\midrule
		32 &    8 &  3711  &    1.01 & 1.43464e6 & 8.1725e0-11 & 21168.2  \\
		64 &    8 &  3103  &    1.67 & 1.34994e6 & 9.96217e-11 & 24282.7  \\
		128 &    8 &  1087  &    0.96 & 1.36792e6 & 6.42651e-11 & 38060.6  \\
		256 &   16 &  4110  &    6.03 & 1.12967e6 & 9.84629e-11 &  8838.4  \\
		512 &   32 &  8016  &   31.40 & 1.02804e6 & 9.57865e-11 &  3777.2  \\
		1024 &   64 &  4819  &   43.77 & 9.48856e5 & 7.17141e-11 &  2010.9  \\
		2048 &  128 &  6689  &  152.28 & 8.89106e5 & 9.08372e-11 &  1086.8  \\
		4096 &  256 &  7008  &  483.05 & 8.39067e5 & 8.54805e-11 &  1034.9  \\
		8192 &  512 & 11402  & 2488.06 & 7.92406e5 & 5.96074e-11 &   716.4  \\
		\bottomrule
	\end{tabular}
	\caption{Computations with different discretizations for the SIR problem}
	\label{tab:SIR_times}
\end{table}
In \cref{tab:SIR_times},
we give further results on different time discretization levels.
Similarly to the decay problem, we see that $C_{\text{switch}}$ decreases
as the time discretization parameter $N_T$ increases.

\section{Outlook}
\label{sec:outlook}
In this section we present some interesting questions and ideas that arise from
our work.

\subsubsection*{Switching point optimization between grid points}
We have seen in \cref{section:numerics}, that we do not have the
desired convergence with respect to the part of the criticality measure
$C_{\text{switch}}$ that relates to the optimality of the switching points.
At the same time, we showed for the example problem that
the method indeed numerically converges given the discretization. This suggests
to refine the discretization, but the algorithm presented in \cref{subsec:solution_subproblem} 
that is used to solve the model problems is expensive. This issue occurs even
if $a_i = b_i$.

We suggest to investigate a post processing
step in form of a switching point optimization that operates between grid points. This could be realized by
changing the Euler step from
\begin{equation*}
	y_{j+1} = y_j + (t_{j+1} - t_{j}) f(t_j, y_j, u_{j+1})
\end{equation*}
to
\begin{equation*}
	y_{j+1} = y_j + h (t_{j+1} - t_{j}) f(t_j, y_j, u_{j+1})
\end{equation*}
with $h \in [0,1]$.
The new optimization variable $h$ describes a switch strictly between
the time points $t_j$ and $t_{j+1}$.
This essentially means, that intermediate values can be expressed as a simple proportional
scaling when applying the forward Euler method to solve the state
and adjoint equations required for the computation of $F(u)$ and $\nabla F(u)$. 
Utilizing this principle, a gradient descent method can be
deployed, that moves switches between time points by introducing the appropriate
scaling. The main benefit is that the computed gradient remains numerically
exact. In first experiments, we could confirm convergence with respect to
$C_{\text{switch}}$ for \eqref{eq:SIR} and even the features of
the cheap pricing function were correctly resolved. 

The second benefit is that this approach yields a proper termination criterion which
can significantly shorten the computational time, since one does not have to use the
expensive inner solver until no further step is taken.

\subsubsection*{Generalizations of the control structure}
The structure of the objective was mainly motivated by control problems arising
for pandemic models, where the decision between activating or deactivating
a measure with its associated switching cost is highly important.
The proposed structure has two natural extensions.

First would be the inclusion of a time
dependence of $g_i$, $a_i$ and $b_i$, where we suspect that one then needs at least Lipschitz
continuity with respect to the time variable. We have to emphasize that the
proofs do not extend trivially to this case.

The second question would be regarding multiple activity sets $u_i(t) \in
[a_i^1,b_i^1]\cup [a_i^2,b_i^2] \cup \cdots \cup [a_i^K,b_i^K]$ with 
$a_i^1 \le b_i^1 < a_i^2 \le b_i^2 < \cdots < a_i^K \le b_i^K$ and a cost
function, where we only assume strong convexity for each
$g_i|_{[a_i^k,b_i^k]}$.
Note that this is related to the concept of multibang control,
see, e.g., \cite{ClasonKunisch2016,ClasonTamelingWirth2018}.
We expect that
the algorithm extends naturally to this case, while it is less obvious what to
expect for the jumping value $u^\star$. We conjecture that a jump between two
activity sets can only occur from and to the control value where
$g_i|_{[a_i^k,b_i^k]}$ and $g_i|_{[a_i^l,b_i^l]}$ share the same tangent
and only if this tangent lies below $g_i|_{[a_i^j,b_i^j]}$ for all $i < j < k$.

\subsubsection*{Combination of multiple oracles}
Lastly, \cref{sec:full_problem} poses a question regarding the abstract trust-region
method. Large parts of the proofs involved splitting the argument into the
proximal gradient part and the part relating to the switching point
optimization. However, the proofs could not be carried out in full separation,
especially the verification of the stability assumption
\cref{asmp:ATRM}\eqref{eq:stability}. 

It would be interesting to investigate under which
conditions for the underlying problem two separate oracles, which
individually satisfy \cref{asmp:ATRM}, can be combined.
We suggest to focus on a combined criticality measure which
is simply the maximum of the two individual components.

We suspect, that a strong assumption of interdependence is
required with respect to the Lipschitz continuity, if both oracles are
evaluated in each step and share the same trust-region radius. In contrast, we expect
significantly less strict assumptions if the choice of the deployed oracle is
made dependent on the currently dominating part of the criticality measure.

It is highly likely, that the convergence of \cref{alg:TRM_outer_solver_2} can be
proven by combining our results for the proximal gradient part in
\cref{sec:convex_problem} directly with the results from \cite[Section~5]{Manns2024}
utilizing such an abstract combination result for separate oracle
functions. 

\appendix
\section{Proof of \texorpdfstring{\hyperref[thm:atrm_convergence]{Theorem~2.2}}{Theorem 2.2}}
\label{sec:appendix}
We prepare the proof of \cref{thm:atrm_convergence}
with a simple lemma.
\begin{lemma}
	\label{lem:ATRM_1}
	Let \cref{asmp:ATRM} be valid.
	There exists a nondecreasing function $\Dmax_{\pred} \colon (0,\infty) \to (0,\infty)$
	such that for all $x \in X$ with $C(x) > 0$ and trust-region radii
	$\Delta \in (0,\Dmax_a(x))$
	with
	$\Delta \le \Dmax_{\pred}(C(x))$
	we have
	\begin{enumerate}
		\item [1.](sufficient decrease):\label{lem:ATRM_predictability} 
			\begin{equation*}
				\pred(x,\Delta) \ge \frac{f\left(C(x)\right)}{2}\Delta,
			\end{equation*}
			i.e., for every nonstationary $x$ there exists a trust-region radius that guarantees
			a positive prediction,
		\item [2.]\label{lem:ATRM_accuracy}(accuracy):        \begin{equation*}
				\ared(x,\Delta) \ge \eta\,\pred(x,\Delta).
			\end{equation*}
			In particular, if $\hat x = T(x,\Delta)$ is not accepted in \cref{alg:ATRM},
			we get a lower bound for the trust-region radius $\Delta \ge \min\set{\Dmax_{\pred}(C(x)),\Dmax_a(x)}$.
	\end{enumerate}
\end{lemma}
\begin{proof}
	From straightforward calculations, we obtain
	\begin{equation*}
		f(C(x)) \Delta - c\Delta^s
		\ge
		\frac{1}{2} f(C(x)) \Delta
		\qquad
		\forall
		\Delta \in \bracks*{0, \parens*{\frac{f(C(x))}{2 c}}^{1/(s-1)}}
		.
	\end{equation*}
	Thus, the assertion follows
	from
	\cref{asmp:ATRM}\eqref{eq:sufficient_decrease} and \eqref{eq:accuracy}
	with
	\begin{equation*}
		\Dmax_{\pred}(\zeta)
		\coloneqq
		\parens*{\frac{ f(\zeta)}{2 c}}^{1/(s-1)}
		\qquad
		\forall \zeta > 0
		.
	\end{equation*}
\end{proof}

Now we are in position to prove \cref{thm:atrm_convergence}.
\begin{proof}[Proof of \cref{thm:atrm_convergence}]
	\cref{alg:ATRM} terminates at $x_n$ if and only if $C(x_n)=0$. For the
	remainder of the proof, we assume that the sequence of iterates $(x_n)$ is infinite.

	First, we argue that for each $x_n$ any trust-region radius below a certain
	threshold leads to an accepted step (for clarification of the term compare with
	\cref{rem:ATRM_accuracy}). If $\Delta_n < \Dmax_a(x_n)$ then \cref{lem:ATRM_1}
	yields the threshold $\Dmax_{\pred}(C(x_n))$. If $\Delta_n \ge \Dmax_a(x_n)$ the threshold is
	$\Dmax_b$ by \cref{asmp:ATRM}\eqref{eq:accuracy_E}, especially if
	$\Dmax_a(x_n) = 0$. Hence, at most finitely many scalings with
	$\gamma_1 < 1$ lead to an accepted step.
	Thus, \cref{alg:ATRM} performs infinitely many successful steps if it does
	not terminate.  We continue with proof by contradiction and assume that
	there exists a subsequence
	$(x_{n_k})$ of successful iterations with $C(x_{n_k})>2\varepsilon>0$.
	Using this subsequence, we show that we get a certain decrease of the objective $J$
	and this yields a contradiction, since $J$ is assumed to be bounded from below.

	The radius $\Dmax_{\pred}(\varepsilon)$ from \cref{lem:ATRM_1} and the assumed monotonicity of $f$ yield
	\begin{equation}
		\label{eq:ATRM_pred_est_simple}
		\ared(x,\Delta) \ge \eta\,\pred(x,\Delta) \qquad
		\land\qquad   \pred(x,\Delta) \ge \frac{f(\varepsilon)}{2}\Delta 
	\end{equation}
	whenever
	$\Delta \in (0,\Dmax_a(x))$,
	$\Delta \le \Dmax_{\pred}\left(\varepsilon\right)$,
	and $C(x)>\varepsilon$. We define
	\begin{equation}
		\label{eq:Dmax}
		\Dmax \coloneqq
		\min\left\{\Dmax_b,\Dmax_c,\Dmax_{\pred}\left(\varepsilon\right),\frac{\varepsilon}{L}\right\}.
	\end{equation}

	Next we set $i \coloneqq n_1$.
	First, we investigate the case $\Delta_{i} \ge \Dmax$.
	We define
	$\mathcal{I}_1\coloneqq \{i\}$ and use a further distinction by cases.
	\begin{itemize}
		\item
			If
			$\Delta_i \ge \Dmax_a(x_{i})$ and $\Dmax_a(x_i) \le \Dmax_b$, we use
			\cref{asmp:ATRM}\eqref{eq:accuracy_E} together with the monotonicity
			of the prediction to estimate
			\begin{equation}
				\label{eq:pred_x_i_E}
				\pred(x_{i},\Delta_{i})
				\overset{\eqref{eq:monotonicity}}{\ge}
				\pred(x_{i},\Dmax_a(x_{i})) \overset{\eqref{eq:accuracy_E}}{\ge} \delta. 
			\end{equation}
		\item
			If
			$\Delta_i \ge \Dmax_a(x_{i})$ and $\Dmax_a(x_i) > \Dmax_b$
			then \eqref{eq:Dmax} implies $\Dmax\le \Dmax_b < \Dmax_a(x_i)$.
			Thus,
			\begin{equation}
				\label{eq:pred_x_i}
				\pred(x_{i},\Delta_{i})\overset{\eqref{eq:monotonicity}}{\ge}
				\pred(x_{i},\Dmax) \overset{\eqref{eq:ATRM_pred_est_simple}}\ge
				\frac{f(\varepsilon)}{2}\Dmax \eqqcolon \delta_1 > 0.
			\end{equation}
		\item
			If $\Dmax_a(x_{i}) > \Delta_i$,
			we get
			$\Dmax \le \Delta_i < \Dmax_a(x_i)$
			and we can argue as in
			\eqref{eq:pred_x_i}.
	\end{itemize}

	In the other case $\Delta_{i} < \Dmax$, we define $\mathcal{I}_1
	\coloneqq\{i,\dots,K\}$, where $K \ge i$ is the smallest index such that at least
	one of the following conditions is satisfied:
	\begin{enumerate}
		\item $\pred(x_K,\Delta_K) \ge \delta$,
		\item $\sum_{l=i}^{K+1}\Delta_l \ge \Dmax$,
		\item $\ared(x_{K+1},\Delta_{K+1}) < \eta\,\pred(x_{K+1},\Delta_{K+1})$,
			i.e. iteration $K+1$ is not successful.
	\end{enumerate}
	We briefly argue that $K$ is finite. By construction, $\mathcal{I}_1$
	encodes only successful iterations for which the trust-region radius is
	nondecreasing. Condition (ii) leads to the bound $K\le i+\frac{\Dmax}{\Delta_i} - 1<\infty$.
	By definition, one of the conditions (i)--(iii) is always valid at $K$. Note that
	the construction directly yields $\sum_{l=i}^{K}\Delta_l <
	\Dmax$, since $\Delta_i \le \Dmax$ was assumed.
	In particular, $\Delta_l < \Dmax$ for $l = i,\dots,K$. Next, we perform
	a distinction by cases on which of these conditions hold at $K$.

	\textbf{Case (i):} By construction, the iteration $K$ is successful and
	$\pred(x_K,\Delta_K)\ge\delta > 0$.

	\textbf{Case $\mathbf{\neg(i)}$:} \cref{asmp:ATRM}\eqref{eq:accuracy_E}
	together with $\neg(i)$ implies $\Delta_l\notin
	[\Dmax_a(x_l),\Dmax_b]$ for all $l\in \mathcal{I}_1$.
	Since we already know
	$\Delta_l\le\Dmax \le \Dmax_b$, we get
	\begin{equation}
		\Delta_l < \Dmax_a(x_l)\qquad
		\forall l=i,\dots,K
		\label{eq:neg_i_1}.
	\end{equation}
	Similarly, \cref{asmp:ATRM}\eqref{eq:stability_E} combined with $\neg (i)$
	and
	$\Delta_l\le\Dmax\le\Dmax_c$
	implies
	\begin{equation*}
		\abs{C(x_{l+1})-C(x_l)}\le
		L\Delta_l\qquad
		\forall l=i,\dots,K.
	\end{equation*}
	As a consequence, for any $l = i,\ldots,K$, we can estimate
	\begin{equation}
		\label{eq:C_est}
		C(x_{l+1}) \ge C(x_l) - L\Delta_{l} \ge C(x_i)
		- L\sum_{k=i}^{l}\Delta_{k} \ge C(x_i) - L\Dmax \ge 2\varepsilon
		- L\frac{\varepsilon}{L} = \varepsilon.  
	\end{equation}
	\textbf{Case (ii)$\land \neg$(i):} Since the iterations $l=i,\ldots,K$ are successful,
	$\Delta_{K+1} \le \gamma_2 \Delta_K$ implies
	\begin{equation*}
		(1+\gamma_2)\sum_{l=i}^{K} \Delta_l \ge \sum_{l=i}^{K+1}\Delta_l
		+ \underbrace{\gamma_2\sum_{l=i}^{K-1}\Delta_l}_{\ge 0} \ge
		\Dmax,
	\end{equation*}
	thus $\sum_{l=i}^{K}\Delta_l \ge \frac{\Dmax}{1+\gamma_2}$.
	From \eqref{eq:C_est} we get $C(x_l)\ge \varepsilon$. 
	Due to \eqref{eq:neg_i_1} we can apply \eqref{eq:ATRM_pred_est_simple} and get
	\begin{equation*}
		\sum_{l=i}^{K}\pred(x_{j},\Delta_{j}) \ge
		\frac{f(\varepsilon)}{2}\sum_{l=i}^{K}\Delta_{l}\ge
		\frac{f(\varepsilon)}{2}\frac{\Dmax}{1+\gamma_2}\eqqcolon
		\delta_2>0
		.
	\end{equation*}
	\textbf{Case $\neg$(ii)$\wedge \neg$(i):}
	From $\neg$(ii), we get $\Delta_{K+1} < \Dmax$.
	In case that, additionally, $\Delta_{K+1} < \Dmax_{a}(x_{K+1})$,
	\eqref{eq:ATRM_pred_est_simple} implies acceptance of iteration $K + 1$.
	Otherwise, $\Delta_{K+1} \ge \Dmax_{a}(x_{K+1})$
	and \eqref{eq:accuracy_E} yields acceptance as well.
	This implies
	\begin{equation*}
		(iii) \Longrightarrow (i) \lor (ii),
	\end{equation*}
	i.e.,
	the distinction by cases is complete. In all cases, we have shown that
	\begin{equation*}
		\sum_{l \in I_1} \pred(x_l, \Delta_l)
		\ge
		\min\{\delta,\delta_1,\delta_2\}>0
		.
	\end{equation*}

	We proceed inductively and define $I_m$, $m > 1$,
	by repeating the above construction with $i \coloneqq n_k$,
	where $k$ is the smallest index such that $n_k \not\in I_{m-1}$.
	This allows for the final estimate
	\begin{align*}
		J(x_1) - \lim_{n\rightarrow\infty} J(x_n) &=
		\sum_{n=1}^{\infty}J(x_{n})-J(x_{n+1})\\
		&\ge \sum_{m=1}^{\infty}\sum_{l\in I_m} J(x_l) - J(x_{l+1}) \\
		&= \sum_{m=1}^{\infty}\sum_{l\in I_m}\ared(x_l,\Delta_l)\\
		&\ge \sum_{m=1}^{\infty}\sum_{l\in I_m}\eta\,\pred(x_l,\Delta_l)\\
		&\ge \sum_{m=1}^{\infty}\eta \min\{\delta,\delta_1,\delta_2\} = \infty.
	\end{align*}
	For the first inequality, we use that the trust-region method only accepts steps
	that decrease the objective and in the third line we get the equality
	exactly from the restriction to successful iterates. In conclusion, this violates the assumption of $J$ being bounded from below leading
	to a contradiction.
	Thus, no such subsequence $x_{n_k}$ with
	$C(x_{n_k})>2\varepsilon$ can exist, finishing the proof.
\end{proof}
\begin{remark}
	\label{rem:ATRM_conv_proof}
	The proof of \cref{thm:atrm_convergence} is constructive. The main ingredient is that
	for each criticality measure there exists a corresponding trust-region radius
	that guarantees sufficient improvement of the objective. This is described in
	detail for the first major case $\Delta_i > \Dmax$. It may however
	happen, that after an accepted step the next criticality measure is large,
	while the trust region radius is too small in relation. The construction of the second part
	of the proof uses the assumed Lipschitz continuity to ensure that the
	criticality measure remains large enough, while the algorithm expands the
	trust region. In this process it is proven that an unsuccessful iterate can
	only occur (condition $(iii)$) after a successful step with sufficient decrease
	was performed.
\end{remark}
\begin{remark}
	\label{rem:ATRM_conv_with_Delta_bound}
	As an interesting case, we discuss the situation in which
	$\Dmax_a(x) \equiv \Dmax_a > 0$, i.e., it does not depend on the point $x \in X$.
	Consequently, one could choose $\Dmax_b \in (0, \Dmax_a)$ arbitrarily
	and \eqref{eq:accuracy_E} would be automatically satisfied.
	If, additionally, there exists a bound
	$\Dmax_{\text{acc}}>0$ such that
	\begin{equation*}
		\ared(x, \Delta) \ge \eta \pred(x, \Delta)
		\qquad\forall x \in X, \Delta \in (0, \Dmax_{\text{acc}}),
	\end{equation*}
	then all trust-region radii are bounded from below by $\gamma_1\Dmax_{\text{acc}}$.
	Now, it can be checked that
	\cref{asmp:ATRM}\eqref{eq:stability} and \eqref{eq:stability_E}
	are no longer needed in the proof of \cref{thm:atrm_convergence}.
	Thus, we can define
	\begin{equation}
		\label{eq:Dmax2}
		\Dmax \coloneqq
		\min\left\{\Dmax_{\pred}(\varepsilon),\frac{\Dmax_{\text{acc}}}{\gamma_2},\Dmax_a\right\}
	\end{equation}
	and \eqref{eq:ATRM_pred_est_simple} holds for all $\Delta \in (0, \Dmax]$.
	By monotonicity, we also get $\pred(x, \Delta) \ge \frac{f(\varepsilon)}{2} \Dmax$
	for all $\Delta > \Dmax$.
	This leads to the estimate
	\begin{equation*}
		\ared(x_{n_k}, \Delta_{n_k})
		\ge
		\eta
		\pred(x_{n_k}, \Delta_{n_k})
		\ge
		\eta\frac{f(\varepsilon)}{2} \Dmax
		.
	\end{equation*}
	This is enough to argue as in the final estimate in the proof of \cref{thm:atrm_convergence}.
	In particular, all the distinctions by cases can be skipped.

	Note that this situation appears precisely
	in \cref{sec:convex_problem} for the proximal gradient algorithm,
	see \cref{lem:OC_conv_predictability} and \cref{rem:accuracy}.
\end{remark}
%%fakesection: Bib
\bibliographystyle{jnsao}
\bibliography{iioc_prox_TV}

\end{document}